\documentclass[11pt]{article}

\usepackage{graphicx}
\usepackage{amsthm,amsmath,amscd,amssymb}
\usepackage[all]{xy}

\newcommand{\N}{\mathbb{N}}
\newcommand{\Z}{\mathbb{Z}}

\title{Combinatorics of branchings in higher dimensional automata}
\author{Philippe Gaucher\\Institut de Recherche Math\'ematique Avanc\'ee\\
ULP et CNRS\\ 7 rue Ren\'e Descartes\\67084 Strasbourg Cedex\\France\\gaucher@irma.u-strasbg.fr}

\date{June 2001}

\newtheorem{thm}{Theorem}[section] 
\newtheorem{prop}[thm]{Proposition}
\newtheorem{lem}[thm]{Lemma}

\newtheorem{cor}[thm]{Corollary}
\newtheorem{rem}[thm]{Remark}

\newtheorem{defn}{Definition}[section]

\newcommand{\be}{\begin{equation}}
\newcommand{\ee}{\end{equation}}
\newcommand{\bea}{\begin{eqnarray}}
\newcommand{\eea}{\end{eqnarray}}
\newcommand{\beas}{\begin{eqnarray*}}
\newcommand{\eeas}{\end{eqnarray*}}

\newcommand{\F}{\mathcal{F}}

\newcommand{\C}{\mathcal{C}}
\newcommand{\D}{\mathcal{D}}

\newcommand{\B}{\mathcal{B}}

\newcommand{\bd}{\begin{defn}}
\newcommand{\ed}{\end{defn}}
\newcommand{\bcd}{\begin{defn}}
\newcommand{\ecd}{\end{defn}}
\newcommand{\blem}{\begin{lem}}
\newcommand{\elem}{\end{lem}}
\newcommand{\bp}{\begin{prop}}
\newcommand{\ep}{\end{prop}}
\newcommand{\bcor}{\begin{cor}}
\newcommand{\ecor}{\end{cor}}
\newcommand{\bth}{\begin{thm}}
\renewcommand{\eth}{\end{thm}}
\newcommand{\bi}{\begin{enumerate}}
\newcommand{\ei}{\end{enumerate}}
\newcommand{\br}{\begin{rem}}
\newcommand{\er}{\end{rem}}
\newcommand{\bpf}{\begin{proof}}
\newcommand{\epf}{\end{proof}}

\newcommand{\p}\times

\newcommand{\iso}{\cong}

\newcommand{\de}{\partial}

\renewcommand{\P}{\mathbb{P}}

\newcommand{\fl}[1]{\ar@{->}[l]_{#1}}
\newcommand{\fr}[1]{\ar@{->}[r]^{#1}}
\newcommand{\fd}[1]{\ar@{->}[d]_{#1}}
\newcommand{\fu}[1]{\ar@{->}[u]|{#1}}
\newcommand{\f}[2]{\ar@{->}[#1]|{#2}}
\newcommand{\ff}[2]{\ar@2{->}[#1]|{#2}}

\newcommand{\coin}[2]{\begin{picture}(11,11)(0,0)
\put(2,0){\vector(1,0){10}}\put(13,0){#1}
\put(2,0){\vector(0,1){10}}\put(2,13){#2}
\end{picture}}

\newcommand{\degv}{\begin{picture}(11,11)(0,0)\linethickness{1pt}
\put(2,-1){\line(0,1){8}}\put(9,-1){\line(0,1){8}}
\end{picture}}

\newcommand{\degh}{\begin{picture}(11,11)(0,0)\linethickness{1pt}
\put(2,7){\line(1,0){7}}\put(2,-1){\line(1,0){7}}
\end{picture}}

\newcommand{\gamp}{\begin{picture}(11,11)(0,0)\linethickness{1pt}
\put(2,-1){\line(0,1){8}}\put(2,-1){\line(1,0){7}}
\end{picture}}

\newcommand{\gamn}{\begin{picture}(11,11)(0,0)\linethickness{1pt}
\put(9,-1){\line(0,1){8}}\put(2,7){\line(1,0){7}}
\end{picture}}

\newcommand{\car}{\begin{picture}(11,11)(0,0)\linethickness{1pt}
\put(2,-1){\line(0,1){8}}\put(2,-1){\line(1,0){7}}
\put(9,-1){\line(0,1){8}}\put(2,7){\line(1,0){7}}
\end{picture}}

\newcommand{\HR}{H\!R}
\newcommand{\CR}{C\!R}
\newcommand{\HF}{H\!F}
\newcommand{\CF}{C\!F}

\newcommand{\tr}{tr}
\newcommand{\AB}{{\!\!}<{\!\!}A,B{\!\!}>{\!\!}}

\begin{document}

\maketitle

\begin{abstract}
  We explore the combinatorial properties of the branching areas of
  execution paths in higher dimensional automata. Mathematically, this
  means that we investigate the combinatorics of the negative corner
  (or branching) homology of a globular $\omega$-category and the
  combinatorics of a new homology theory called the reduced branching
  homology. The latter is the homology of the quotient of the
  branching complex by the sub-complex generated by its thin elements.
  Conjecturally it coincides with the non reduced theory for higher
  dimensional automata, that is $\omega$-categories freely generated
  by precubical sets. As application, we calculate the branching
  homology of some $\omega$-categories and we give some invariance
  results for the reduced branching homology. We only treat the
  branching side. The merging side, that is the case of merging areas
  of execution paths is similar and can be easily deduced from the
  branching side.  \footnote{2000 AMS classification : 55U}
\end{abstract}

\section{Introduction}

After \cite{Pratt,HDA}, one knows that it is possible to model higher
dimensional automata (HDA) using precubical sets
(Definition~\ref{def_cubique}). In such a model, a $n$-cube
corresponds to a $n$-transition, that is the concurrent execution of
$n$ $1$-transitions. This theoretical idea would be  implemented later.
Indeed a CaML program translating programs in Concurrent Pascal into a
text file coding a precubical set is presented in
\cite{cridlig96implementing}. At this step, one does not yet consider
cubical sets with or without connections since the degenerate elements
have no meaning at all from the point of view of computer-scientific
modeling (even if in the beginning of \cite{Gau}, the notion of cubical sets
is directly introduced by intellectual reflex).

In \cite{HDA}, the following fundamental observation is made :
\textit{given a precubical set $(K_n)_{n\geq 0}$ together with its
  two families of face maps $(\de_i^\alpha)$ for $\alpha\in\{-,+\}$,
  then both chain complexes $(\Z K_*,\de^\alpha)$, where $\Z X$ means
  the free abelian group generated by X and where
  $\de^\alpha=\sum_{i}(-1)^{i+1}\de_i^\alpha$, give rise to two
  homology theories $H^\alpha_*$ for $\alpha\in\{-,+\}$ whose
  non-trivial elements model the branching areas of execution paths
  for $\alpha=-$ and the merging areas of execution paths for
  $\alpha=+$ in strictly positive dimension. Moreover the group $H_0^-$ (resp.
  $H_0^+$) is the free abelian group generated by the final states
  (resp.  the initial states) of the HDA.}

\begin{figure}
\begin{center}
\includegraphics[width=7cm]{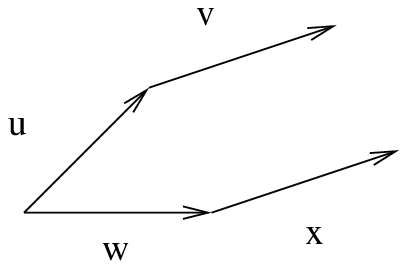}
\end{center}
\caption{A $1$-dimensional branching area}
\label{ex1}
\end{figure}

Consider for instance the $1$-dimensional HDA of Figure~\ref{ex1}.
Then $u-w$ gives rise to a non-trivial homology class which
corresponds to the branching which is depicted.

\begin{figure}
\begin{center}
\includegraphics[width=3cm]{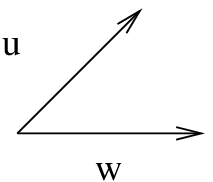}
\end{center}
\caption{A $1$-dimensional branching area}
\label{ex11}
\end{figure}

Then the first problem is that the category of precubical sets is
not appropriate to identify the HDA of Figure~\ref{ex1} with that of
Figure~\ref{ex11} because there is no morphism between them preserving 
the initial state and both  final states.

No matter : it suffices indeed to work with the category of
precubical sets endowed with the $+_i$ cubical composition laws 
satisfying the axioms of Definition~\ref{cubicalomega} and with the
morphisms obviously defined. Now for any $n\geq 1$, there are $n$
cubical composition laws $+_1,\dots,+_n$ representing the
concatenation of $n$-cubes in the $n$ possible directions. Let 
$X=u+_1 v$ and $Y:= w+_1 x$. Then there is a unique morphism $f$ in this new
category of HDA from the HDA of Figure~\ref{ex11} to the HDA of
Figure~\ref{ex1} such that $f:u\mapsto X$ and $f:w\mapsto Y$. However
$f$ is not invertible in the category of precubical sets equipped
with cubical composition laws because there still does not exist any morphism 
from the HDA of Figure~\ref{ex1} to the HDA of
Figure~\ref{ex11}.

To make $f$ invertible (recall that we would like to find a category
where both HDA would be isomorphic), it still remains the possibility
of formally adding inverses by the process of localization of a
category with respect to a collection of morphisms.  However a serious
problem shows up : the non-trivial cycles $u-w$ and $X-Y$ of
Figure~\ref{ex1} give rise to two distinct homology classes although
these two distinct homology classes correspond to the same branching
area. Indeed there is no chain in dimension $2$ (i.e.
$K_2=\emptyset$), so no way to make the required identification !

\begin{figure}
\begin{center}
\includegraphics[width=14cm]{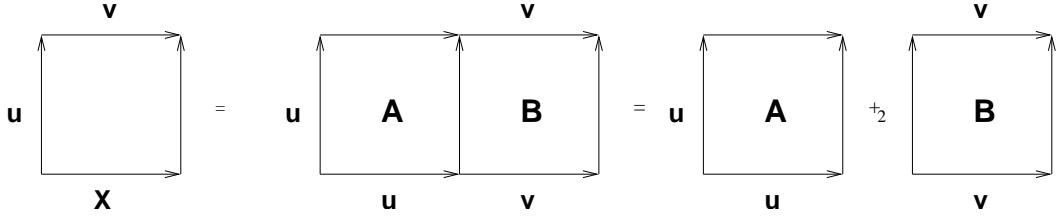}
\end{center}
\caption{Identifying $u+_1v$ and $u$}
\label{identificationbythin}
\end{figure}

This means that something must be added in dimension $2$, but without 
creating additional homology classes.  Now
consider Figure~\ref{identificationbythin}. The element $A$ must be
understood as a \textit{thin} $2$-cube such that, with our convention
of orientation, $\de_1^-A=u$, $\de_2^-A=u$,
$\de_1^+A=\epsilon_1\de_1^-v$, $\de_2^+A=\de_2^-B=\epsilon_1\de_1^-v$.
And the element $B$ must be understood as another \textit{thin}
$2$-cube such that $\de_2^-B=\epsilon_1\de_1^-v$,
$\de_2^+B=\epsilon_1\de_1^+v$ and $\de_1^-B=\de_1^+B=v$. In such a
situation, $\de^-(A +_2 B)=u+_1v -u$ therefore $u+_1v$ and $u$ become
equal in the first homology group $H_1^-$. By adding this kind of thin
$2$-cubes to the chain complex $(\Z K_*,\de^-)$, one can then identify the
two cycles $u-w$ and $X-Y$. One sees that there are two kinds of thin
cubes which are necessary to treat the branching case. The first kind
is well-known in cubical set theory : this is for example
$B=\epsilon_1 v$ or $\de_1^+A=\epsilon_1\de_1^-v$.  The second kind is
for example $A$ which will be denoted by $\Gamma_1^- u$ and which
corresponds to extra-degeneracy maps as defined in \cite{Brown_cube}.

To take into account the symmetric problem of merging areas of
execution paths, a third family $\Gamma_i^+$ of degeneracy maps will be necessary. In this 
paper, we will only treat the case of branchings. The case of mergings is similar and 
easy to deduce from the branching case.
The solution presented in this paper to overcome the above problems is then 
as follows : 
\begin{itemize}
\item One considers the free globular $\omega$-category $F(K)$
  generated by the precubical set $K$ : it is obtained by
  associating to any $n$-cube $x$ of $K$ a copy of the free globular
  $\omega$-category  $I^n$ generated by the faces of the $n$-cube 
(paragraph~\ref{In}) ; the faces of this $n$-cube are denoted 
by $(x;k_1\dots k_n)$ ; one takes the direct sum of all these cubes and one 
takes the quotient by the relations 
\[(\de_i^\alpha y;k_1\dots k_n)\sim (y;k_1\dots k_{i-1}\alpha k_i\dots k_n)\]
for any $y\in K_{n+1}$, $\alpha\in\{-,+\}$ and $1\leq i\leq n+1$.
\item Then we take its cubical singular nerve
  $\mathcal{N}^\square(F(K))$ (which is equal also to the free cubical
  $\omega$-category generated by $K$) ; the required thin elements
  above described (the three families $\epsilon_i$, $\Gamma_i^-$ and
  $\Gamma_i^+$) do appear in it as components of the algebraic
  structure of the cubical nerve (Definition~\ref{cubicalomega} and
  Definition~\ref{cubical_singular}).
\item The branching homology of $F(K)$ (Definition~\ref{def_orientee}) 
is the solution for both following reasons : 
\begin{enumerate}
\item Let $x$ and $y$ be two $n$-cubes of the cubical nerve which are in the
branching complex. If $x+_j y$ exists for some $j$ with $1\leq j\leq n$, then 
$x$ and $x+_j y$ are equal modulo elements in the chain complex generated by 
the thin elements (Theorem~\ref{folding_plus}) ;
\item The chain complex generated by the thin elements is conjecturally acyclic 
in this situation, and so it does not create non-trivial homology classes 
(Conjecture~\ref{thin}).
\end{enumerate}
\end{itemize}

We have explained above the situation in dimension $1$. The
$2$-dimensional case is depicted in Figure~\ref{ex3}. Additional
explanations are available at the end of
Section~\ref{folding_composition}.

The branching homology (or negative corner homology) and the merging
homology (or positive corner homology) were already introduced in
\cite{Gau}.  This invariance with respect to the cubifications of the
underlying HDA was already suspected for other reasons. The branching
and merging homology theories are the solution to overcome the
drawback of Goubault's constructions.

There are three key concepts in this paper which are not so
common in the general literature and which we would like to
draw to the reader's attention.

\begin{enumerate}
\item the extra structure of {\em connections} $\Gamma^\pm$ on
cubical sets, which allow extra degenerate elements in which
adjacent faces coincide. This structure was first introduced
in \cite{Brown_cube}.
\item   the notion of {\em folding operator}. This was introduced in
the groupoid context in \cite{Brown_cube}, to fold down a cube to an
element in a  crossed complex, and in the category context
in \cite{phd-Al-Agl} to fold down a cube to an element in a globular
category. Properties of this folding operator are further
developed in \cite{equiv_glob_cub}. This we call the `usual folding operator'.
\item  the notion of {\em thin cube}, namely a multiple composition of
cubes of the form $\epsilon_i y$ or $\Gamma^\pm z$. A crucial
result is that these are exactly the elements which fold
down to 1 in the contained globular category.
\end{enumerate}

So there are many ways of choosing a cycle in the branching complex
for a given homology class, i.e.  a given branching area, according to
the choice of the cubification of the considered HDA. This possibility
of choice reveals an intricate combinatorics. The most appropriate
tool constructed in the mathematical theory of cubical sets to study
this combinatorics is not relevant here. The machinery of folding
operators \cite{Brown_cube,phd-Al-Agl} does not work indeed for the
study of the branching homology because the usual folding operators
are not internal to the branching chain complex (see
Section~\ref{folding_op}).  The core of this paper is the proposal of
a new folding operator adapted for the study of the branching complex
(Section~\ref{folding_op2}). This operator enables us to deduce
several results on the reduced branching homology, the latter being
obtained by taking the quotient of the former by the sub-complex
generated by its thin elements. This sub-complex is conjecturally
acyclic for a wide variety of $\omega$-categories, including that
freely generated by a precubical set or a globular set
(Conjecture~\ref{thin}). Our main result is that the negative folding
operator induces the identity map on the reduced branching complex
(Corollary~\ref{identification}). Using some relations between the
branching homology of some particular $\omega$-categories and the
usual simplicial homology of some associated $\omega$-categories
(Theorem~\ref{calcul}), the behaviour of the composition maps (the
globular and the cubical ones) modulo thin elements is completely
studied (Section~\ref{folding_composition}). All these results lead us
to a question about the description of the reduced branching complex
using globular operations by generators and relations
(Proposition~\ref{coin_formel} and Question~\ref{formelle}) and to two
invariance results for the reduced branching homology
(Proposition~\ref{invariance} and Theorem~\ref{invariance2}).

This paper is organized as follows. Section~\ref{prel} recalls some
important notations and conventions for the sequel. In
Section~\ref{corner}, the branching homology and the reduced
branching homology are introduced. In Section~\ref{matrix}, the
matrix notations for connections and degeneracies are described. Next
in Section~\ref{relation_usual_nerve}, the branching homology of some
particular $\omega$-categories (the $\omega$-categories of length at
most $1$) is completely calculated in terms of the usual simplicial
homology.  In Section~\ref{construct_folding}, the negative folding
operators are introduced. In Section~\ref{elem_move}, the negative
folding operators are decomposed in terms of elementary moves. In
Section~\ref{id_phi}, we prove that each elementary move appearing in
the decomposition of the folding operators induces the identity map on
the reduced branching complex. Therefore the folding operators
induce the identity map as well. In Section~\ref{folding_composition},
the behaviour of the cubical and globular composition laws in the
reduced branching complex is completely studied. In the following
Section~\ref{diff_section}, some facts about the differential map in
the reduced branching complex are exposed. In the last
Section~\ref{invariance_result}, some invariance results for the
reduced branching homology are exposed and the reduced branching homology is
calculated for some simple globular $\omega$-categories.

\section{Preliminaries : cubical set, globular and cubical category}\label{prel}

Here is a recall of some basic definitions, in order to make precise
some notations and some conventions for the sequel.

\bd\label{def_cubique}\cite{Brown_cube} \cite{cube}
A \textit{cubical set} consists of a family of
sets $(K_n)_{n\geqslant 0}$, of a family of face maps
$\xymatrix@1{{K_n}\fr{\de_i^\alpha} &{K_{n-1}}}$ for $\alpha\in\{-,+\}$
and of a family of degeneracy maps
$\xymatrix@1{K_{n-1}\fr{\epsilon_i}&K_{n}}$ with $1\leqslant i\leqslant n$ 
which satisfy the following relations
\begin{enumerate}
\item $\de_i^\alpha \de_j^\beta = \de_{j-1}^\beta \de_i^\alpha$
for all $i<j\leqslant n$ and $\alpha,\beta\in\{-,+\}$ (called sometimes 
the cube axiom)
\item $\epsilon_i\epsilon_j=\epsilon_{j+1}\epsilon_i$
for all $i\leqslant j\leqslant n$
\item $\de_i^\alpha \epsilon_j=\epsilon_{j-1}\de_i^\alpha$       for
  $i<j\leqslant n$
and $\alpha\in\{-,+\}$
\item $\de_i^\alpha \epsilon_j=\epsilon_{j}\de_{i-1}^\alpha$   for
  $i>j\leqslant n$ and
$\alpha\in\{-,+\}$
\item $\de_i^\alpha \epsilon_i=Id$
\end{enumerate}
A family $(K_n)_{n\geq 0}$ only equipped with a family of face maps
$\de_i^\alpha$ satisfying the same axiom as above is called a 
\textit{precubical set}. An element of $K_0$ will be sometimes called 
a \textit{state}, or a $0$-cube and an element of $K_n$ a $n$-cube, 
or a $n$-dimensional cube.
\ed

\bd Let $(K_n)_{n\geq 0}$ and $(L_n)_{n\geq 0}$ be two cubical
(resp. precubical) sets. Then a morphism $f$ from 
$(K_n)_{n\geq 0}$ to $(L_n)_{n\geq 0}$ is a family $f=(f_n)_{n\geq 0}$ 
of set maps $f_n:K_n\rightarrow L_n$ such that $f_n\de_i^\alpha =
\de_i^\alpha f_n$ and $f_n\epsilon_i = \epsilon_i f_n$
(resp. $f_n\de_i^\alpha = \de_i^\alpha f_n$) for any $i$. The
corresponding category of cubical sets is isomorphic to the
category of pre-sheaves $Sets^{\square^{op}}$ over a small
category $\square$. The corresponding category of precubical
sets is isomorphic to the category of pre-sheaves
$Sets^{{\square^{pre}}^{op}}$ over a small category
$\square^{pre}$.  \ed

\bd\label{omega_categories}\cite{Brown-Higgins0} \cite{oriental} \cite{Tensor_product}
A (globular) \textit{$\omega$-category} is a set $A$
endowed with two families of maps $(d_n^-=s_n)_{n\geqslant 0}$ and
$(d_n^+=t_n)_{n\geqslant 0}$ from $A$ to $A$ and with a family of partially
defined 2-ary operations $(*_n)_{n\geqslant 0}$ where for any
$n\geqslant 0$, $*_n$ is a map from $\{(a,b)\in A\p A,
t_n(a)=s_n(b)\}$ to $A$ ($(a,b)$ being carried over $a *_n b$) which
satisfy the following axioms for all $\alpha$ and $\beta$ in
$\{-,+\}$ :

\begin{enumerate}
\item $d_m^\beta d_n^\alpha x=
\left\{\begin{CD}d_m^\beta x \hbox{  if $m<n$}\\  d_n^\alpha x \hbox{  if $m\geqslant n$}
\end{CD}\right.$
\item $s_n x *_n x= x *_n t_n x = x$
\item if  $x *_n y$ is well-defined, then  $s_n(x *_n y)=s_n x$, $t_n(x *_n y)=t_n y$
and for  $m\neq n$, $d_m^\alpha(x *_n y)=d_m^\alpha x *_n d_m^\alpha y$
\item as soon as the two members of the following equality exist, then
$(x *_n y) *_n z= x *_n (y *_n z)$
\item if $m\neq n$ and if the two members of the equality make sense, then
$(x *_n y)*_m (z*_n w)=(x *_m z) *_n (y *_m w)$
\item for any  $x$ in $A$, there exists a natural number $n$ such that $s_n x=t_n x=x$
(the smallest of these numbers is called the dimension of $x$ and is denoted
by $dim(x)$).
\end{enumerate}
\ed

A globular set is a set $A$ endowed with two families of maps
$(s_n)_{n\geqslant 0}$ and $(t_n)_{n\geqslant 0}$ satisfying the same
axioms as above
\cite{petittoposglob,Penon_weakcat,Batanin_1}. We call $s_n(x)$
the $n$-source of $x$ and $t_n(x)$ the $n$-target of
$x$.

\subsection*{Notation} The category of all $\omega$-categories (with the obvious
morphisms) is denoted by $\omega Cat$. The corresponding morphisms are
called $\omega$-functors. The set of $n$-dimensional morphisms of $\C$
is denoted by $\C_n$. The set of morphisms of $\C$ of dimension lower
or equal than $n$ is denoted by $\tr_n\C$. The element of $\C_0$ will be 
sometimes called \textit{states}. An \textit{initial state} 
(resp. \textit{final state}) of $\C$ is a $0$-morphism 
$\alpha$ such that $\alpha =s_0x$ (resp. $\alpha =t_0x$) implies $x=\alpha$.

\bd\cite{Brown_cube,phd-Al-Agl}\label{cubicalomega}
A \textit{cubical $\omega$-category} consists of a cubical set
$$((K_n)_{n\geqslant 0},\de_i^\alpha,\epsilon_i)$$
together with two
additional families of degeneracy maps called connections
$$\Gamma_i^\alpha:\xymatrix{K_{n}\fr{} & K_{n+1}} $$ with
$\alpha\in\{-,+\}$, $n\geqslant 1$ and $1\leqslant i \leqslant n$
and a family
of associative operations $+_j$ defined on $\{(x,y)\in K_n\p K_n,
\de^+_j x=\de_j^- y\}$ for $1\leqslant j\leqslant n$ such that
\begin{enumerate}
\item $\de_i^\alpha\Gamma_j^\beta=\Gamma^\beta_{j-1}\de_i^\alpha$
  for all $i<j$ and all $\alpha,\beta\in\{-,+\}$
\item $\de_i^\alpha\Gamma_j^\beta=\Gamma^\beta_{j}\de_{i-1}^\alpha$
  for all $i>j+1$ and all $\alpha,\beta\in\{-,+\}$
\item $\de_j^\pm \Gamma_j^\pm = \de_{j+1}^\pm \Gamma_j^\pm =Id$
\item $\de_j^\pm \Gamma_j^\mp = \de_{j+1}^\pm \Gamma_j^\mp= \epsilon_{j}\de_{j}^\pm $
\item $\Gamma_i^\pm \Gamma_j^\pm =\Gamma_{j+1}^\pm
\Gamma_i^\pm $ if $i\leqslant j$\
\item $\Gamma_i^\pm \Gamma_j^\mp= \Gamma_{j+1}^\mp\Gamma_i^\pm $  if $i<j$
\item $\Gamma_i^\pm \Gamma_j^\mp= \Gamma_{j}^\mp\Gamma_{i-1}^\pm $  if $i>j+1$
\item $\Gamma_i^\pm \epsilon_j=\epsilon_{j+1}\Gamma_i^\pm $ if $i<j$
\item $\Gamma_i^\pm \epsilon_j=\epsilon_i\epsilon_i$  if $i=j$
\item $\Gamma_i^\pm \epsilon_j=\epsilon_j\Gamma_{i-1}^\pm $ if $i>j$
\item $(x+_j y)+_j z=x+_j(y+_j z)$
\item $\de^-_j(x +_j y)=\de^-_j(x)$
\item $\de^+_j(x +_j y)=\de^+_j(y)$
\item $\de^\alpha_i(x +_j y)=\left\{
\begin{CD}
&&\de^\alpha_i(x) +_{j-1} \de^\alpha_i(y)\hbox{ if $i<j$}\\
&&\de^\alpha_i(x) +_{j} \de^\alpha_i(y)\hbox{ if $i>j$}\\ \end{CD}
\right.$
\item $(x +_i y) +_j (z +_i t)=(x +_j z) +_i (y +_j t)$.
\item $\epsilon_i(x +_j y)=\left\{
\begin{CD}&&\epsilon_i(x) +_{j+1} \epsilon_i(y) \hbox{ if $i\leqslant j$}\\
&&\epsilon_i(x) +_{j} \epsilon_i(y) \hbox{ if $i>j$}\\
\end{CD}\right.$
\item $\Gamma_i^\pm (x +_j y)=\left\{
\begin{CD}&&\Gamma_i^\pm (x) +_{j+1} \Gamma_i^\pm (y)\hbox{ if $i<j$}\\
&&\Gamma_i^\pm (x) +_{j} \Gamma_i^\pm (y)\hbox{ if $i> j$}
\end{CD}\right.$
\item If $i=j$, $\Gamma_i^- (x +_j y)=\left[\begin{array}{cc}
\epsilon_{j+1}(y) & \Gamma_j^- (y)\\
\Gamma_j^- (x) & \epsilon_j(y)\\
\end{array}\right] \coin{j+1}{j}$
\item If $i=j$, $\Gamma_i^+ (x +_j y)=\left[\begin{array}{cc}
\epsilon_{j}(x) & \Gamma_j^+ (y)\\
\Gamma_j^+ (x) & \epsilon_{j+1}(x)\\
\end{array}\right] \coin{j+1}{j}$
\item $\Gamma_j^+ x +_{j+1} \Gamma_j^- x = \epsilon_j x$ and
$\Gamma_j^+ x +_{j} \Gamma_j^- x = \epsilon_{j+1} x$
\item $\epsilon_i \de_i^- x +_i x = x +_i \epsilon_i \de_i^+ x = x$
\end{enumerate}
The corresponding category with the obvious morphisms is denoted by
$\infty Cat$.
\ed

Without further precisions, the word $\omega$-category is always
supposed to be taken in the sense of globular $\omega$-category.  In
\cite{equiv_glob_cub}, it is proved that the category of cubical
$\omega$-categories and the category of globular $\omega$-categories
are equivalent.

\subsection*{Notation} 
If $S$ is a set, the free abelian group generated by $S$ is denoted by
$\Z S$. By definition, an element of $\Z S$ is a formal linear
combination of elements of $S$.

\bd\cite{Gau} Let $\C$ be an $\omega$-category. Let $\C_n$ be the set of
$n$-dimensional morphisms of $\C$. Two $n$-morphisms $x$ and $y$ are
\textit{homotopic} if there exists $z\in \Z\C_{n+1}$ such that $s_n
z-t_n z= x-y$. This property is denoted by $x\sim y$.  \ed

We have already observed in \cite{Gau} that the corner homologies do
not induce functors from $\omega Cat$ to the category of abelian
groups. A notion of non-contracting $\omega$-functors was required.

\bd\label{noncontractant}\cite{Gau} Let $f$ be an $\omega$-functor from $\C$ to
$\D$.  The morphism $f$ is \textit{non-contracting} if for any
$1$-dimensional $x\in \C$, the morphism $f(x)$ is a $1$-dimensional
morphism of $\D$. \ed

The theoretical developments of this paper and future works in progress entail
the following definitions too.

\bd Let $\C$ be an $\omega$-category. Then $\C$ is \textit{non-contracting}
if and only if for any $x\in \C$ of strictly positive dimension, $s_1x$ and
$t_1x$ are $1$-dimensional (they could be a priori $0$-dimensional as well).
\ed

A justification of this definition among a lot of them is that if
$\C$ is an $\omega$-category which is not non-contracting, then there
exists a morphism $u$ of $\C$ such that $dim(u)>1$ and such that for
instance $s_1u$ is $0$-dimensional.  For example consider the
two-element set $\{A,\alpha\}$ with the rules
$s_1A=t_1A=s_0A=t_0A=\alpha$ and $s_2A=t_2A=A$. This defines an 
$\omega$-category which
is not non-contracting. Then $A$ is
$2$-dimensional though $s_1A$ and $t_1A$ are $0$-dimensional. And in
this situation $\square_2^-(A)$ defined in Section~\ref{folding_op2}
is not an element of the branching nerve, and therefore for that
$\C$, the morphism $\CF^-_2(\C)$ (see Proposition~\ref{coin_formel})
to $\CR^-_2(\C)$ is not defined.

\subsection*{Notation} 
The category of non-contracting  $\omega$-categories with the
non-contracting $\omega$-functors is denoted by $\omega Cat_1$.

If $f$ is a non-contracting $\omega$-functor from $\C$ to $\D$, then
for any morphism $x\in \C$ of dimension greater than $1$, $f(x)$ is of
dimension greater than one as well. This is due to the equality
$f(s_1 x)=s_1 f(x)$.

All globular $\omega$-categories that will appear in this work will be non-contracting.

\section{Reduced branching homology}\label{corner}

\subsection{The globular $\omega$-category $I^n$}\label{In} 
We need first to describe precisely the $\omega$-category associated
to the $n$-cube. 
Set $\underline{n}=\{1,...,n\}$ and let $\underline{cub}^n$ be the set of maps
from $\underline{n}$ to $\{-,0,+\}$ (or in other terms the set of
words of length $n$ in the alphabet $\{-,0,+\}$). We say that an
element $x$ of $\underline{cub}^n$ is of dimension $p$ if $x^{-1}(0)$ is a set
of $p$ elements. The set $\underline{cub}^n$ is supposed to be graded by the
dimension of its elements. The set $\underline{cub}^0$ is the set of maps from
the empty set to $\{-,0,+\}$ and therefore it is a singleton. Let
$y\in \underline{cub}^i$. Let $r_y$ be the map from $(\underline{cub}^n)_i$ to
$(\underline{cub}^n)_{dim(y)}$ defined as follows (with $x\in (\underline{cub}^n)_i$)
: for $k\in\underline{n}$, $x(k)\neq 0$ implies $r_y(x)(k)=x(k)$ and
if $x(k)$ is the $l$-th zero of the sequence $x(1),...,x(n)$, then
$r_y(x)(k)=y(\ell)$. If for any $\ell$ between $1$ and $i$,
$y(\ell)\neq 0$ implies $y(\ell)=(-)^\ell$, then we set
$b_y(x):=r_y(x)$. If for any $\ell$ between $1$ and $i$, $y(\ell)\neq
0$ implies $y(\ell)=(-)^{\ell+1}$, then we set $e_y(x):=r_y(x)$. We
have

If $x$ is an element of $\underline{cub}^n$, let us denote by $R(x)$ the
subset of $\underline{cub}^n$ consisting of $y\in \underline{cub}^n$ such that
$y$ can be obtained from $x$ by replacing some occurrences of $0$
in $x$ by $-$ or $+$. For example, $-00++-\in R(-000+-)$ but
$+000+-\notin R(-000+-)$. If $X$ is a subset of $\underline{cub}^n$, then
let $R(X)=\bigcup_{x\in X}R(x)$. Notice that $R(X\cup Y)=R(X)\cup
R(Y)$.

\bth There is one and only one $\omega$-category $I^n$ such that
\begin{enumerate}
\item the underlying set of $I^n$ is included in the set of
subsets of $\underline{cub}^n$
\item the underlying set of $I^n$ contains
all subsets like $R(x)$ where $x$ runs over $\underline{cub}^n$
\item all elements of $I^n$ are compositions of $R(x)$ where $x$ runs over $\underline{cub}^n$
\item for $x$ $p$-dimensional with $p\geqslant 1$, one has
\beas
&&s_{p-1}(R(x))=R\left(\{b_y(x),dim(y)=p-1\}\right)\\
&&t_{p-1}(R(x))=R\left(\{e_y(x),dim(y)=p-1\}\right)
\eeas
\item if $X$ and $Y$ are two elements of $I^n$ such that
$t_p(X)=s_p(Y)$ for some $p$, then $X\cup Y\in I^n$ and
$X\cup Y=X *_p Y$.
\end{enumerate}
Moreover, all elements $X$ of $I^n$ satisfy the equality $X=R(X)$.
\eth

The elements of $I^n$ correspond to the loop-free well-formed sub
pasting schemes of the pasting scheme $\underline{cub}^n$ \cite{CPS}
\cite{Crans_Tensor_product} or to the molecules of an $\omega$-complex
in the sense of \cite{MR99e:18008}. The condition ``\textit{$X *_n Y$ exists if
and only if $X\cap Y=t_n X=s_n Y$}'' of \cite{MR99e:18008} is not
necessary here because the situation of \cite{MR99e:18008} Figure 2
cannot appear in a composable pasting scheme.

The map which sends every $\omega$-category $\C$ to
$\mathcal{N}^\square(\C)_*=\omega Cat(I^*,\C)$ induces a functor from
$\omega Cat$ to the category of cubical sets. If $x$ is an
element of $\omega Cat(I^n,\C)$, $\epsilon_i(x)$ is the
$\omega$-functor from $I^{n+1}$ to $\C$ defined by
$\epsilon_i(x)(k_1...k_{n+1})=x(k_1...\widehat{k_i}...k_{n+1})$ for
all $i$ between $1$ and $n+1$ and $\de_i^\alpha(x)$ is the
$\omega$-functor from $I^{n-1}$ to $\C$ defined by
$\de_i^\alpha(x)(k_1...k_{n-1})=x(k_1...k_{i-1}\alpha k_i...k_{n-1})$
for all $i$ between $1$ and $n$.

The arrow $\de_i^\alpha$ for a given $i$ such that $1\leqslant
i\leqslant n$ induces a natural transformation from $\omega
Cat(I^n,-)$ to $\omega Cat(I^{n-1},-)$ and therefore, by Yoneda,
corresponds to an $\omega$-functor $\delta^\alpha_i$ from $I^{n-1}$ to
$I^n$.  This functor is defined on the faces of $I^{n-1}$ by
$\delta^\alpha_i(k_1...k_{n-1})=R(k_1...[\alpha]_i...k_{n-1})$. The
notation $[...]_i$ means that the term inside the brackets is at the
$i$-th place.

\bd\label{cubical_singular} The cubical set
$(\omega Cat(I^*,\C),\de_i^\alpha,\epsilon_i)$
is called the \textit{cubical singular nerve} of the $\omega$-category
$\C$. \ed

\subsection{Remark}\label{2diff} 
For $\alpha\in\{-,+\}$, and $x\in \omega Cat(I^n,\C)$, let
\[\de^\alpha x := \sum_{i=1}^{n}(-1)^{i+1} \de_i^\alpha x\]
Because of the cube axiom, one has $\de^\alpha\circ \de^\alpha=0$.

\bd\label{def_orientee}\cite{Gau} Let $\C$ be a non-contracting 
$\omega$-category. The set of $\omega$-functors $x\in \omega Cat(I^n,\C)$ 
such that for any
$1$-morphism $u$ with $s_0u=-_{n+1}$, $x(u)$ is $1$-dimensional (a priori
$x(u)$ could be $0$-dimensional as well) is denoted by $\omega Cat(I^n,\C)^-$.
Then $$\de^-(\Z\omega Cat(I^{*+1},\C)^-)\subset \Z\omega Cat(I^{*},\C)^-$$
by construction. We  set
$$H_*^-(\C)=H_*(\Z \omega Cat(I^{*},\C)^-,\de^-)$$ and we call
this homology theory the \textit{branching homology} of
$\C$. The cycles are called the \textit{branchings} of $\C$. 
The map $H_*^-$ induces a functor from $\omega
Cat_1$ to $Ab$.  \ed

The definition of $\omega Cat(I^n,\C)^-$ is a little bit different
from that of \cite{Gau}. Both definitions coincide if $\C$ is the free
$\omega$-category generated by a precubical set or a globular set. This
new definition ensures that the elementary moves introduced in
Section~\ref{elem_move} are well-defined on the branching nerve.
Otherwise it is easy to find counterexample, even in the case of a
non-contracting $\omega$-category.

\subsection{Conjecture}\label{thin} (About the thin elements of the branching
complex) \textit{Let $\C$ be a globular $\omega$-category which is either the
free globular $\omega$-category generated by a precubical set or the free
globular $\omega$-category generated by a globular set.  Let $x_i$ be
elements of $\omega Cat(I^n,\C)^- $ and let $\lambda_i$ be natural
numbers, where $i$ runs over some set $I$. Suppose that for any $i$,
$x_i(0_n)$ is of dimension strictly lower than $n$ (one calls it a
thin element). Then $\sum_i \lambda_i x_i$ is a boundary if and only
if it is a cycle.}

The thin elements conjecture is not true in general. Here is a
counterexample.  Consider an $\omega$-category $\C$ constructed by
considering $I^2$ and by dividing by the relations $R(-0)=R(0-)$ and
$R(-0)*_0 R(0+)=R(0-)*_0 R(+0)$.  Then the $\omega$-functor $F\in
\omega Cat(I^2,\C)^-$ induced by the identity functor from $I^2$ to
itself is a thin cycle in the branching homology. One can verify
that this cycle would be a boundary if and only if $R(0+)$ was
homotopic to $R(+0)$ in $\C$. This observation suggests the following
questions.

\bd Let $\C$ be an $\omega$-category. Then the $n$-th composition law
is said to be left regular up to homotopy if and only if for any
morphisms $x$, $y$ and $z$ such that $x*_n y=x *_n z$, then $y\sim
z$. \ed

\subsection{Question}
\textit{Does the thin elements conjecture hold for an $\omega$-category $\C$
such that all composition laws $*_n$ for any $n\geqslant 0$ are left
regular up to homotopy ?}

\subsection{Question}
\textit{How may we characterize the $\omega$-categories for which the thin
elements conjecture holds ?}

\bd Let  $M_n^-(\C)\subset  \Z\omega Cat(I^n,\C)^-$ be the sub-$\Z$-module generated
by the thin
elements (M for ``mince'' which means ``thin'' in French). Set $$\CR_n^-(\C)=\Z\omega
Cat(I^n,\C)^-/(M_n^-(\C)+\de^- M_{n+1}^-(\C))$$ where $M_n^-(\C)+\de^-
M_{n+1}^-(\C)$ is the sub-$\Z$-module of $\Z\omega Cat(I^n,\C)^-$
generated by $M_n^-(\C)$ and the image of $M_{n+1}^-(\C)$ by $\de^-$. The
differential map $\de^-$ induces a differential map
$$\CR_{n+1}^-(\C)\longrightarrow \CR_{n}^-(\C)$$ This chain complex is
called the
\textit{reduced branching complex} of $\C$.  The homology associated to
this chain complex is denoted by $\HR_*^-(\C)$ and is
called the \textit{reduced branching  homology} of $\C$.
\ed

\bp Conjecture~\ref{thin} is equivalent to the following statement : if
$\C$ is the free $\omega$-category generated by a precubical set or by a
globular set, then the canonical map from the branching chain complex to
the reduced branching chain complex of $\C$ is a quasi-isomorphism. \ep

\bpf By the following short exact sequence of
chain complexes $$\xymatrix{0\fr{}& {M_*^-(\C)+\de^- M_{*+1}^-(\C)}\fr{}&
{\Z\omega Cat(I^*,\C)^-}\fr{} & {\CR_{*}^-(\C)}\fr{}&0}$$
the assumption $H_n^-(\C)\iso \HR_n^-(\C)$ for all $n$ is equivalent
to the acyclicity of the chain complex $(M_*^-+\de^- M_{*+1}^-,\de^-)$
(notice that $M_0^-(\C)=M_1^-(\C)=0$).

Now if Conjecture~\ref{thin} holds, then take an element
$x\in M_n^-(\C)+\de^-M_{n+1}^-(\C)$ which is a cycle. Then $x=t_1+\de^- t_2$
where $t_1\in M_n^-(\C)$ and $t_2\in M_{n+1}^-(\C)$. Then $t_1$
is a cycle in $\Z\omega Cat(I^n,\C)^-$ and a linear combination of thin
elements. Therefore $t_1$ is a cycle in  $\Z\omega Cat(I^n,\tr_{n-1}\C)^-$.
By Conjecture~\ref{thin}, $t_1=\de^- t_3$ where
$t_3\in \Z\omega Cat(I^{n+1},\tr_{n-1}\C)^-$.
Therefore $t_1\in \de^-M_{n+1}^-(\C)$. Conversely, suppose that the
sub-complex generated by the thin elements is acyclic. Take a cycle $t$
of $\Z\omega Cat(I^n,\C)^-$ which is a linear combination of thin elements.
Then $t$ is a cycle of $M_n^-(\C)+\de^-M_{n+1}^-(\C)$, therefore
there exists $t_1\in M_{n+1}^-(\C)$ and $t_2\in M_{n+2}^-(\C)$ such
that $t=\de^-(t_1+\de^- t_2)=\de^- t_1$.
\epf

\bd Let $x$ and $y$ be two elements of $\Z\omega Cat(I^n,\C)^-$. Then
$x$ and $y$ are \textit{T-equivalent} (T for thin) if the
corresponding elements in the reduced branching complex are
equal, that means if $x-y\in M_n^-(\C)+\de^- M_{n+1}^-(\C)$.  This defines an
equivalence relation on $\Z\omega Cat(I^n,\C)^-$ indeed. \ed

\section{Matrix notation for higher dimensional composition in the cubical singular nerve}\label{matrix}

There exists on the cubical nerve $\omega Cat(I^*,\C)$ of an
$\omega$-category $\C$ a structure of cubical $\omega$-categories
\cite{Gau} by setting
\beas
&& \Gamma_i^-(x)(k_1\dots k_n)=x(k_1\dots \max(k_i,k_{i+1})\dots k_n)\\
&& \Gamma_i^+(x)(k_1\dots k_n)=x(k_1\dots \min(k_i,k_{i+1})\dots k_n)
\eeas
with the order $-<0<+$ and with the proposition-definition :

\bp\cite{Gau} Let $\C$ be a globular $\omega$-category.
For any strictly positive natural number $n$ and any $j$ between $1$
and $n$, there exists one and only one natural map $+_j$ from the set
of pairs $(x,y)$ of $\mathcal{N}^\square(\C)_n\p
\mathcal{N}^\square(\C)_n$ such that  $\de_j^+(x)=\de_j^-(x)$ to the
set $\mathcal{N}^\square(\C)_n$ which satisfies the following properties :
\beas
&&\de^-_j(x +_j y)=\de^-_j(x)\\
&&\de^+_j(x +_j y)=\de^+_j(x)\\
&&\de^\alpha_i(x +_j y)=\left\{
\begin{CD}
&&\de^\alpha_i(x) +_{j-1} \de^\alpha_i(y)\hbox{ if $i<j$}\\
&&\de^\alpha_i(x) +_{j} \de^\alpha_i(y)\hbox{ if $i>j$}\\ \end{CD}
\right.
\eeas
Moreover, these operations induce a structure of cubical
$\omega$-category on $\mathcal{N}^\square(\C)$.
\ep

The sum $(x +_i y) +_j (z +_i t)=(x +_j z) +_i (y +_j t)$ if there exists
will be denoted by
$$\left[\begin{array}{cc}
x & z\\
y & t\\
\end{array}\right] \coin{j}{i}
$$
and using this notation, one can write

\begin{itemize}
\item If $i=j$, $\Gamma_i^- (x +_j y)=\left[\begin{array}{cc}
\epsilon_{j+1}(y) & \Gamma_j^- (y)\\
\Gamma_j^- (x) & \epsilon_j(y)\\
\end{array}\right] \coin{j+1}{j}$
\item If $i=j$, $\Gamma_i^+ (x +_j y)=\left[\begin{array}{cc}
\epsilon_{j}(x) & \Gamma_j^+ (y)\\
\Gamma_j^+ (x) & \epsilon_{j+1}(x)\\
\end{array}\right] \coin{j+1}{j}$
\end{itemize}

The matrix notation can be generalized to any composition like
\[(a_{11}+_i \dots +_i a_{1n})+_j \dots +_j (a_{m1} +_i \dots +_i a_{mn})\]
whenever the sources and targets of the $a_{ij}$ match up in an
obvious sense (this is not necessarily true). In that case, the above
expression is equal by the interchange law to
\[(a_{11}+_j \dots +_j a_{m1})+_i\dots +_i (a_{1n} +_j \dots +_j a_{mn})\]
and  we can denote the common value by

\[\left[\begin{array}{ccc}
a_{m1}&\dots & a_{mn}\\
\vdots && \vdots \\
a_{11} & \dots & a_{1n}\\
\end{array}\right] \coin{i}{j}\]

In such a matrix, an element like $\epsilon_i x$ is denoted by
$\degh$.
An element like $\epsilon_j x$ is denoted by
$\degv$. In a situation where $i=j+1$,
an element like $\Gamma_j^-(x)$ is denoted by $\gamn$ and an
element like $\Gamma_{j}^+(x)$ is denoted by  $\gamp$.
An element like $\epsilon_j \epsilon_j x=
\epsilon_{j+1} \epsilon_j x$ is denoted by $\car$. With $i=j+1$,
we can verify some of the above formulae :

\[
\Gamma_j^- (x +_j y)=\left[\begin{array}{cc}
\degh&\gamn\\\gamn&\degv\\\end{array}\right]=\left[\begin{array}{cc}
\epsilon_{j+1}(y) & \Gamma_j^- (y)\\
\Gamma_j^- (x) & \epsilon_j(y)\\
\end{array}\right] \]

\[
\Gamma_j^+ (x +_j y)=\left[\begin{array}{cc}
\degv&\gamp\\\gamp&\degh\\\end{array}\right]=
\left[\begin{array}{cc}
\epsilon_{j}(x) & \Gamma_j^+ (y)\\
\Gamma_j^+ (x) & \epsilon_{j+1}(x)\\
\end{array}\right] \]

\bd\cite{Brown_cube}\cite{phd-Al-Agl}
A \textit{$n$-shell} in the cubical singular nerve is a family of
$2(n+1)$ elements $x_i^\pm $ of $\omega Cat(I^n,\C)$ such that
$\de_i^\alpha x_j^\beta= \de_{j-1}^\beta x_i^\alpha$ for $1\leqslant
i< j\leqslant n+1$ and $\alpha,\beta\in\{-,+\}$.  \ed

\bd
The $n$-shell $(x_i^\pm )$ is \textit{fillable}  if
\begin{enumerate}
\item the sets $\{x_i^{(-)^i},1\leqslant i\leqslant n+1\}$ and
$\{x_i^{(-)^{i+1}},1\leqslant i\leqslant n+1\}$ have each one exactly
one non-thin element and if the other ones are thin.
\item
if $x_{i_0}^{(-)^{i_0}}$ and $x_{i_1}^{(-)^{i_1+1}}$ are these two
non-thin  elements then there exists $u\in \C$ such that
$s_n(u)=x_{i_0}^{(-)^{i_0}}(0_n)$ and
$t_n(u)=x_{i_1}^{(-)^{i_1+1}}(0_n)$.
\end{enumerate}
\ed

The following proposition is an analogue of \cite{phd-Al-Agl} Proposition 2.7.3.

\bp\label{remplissage}\cite{Gau} Let $(x_i^\pm )$ be a fillable $n$-shell with
$u$ as above.  Then there exists one and only one element $x$ of
$\omega Cat(I^{n+1},\C)$ such that $x(0_{n+1})=u$, and for $1\leqslant
i\leqslant n+1$, and $\alpha\in\{-,+\}$ such that $\de_i^\alpha x=
x_i^\alpha$.  \ep

Proposition~\ref{remplissage} has a very important consequence concerning
the use of the above notations. In dimension $2$, an expression $A$ like (for
example)

$$\left[\begin{array}{ccccc}
\gamp & \gamn & \car & \degv & \car \\
\degh & x & \degh & y & \degh \\
\car & \gamp & \gamn & \gamp & \degh
\end{array}\right]\coin{2}{1}$$
is necessarily equal to

$$\left[\begin{array}{cc} x & y \\ \degv & \gamp
\end{array}\right]\coin{2}{1}$$ because the labels of the interior are
the same ($A(00)=(x +_2 y)(00)$) and because the shells of $1$-faces are equal 
($\de_1^-A=\de^-_1x$, $\de_1^+A=\de_1^+x +_1 \de_1^+y$, $\de_2^-A= \de_2^-x$, 
$\de_2^+A=\de_1^-y +_1 \de_2^+y$) : the dark lines represent degenerate elements which 
 are like mirrors reflecting rays of light. This is a
fundamental phenomenon to understand some of the calculations of this
work. Notice that $A\neq x +_2 y$ because $\de_1^-A\neq \de_1^-(x +_2 y)$.

All calculations involving these matrix notations are justified
because the Dawson-Par\'e condition holds in $2$-categories due to 
the existence of connections (see \cite{MR94c:18009} and
\cite{MR1694653}). The Dawson-Par\'e condition stands as follows :
\textit{suppose that a square $\alpha$ has a decomposition of one edge
$a$ as $a=a_1+_1 a_2$.  Then $\alpha$ has a compatible composition
$\alpha=\alpha_1 +_i \alpha_2$, i.e.  such that $\alpha_j$ has edge
$a_j$ for $j=1,2$}. This condition can be understood as a coherence
condition which ensures that all ``compatible'' tilings represent the
same object.

Let us mention that these special $2$-dimensional notations for connections and
degeneracies first appeared in \cite{first_use} and in \cite{dim2_2}.

\section{Relation between the simplicial nerve and the branching nerve}\label{relation_usual_nerve}

\bp\label{nerf_coin}\cite{Gau} Let $\C$ be an $\omega$-category and $\alpha\in\{-,+\}$.
We set $$\mathcal{N}^-_{n}(\C)=\omega Cat(I^{n+1},\C)^-$$ and for all
$n\geqslant 0$ and all $0\leqslant i\leqslant n$,
$$\xymatrix{{\de_i:\mathcal{N}^-_{n}(\C)}\fr{}&{\mathcal{N}^-_{n-1}(\C)}}$$
is the arrow $\de^-_{i+1}$, and
$$\xymatrix{{\epsilon_i:\mathcal{N}^-_{n}(\C)}\fr{}&{\mathcal{N}^-_{n+1}(\C)}}$$
is the arrow $\Gamma^-_{i+1}$. We obtain in this way a simplicial
set $$(\mathcal{N}^-_*(\C),\de_i,\epsilon_i)$$ called the
\textit{branching simplicial nerve}
of $\C$. The non normalized complex associated to it gives exactly the
branching homology of $\C$ (in degree greater than or equal to $1$).  The
map $\mathcal{N}^-$ induces a functor from $\omega Cat_1$ to the
category $Sets^{\Delta^{op}}$ of simplicial sets.  \ep

\subsection*{The globular $\omega$-category $\Delta^n$} 
Now let us recall the construction of the $\omega$-category called by
Street the $n$-th oriental \cite{oriental}. We use actually the
construction appearing in \cite{MR92j:18004}. Let $O^n$ be the set of
strictly increasing sequences of elements of $\{0,1,\dots,n\}$. A
sequence of length $p+1$ will be of dimension $p$.  If
$\sigma=\{\sigma_0<\dots <\sigma_p\}$ is a $p$-cell of $O^n$, then we
set $\de_j
\sigma=\{\sigma_0<\dots<\widehat{\sigma_j}<\dots<\sigma_k\}$. If $\sigma$
is an element of $O^n$, let $R(\sigma)$ be the subset of $O^n$
consisting of elements $\tau$ obtained from $\sigma$ by removing some
elements of the sequence $\sigma$ and let $R(\Sigma)=\bigcup_{\sigma\in
\Sigma}R(\sigma)$. Notice that $R(\Sigma \cup T)=R(\Sigma)\cup R(T)$.

\bth
There is one and only one $\omega$-category $\Delta^n$ such that
\begin{enumerate}
\item the underlying set of $\Delta^n$ is included in the set of
subsets of $O^n$
\item the underlying set of $\Delta^n$ contains
all subsets like $R(\sigma)$ where $\sigma$ runs over $O^n$
\item all elements of $\Delta^n$ are compositions of $R(\sigma)$ where $\sigma$ runs over $O^n$
\item for $\sigma$ $p$-dimensional with $p\geqslant 1$, one has
\beas
&&s_{p-1}(R(\sigma))=R\left(\{\de_j\sigma,j\hbox{ is even}\}\right)\\
&&t_{p-1}(R(\sigma))=R\left(\{\de_j\sigma,j\hbox{ is odd}\}\right)
\eeas
\item if $\Sigma$ and $T$ are two elements of $\Delta^n$ such that
$t_p(\Sigma)=s_p(T)$ for some $p$, then $\Sigma\cup T\in \Delta^n$ and
$\Sigma\cup T=\Sigma *_p T$.
\end{enumerate}
Moreover, all elements $\Sigma$ of $\Delta^n$
satisfy the equality $\Sigma=R(\Sigma)$.
\eth

If $\C$ is an $\omega$-category and if $x\in \omega Cat(\Delta^n,\C)$,
then consider the labeling of
the faces of respectively $\Delta^{n+1}$ and $\Delta^{n-1}$ defined by :
\begin{itemize}
\item $\epsilon_i(x)(\sigma_0<\dots<\sigma_{r})=x(\sigma_0<\dots<\sigma_{k-1}<\sigma_{k}-1<\dots<\sigma_{r}-1)$
 if $\sigma_{k-1}<i$ and $\sigma_{k}>i$.
\item $x(\sigma_0<\dots<\sigma_{k-1}<i<\sigma_{k+1}-1<\dots<\sigma_{r}-1)$
if $\sigma_{k-1}<i$, $\sigma_{k}=i$ and $\sigma_{k+1}>i+1$.
\item $x(\sigma_0<\dots<\sigma_{k-1}<i<\sigma_{k+2}-1<\dots<\sigma_{r}-1)$
if $\sigma_{k-1}<i$, $\sigma_{k}=i$ and $\sigma_{k+1}=i+1$.
\end{itemize}
and
\[
\de_i(x)(\sigma_0<\dots<\sigma_{s})=x(\sigma_0<\dots<\sigma_{k-1}<\sigma_{k}+1<\dots<\sigma_{s}+1)
\]
where $\sigma_k,\dots,\sigma_{s}\geq i$ and $\sigma_{k-1}<i$.

It turns out that $\epsilon_i(x)\in \omega Cat(\Delta^{n+1},\C)$ and
$\de_i(x)\in \omega Cat(\Delta^{n-1},\C)$. See \cite{May,Weibel} for
further information about simplicial sets. One has :

\bd\cite{oriental} The simplicial set
$(\omega Cat(\Delta^n,\C),\de_i,\epsilon_i)$ is called the
\textit{simplicial ner\-ve} $\mathcal{N}(\C)$ of the globular $\omega$-category
$\C$. The corresponding homology is denoted by $H_*(\C)$. \ed

\bd Let $\C$ be a non-contracting 
$\omega$-category. By definition, $\C$ is of
length at most $1$ if and only if for any morphisms $x$ and $y$ of
$\C$ such that $x *_0 y$ exists, then either $x$ or $y$ is
$0$-dimensional. \ed

\bth\label{calcul} Let $\C$ be an $\omega$-category  of length at most $1$.
Denote by $\P\C$ the unique $\omega$-category such that its set of $n$-morphisms is
exactly the set of $(n+1)$-morphisms of $\C$ for any $n\geqslant 0$
with an obvious definition of the source and target maps and of the
composition laws. Then one has the isomorphisms
$H_n(\P\C)\iso H_{n+1}^-(\C)$
for $n\geqslant 1$.
\eth

\bpf We give only a sketch of proof.  
By definition, $H_{n+1}^-(\C)=H_n(\mathcal{N}^-(\C))$
for $n\geqslant 1$. Because of the hypothesis on $\C$, every element
$x$ of $\omega Cat(I^{n+1},\C)^-$ is determined by the values of the
$x(k_1\dots k_{n+1})$ where $k_1\dots k_{n+1}$ runs over the set of
words on the alphabet $\{0,-\}$. It turns out that there is a
bijective correspondence between $O^n$ and the word of length $n+1$ on
the alphabet $\{0,-\}$ : if $\sigma_0<\dots<\sigma_p$ is an element of
$O^n$, the associated word of length $n+1$ is the word $m_0\dots
m_{n}$ such that $m_{\sigma_i}=0$ and if
$j\notin\{\sigma_0,\dots,\sigma_p\}$, then $m_j=-$. It is then straightforward
to check  that the simplicial structure of $\mathcal{N}^-(\C)$ is exactly the
same as the simplicial structure of $\omega Cat(\Delta^*,\P\C)$ in
strictly positive dimension~\footnote{The latter point is actually detailed in
\cite{sglob}.}.
\epf

The above proof together with Proposition~\ref{nerf_coin} gives a
new proof of the fact that if $x\in\omega Cat(\Delta^n,\C)$, the
labelings $\de_i(x)$ and $\epsilon_i(x)$ above defined yield
$\omega$-functors from $\Delta^{n-1}$ (resp. $\Delta^{n+1}$) to
$\C$.

\begin{figure}
\begin{center}
\includegraphics[width=10cm]{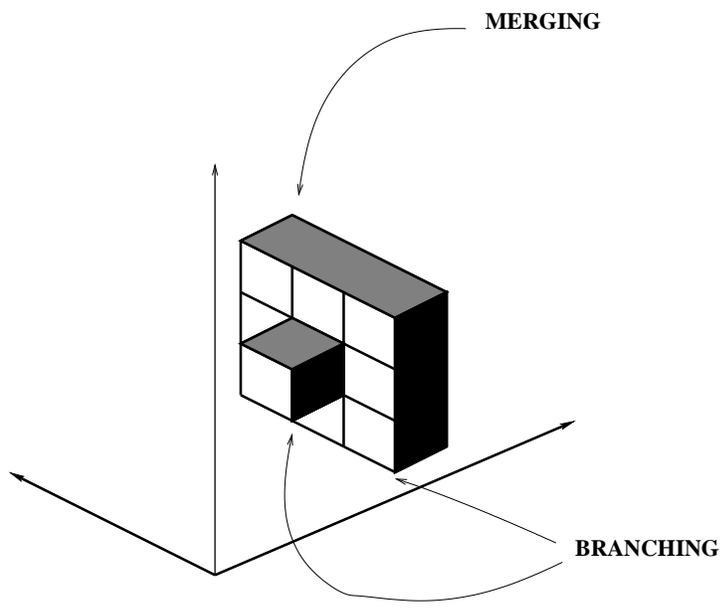}
\end{center}
\caption{A case where branching and merging homologies are not equal in dimension $2$}
\label{counterexample}
\end{figure}

Notice that the above proof also shows that $H_n(\P\C)\iso H_{n+1}^+(\C)$ 
where $H_*^+$ is the merging homology functor~\footnote{Like the branching 
nerve, the definition of the merging nerve needs to be slightly change, with 
respect to the definition given in \cite{Gau}. The correct definition is : an 
$\omega$-functor $x$ from $I^n$ to a non-contracting $\omega$-category $\C$ 
belongs to $\omega Cat(I^n,\C)^+$ if and only if for any $1$-morphism 
$\gamma$ of $I^n$ such that $t_0(\gamma)=R(+_n)$, then $x(\gamma)$ is a 
$1$-dimensional morphism of $\C$.} This means that for an $\omega$-category of length
at most $1$, $H_{n+1}^-(\C)\iso H_{n+1}^+(\C)$ for any $n\geqslant
1$. In general, this isomorphism is false as shown by
Figure~\ref{counterexample}. The precubical set we are considering in
this figure is the complement of the depicted obstacle.  Its branching
homology is $\Z\oplus \Z$ in dimension two, and its merging
homology is  $\Z$ in the same dimension.

The result $H_n(\P\C)\iso H_{n+1}^-(\C)\iso H_{n+1}^+(\C)$ for $\C$
of length at most one and for $n\geqslant 1$ also suggests that the
program of constructing the analogue in the computer-scientific
framework of usual homotopy invariants is complete for this kind of
$\omega$-categories. The simplicial set $\mathcal{N}(\P\C)$ together
with the graph obtained by considering the $1$-category generated by
the $1$-morphisms of $\C$ up to homotopy contain indeed all the
information about the topology of the underlying
automaton. Intuitively the simplicial set $\mathcal{N}(\P\C)$ is an
\textit{orthogonal section} of the automaton. Theorem~\ref{calcul} 
suggests that non-contracting $\omega$-categories of length at most one 
play a particular role
in this theory.  This idea will be deepened  in future works.

\bcor\label{astuce_calcul} With the same notation,
if $\P\C$ is the free globular $\omega$-category generated by a
composable pasting scheme in the sense of \cite{CPS}, then
$H_{n+1}^-(\C)$ vanishes for $n\geqslant 1$. \ecor

\bpf By \cite{MR99e:18008} Corollary 4.17 or by \cite{MR92j:18004}
Theorem 2.2, the simplicial nerve of the $\omega$-category of any
composable pasting scheme is contractible.  \epf

\bcor Let $2_p$ be the free $\omega$-category generated by a
$p$-morphism.  For any $p\geqslant 1$ and any $n\geqslant 1$,
$H^-_n(2_p)=0$.
\ecor

\bpf It is obvious for $n=1$ and for $n\geqslant 2$, $H^-_n(2_p)\iso
H_{n-1}(\P 2_p)$. But $\P 2_p=2_{p-1}$, therefore it suffices to
notice that the $(p-1)$-simplex is contractible. \epf

\bcor For any $n\geqslant 1$, let $G_n\AB$ be the $\omega$-category
generated by two $n$-mor\-phisms $A$ and $B$ satisfying
$s_{n-1}(A)=s_{n-1}(B)$ and $t_{n-1}(A)=t_{n-1}(B)$. Then
$$H_p^-(G_n\AB)=0$$ for $0<p<n$ or $p>n$ and
$$H_0^-(G_n\AB)=H_n^-(G_n\AB)=\Z.$$
\ecor

\bpf It suffices to calculate the simplicial homology of a simplicial
set homotopic to a $(n-1)$-sphere.  \epf

Let $S$ be a composable pasting scheme (see \cite{CPS} for the
definition and \cite{MR92j:18004} for additional explanations). A
reasonable conjecture is that the branching homology of the free
$\omega$-category $Cat(S)$ generated by any composable pasting scheme
$S$ vanishes in strictly positive dimension. By Conjecture~\ref{zero},
it would suffice for a given composable pasting scheme $S$ to
calculate the branching homology of the \textit{bilocalization} $Cat(S)[I,F]$ of
$Cat(S)$ with respect to its initial state $I$ and its final state
$F$, that is the sub-$\omega$-category of $Cat(S)$ which
consists of the $p$-morphisms $x$ with $p\geq 1$
such that $s_0x\in I$ and $t_0x\in F$ and of the $0$-morphism $I$ and $F$.
The question of the calculation of
$$H_{p+1}^-(Cat(S)[I,F])\iso H_p(\P Cat(S)[I,F])$$
for $p\geqslant 1$ seems
to be related to the existence of what Kapranov and Voevodsky call the
derived pasting scheme of a composable pasting scheme
\cite{MR92j:18004}.  It is in general not true that $\P Cat(S)[I,F]$
(denoted by $\Omega Cat(S)$ in their article) is the free
$\omega$-category generated by a composable pasting scheme. But we may
wonder  whether there is a ``free cover'' of $\Omega Cat(S)$
by some $Cat(T)$ for some composable pasting scheme $T$. This $T$
would be the derived pasting scheme of $S$.

As for the $n$-cube $I^n$, its derived pasting scheme is the
composable pasting scheme generated by the permutohedron
\cite{MR34:6767,MR81m:55010,symetrique_cube}.  Therefore one has

\bp\label{prem_cit_perm}
Denote by $I^n[-_n,+_n]$ the bilocalization of $I^n$ with respect to
its initial state $-_n$ and its final state $+_n$,  Then for all $p\geqslant 1$,
$H_p^-(I^n[-_n,+_n])=0$.
\ep

\bpf It is clear that
$H_1^-(I^n[-_n,+_n])=0$. For $p\geqslant 2$,
$H_p^-(I^n[-_n,+_n])\iso H_{p-1}(\Omega I^n)$
by Theorem~\ref{calcul}. But $\Omega I^n$ is the free
$\omega$-category generated by the permutohedron, and with
Corollary~\ref{astuce_calcul}, one gets $H_p^-(I^n[-_n,+_n])=0$ for
$p\geqslant 2$.  \epf

By filtrating the $1$-morphisms of $I^n$ by their length, it is possible
to construct a spectral sequence abutting to the branching homology of
$I^n$. More precisely a $1$-morphism $x$ is of length $\ell(x)$ if
$x=R(x_1)*_0\dots *_0 R(x_{\ell(x)})$ where
$x_1,\dots,x_{\ell(x)}\in (\underline{cub}^n)_1$. Now let
$F_p \omega Cat(I^*,I^n)^-$ be the subset of $x\in \omega Cat(I^*,I^n)^-$
such that for any $k_1\dots k_*\in (\underline{cub}^*)_1$ such that
$+\in\{k_1,\dots, k_*\}$, $\ell(x(k_1\dots k_*))\leqslant p$. Then
one gets a filtration on the branching complex of $I^n$ such that
$$F_{-1}\Z\omega Cat(I^*,I^n)^-\subset F_{0}\Z\omega Cat(I^*,I^n)^-\subset \dots \subset F_{n}\Z\omega Cat(I^*,I^n)^-$$
with
\beas
F_{-1}\Z\omega Cat(I^*,I^n)^-&=&0\\
F_{0}\Z\omega Cat(I^*,I^n)^-&=&\Z\omega Cat(I^*,I^n(-_n,+_n))^-\\
F_{n}\Z\omega Cat(I^*,I^n)^-&=&\Z\omega Cat(I^*,I^n)^-.
\eeas
One
has $E_{pq}^1=H_{p+q}(F_p \Z\omega Cat(I^*,I^n)^-/F_{p-1}\Z\omega
Cat(I^*,I^n)^-)\Longrightarrow H_{p+q}^-(I^n)$. By
Proposition~\ref{prem_cit_perm}, $E_{0q}=0$ if $q\neq 0$ and
$E_{00}=\Z$.

The above spectral sequence probably plays a role in the following
conjecture :

\subsection{Conjecture}\label{zero}
\textit{Let $\C$ be a finite $\omega$-category (that is such that  the
underlying set is finite). Let $I$ be the set of initial states of
$\C$ and let $F$ be the set of final states of $\C$ (then
$H_0^-(\C)=H_0^-(\C[I,F])=\Z F$). If for any $n>0$,
$H_n^-(\C[I,F])=0$, then for any $n>0$, $H_n^-(\C)=0$.}

By \cite{MR92j:18004}, $\Omega \Delta^n=I^{n-1}$, therefore the
vanishing of the branching homology of $I^{n-1}$ in strictly positive
dimension and Conjecture~\ref{zero} would enable to establish that
$H_p^-(\Delta^n)=0$ for $p>0$ and for any $n$.

\section{About folding operators}\label{construct_folding}

The aim of this section is to introduce an analogue in our framework
of the usual folding operators in cubical $\omega$-categories. First we
show how to recover the usual folding operators in our context.

The notations $\square_0$ or $\square_0^-$ (resp. $\square_1$ or
$\square_1^-$) correspond to the canonical map from $\C_0$ to $\omega
Cat(I^0,\C)$ (resp. from $\tr_1 \C$ to $\omega Cat(I^1,\C)$). Now let
us recall the construction of the operators $\square_n^-$ of
\cite{Gau}.

\bp\cite{Gau}\label{negfold} Let $\C$ be an $\omega$-category and let $n\geqslant 1$.
There exists one and only one natural map $\square_n^-$ from
$\tr_n\C$ to $\omega Cat(I^n,\C)$ such that the following
axioms hold :
\begin{enumerate}
\item one has $ev_{0_n}\square_n=Id_{\tr_n\C}$ where $ev_{0_n}(x)=x(0_n)$
is the label of the interior of $x$.
\item if $n\geqslant 3$ and $1\leqslant i\leqslant n-2$, then
$\de_i^\pm \square_n^-= \Gamma_{n-2}^- \de_i^\pm \square^-_{n-1} s_{n-1}$.
\item if $n\geqslant 2$ and $n-1\leqslant i \leqslant n$, then
$\de_i^- \square_n^-=\square_{n-1}^- d_{n-1}^{(-)^i}$ and
$\de_i^+ \square_n^-= \epsilon_{n-1} \de_{n-1}^+ \square_{n-1}^- s_{n-1}$.
\end{enumerate}
Moreover for $1\leqslant i \leqslant n$, we have
$\de_i^\pm \square_n^- s_n= \de_i^\pm \square_n^- t_n$ and if $x$ is
of dimension greater or equal than $1$, then
$\square_n^-(x)\in \omega Cat(I^n,\C)^-$.
\ep

\subsection{The usual folding operators}\label{folding_op}

One defines a natural map $\square_n$ from $\C_n$ to $\omega
Cat(I^n,\C)$ by induction on $n\geqslant 2$ as follows (compare with 
Proposition~\ref{negfold}).

\bp For any natural number $n$ greater or equal than $2$, there exists
a unique natural map $\square_n$ from $\C_n$ to $\omega Cat(I^n,\C)$ such
that
\begin{enumerate}
\item the equality $\square_n(x)(0_n)=x$ holds.
\item one has $\de_1^\alpha \square_n=\square_{n-1} d_{n-1}^{(-)^\alpha}$ for $\alpha=\pm $.
\item for $1<i\leqslant n$, one has $\de_i^\alpha \square_n = \epsilon_1
\de_{i-1}^\alpha \square_{n-1} s_{n-1}$.
\end{enumerate}
Moreover for $1\leqslant i \leqslant n$, we have $\de_i^\pm \square_n s_n u=
\de_i^\pm \square_n t_nu$ for any $(n+1)$-morphism $u$.
\ep

\bpf The induction equations define a fillable $(n-1)$-shell
(see Proposition~\ref{remplissage}). \epf

\bp For all $n\geqslant 0$, the evaluation map $ev_{0_n}:x\mapsto x(0_n)$
from $\omega Cat(I^n,\C)$ to $\C$ induces a bijection from $\gamma
\mathcal{N}^\square(\C)_n$ to $\tr_n\C$ where $\gamma$ is the functor
defined in \cite{phd-Al-Agl}. \ep

\bpf Obvious for $n=0$ and $n=1$. Recall that $\gamma$ is defined by
$$(\gamma G)_n=\{x\in G_n,\de_j^\alpha x\in
\epsilon_1^{j-1}G_{n-j}\hbox{ for $1\leqslant j\leqslant n$},
\alpha=0,1\}$$
Let us suppose that $n\geqslant 2$ and let us proceed by induction on
$n$. Since $ev_{0_n} \square_n (u)=u$ by the previous proposition,
then the evaluation map $ev$ from $\gamma\mathcal{N}^\square(\C)_n$ to
$\tr_n\C$ is surjective. Now let us prove that $x\in \gamma
\mathcal{N}^\square(\C)_n$ and $y\in \gamma
\mathcal{N}^\square(\C)_n$ and $x(0_n)=y(0_n)=u$ imply $x=y$.  Since
$x$ and $y$ are in $\gamma \mathcal{N}^\square(\C)_n$, then one sees
immediately that the four elements $\de_1^\pm x$ and $\de_1^\pm y$ are
in $\gamma \mathcal{N}^\square(\C)_{n-1}$. Since all other
$\de_i^\alpha x$ and $\de_i^\alpha y$ are thin, then
$\de_1^-x(0_{n-1})=\de_1^-y(0_{n-1})=s_{n-1}u$ and
$\de_1^+x(0_{n-1})=\de_1^+y(0_{n-1})=t_{n-1}u$.  By induction
hypothesis, $\de_1^-x=\de_1^-y=\square_{n-1}(s_{n-1}u)$ and
$\de_1^+x=\de_1^+y=\square_{n-1}(t_{n-1}u)$. By hypothesis, one can
set $\de_j^\alpha x=\epsilon_1^{j-1} x^\alpha_j$ and $\de_j^\alpha
y=\epsilon_1^{j-1} y^\alpha_j$ for $2 \leqslant j\leqslant n$. And one
gets $x_j^\alpha=(\de_1^\alpha)^{j-1} \de_j^\alpha
x=(\de_1^\alpha)^{j}x= (\de_1^\alpha)^{j}y=y_j^\alpha$. Therefore
$\de_j^\alpha x=\de_j^\alpha y$ for all $\alpha\in\{-,+\}$ and all
$j\in[1,\dots,n]$. By Proposition~\ref{remplissage}, one gets $x=y$.
\epf

The above proof shows also that the map which associates to any cube
$x$ of the cubical singular nerve of $\C$ the cube
$\square_{dim(x)}(x(0_{dim(x)}))$ is exactly the usual folding
operator as exposed in \cite{phd-Al-Agl}.

Unfortunately, these important operators are not internal to the
branching complex, due to the fact that an $n$-cube $x$ of the
cubical singular nerve is in the branching complex if and only
for any $1$-morphism $\gamma$ of $I^n$ starting from the initial state
$-_n$ of the $n$-cube, $x(\gamma)$ is $1$-dimensional (see
Definition~\ref{def_orientee}). But for example
$(\square_n(x(0_n)))(-\dots - 0)$ is $0$-dimensional.

\subsection{The negative folding operators}\label{folding_op2}

The idea of the negative folding operator $\Phi_n^-$ is
to ``concentrate'' a $n$-cube $x$ of the cubical singular nerve of an
$\omega$-category $\C$ to the faces $\delta_{n-1}^-(0_{n-1})$ and
$\delta_{n}^-(0_{n-1})$. Hence the following definition.

\bd Set $\Phi_n^-(x)=\square_n^-(x(0_n))$.  This operator is called the
$n$-dimensional negative folding operator.  \ed

It is clear that $x\in\omega Cat(I^n,\C)^-$ implies $\Phi_n^-(x)\in
\omega Cat(I^n,\C)^-$. Therefore $\Phi_n^-$ yields a map from $\omega
Cat(I^n,\C)^-$ to itself.

Since $\de_{n-1}^-\square_n^-=\square_{n-1}^- d_{n-1}^{(-)^{n-1}}$ and
$\de_{n}^-\square_n^-=\square_{n-1}^- d_{n-1}^{(-)^{n}}$, the effect
of $\square_n^-(x(0_n))$ is indeed to concentrate the faces of the
$n$-cube $x$ on the faces $\delta_{n-1}^-(0_{n-1})$ and
$\delta_{n}^-(0_{n-1})$. All the $(n-1)$-cubes
$\de_i^\alpha\square_n^-(x)$ for $(i,\alpha)\notin \{(n-1,-),(n,-)\}$
are thin. Of course there is not only one way of concentrating the faces
of $x$ on $\delta_{n-1}^-(0_{n-1})$ and $\delta_{n}^-(0_{n-1})$. But
in some way, they are all equivalent in the branching complex
(Corollary~\ref{identification}). We could decide also to concentrate
the $n$-cubes for $n\geqslant 2$ on the faces $\delta_{1}^-(0_{n-1})$
and $\delta_{2}^-(0_{n-1})$, or more generally to concentrate the
$n$-cubes on the faces $\delta_{p(n)}^-(0_{n-1})$ and
$\delta_{q(n)}^-(0_{n-1})$ where $p(n)$ and $q(n)$ would be  integers of
opposite parity for all $n\geqslant 2$. Let us end this section by
explaining precisely the structure of all these choices.

In an $\omega$-category, recall that $d_n^-=s_n$, $d_n^+=t_n$ and by
convention, let $d_\omega^-=d_\omega^+=Id$. All the usual axioms of
globular $\omega$-categories remain true with this convention and the partial
order $n<\omega$ for any natural number $n$.

If $x$ is an element of an $\omega$-category $\C$, we denote by
$\langle x \rangle$ the $\omega$-category generated by $x$. The underlying set
of $\langle x \rangle$ is $\{s_n x,t_n x,n\in\N\}$. We denote by $2_n$ any
$\omega$-category freely generated by one $n$-dimensional element.

Let $R(k_1\dots k_n)\in I^n$ a face. Denote by $ev_{k_1\dots k_n}$
the natural transformation from $\omega Cat(I^n,-)$ to $\tr_\omega$ which
maps $f$ to $f(R(k_1\dots k_n))$.

\bd Let $n\in\N$. Recall that $\tr_n$ is the forgetful functor from
$\omega$-categories to sets which associates to any $\omega$-category
its set of morphisms of dimension lower or equal than $n$ and let
$i_n$ be the inclusion functor from $\tr_{n-1}$ to $\tr_n$.  We call
cubification of dimension $n$, or $n$-cubification a natural
transformation $\square$ from $\tr_n$ to $\omega Cat (I^n,-)$. If moreover,
$ev_{0_n}\square=Id$, we say that the cubification is thick.  \ed

We  see immediately that $\square_n^-$, $\square_n$ (and $\square_n^+$
of \cite{Gau}) are examples of thick $n$-cubifica\-tions. By Yoneda the set of
$n$-cubifications is in bijection with the set of
$\omega$-functors from $I^n$ to $2_n$. So for a given $n$, there is a
finite number of $n$-cubifications.

\bp\label{transf} Let $f$ be a natural transformation from $\tr_m$ to $\tr_n$ with
$m,n\in \N\cup\{\omega\}$. Then there
exists $p\leq m$ and $\alpha\in\{-,+\}$ such that $f=d_p^\alpha$.
And necessarily, $p\leq I\!n\!f(m,n)$.
\ep

\bpf  Denote by
 $$\begin{CD}{\!\!}<{\!\!}A{\!\!}>{\!\!}=2_n @>g>> {\!\!}<{\!\!}B{\!\!}>{\!\!}=2_m\end{CD}$$
the $\omega$-functor
which corresponds to $f$ by Yoneda. Then $g(A)=d_p^\alpha(B)$ for some $p$
and some $\alpha$. And necessarily, $p\leq \min(m,n)$ (where the notation $min$
means the smallest element).
\epf

\bcor Let $\square$ be a $n$-cubification with $n\geq 1$ a natural
number. Then for any $i$ with $1\leq i\leq n$, $\de_i^\pm \square
s_n=\de_i^\pm \square t_n$.
\ecor

\bpf We have
$$\de_i^\pm \square s_n x(l_1\dots l_{n-1})= ev_{l_1\dots [\pm]_i\dots l_{n-1}} \square s_n(x)$$
But $ev_{l_1\dots [\pm]_i\dots l_{n-1}} \square$ is a natural transformation
from $\tr_n$ to $\tr_{n-1}$. By Proposition~\ref{transf}, we get
$$\de_i^\pm \square s_n x(l_1\dots l_{n-1})
=ev_{l_1\dots [\pm]_i\dots l_{n-1}} \square t_n(x)
= \de_i^\pm \square t_n x(l_1\dots l_{n-1})$$
\epf

We arrive at a theorem which explains the structure of all cubifications :

\bth\label{unicite} Let $\square$ be a thick $n$-cubification and let $f$ 
be an $\omega$-functor
from $I^{n+1}$ to $I^n$ such that $f(R(0_{n+1}))=R(0_n)$.
Denote by
$f^*$ the corresponding natural transformation from $\omega Cat(I^n,-)$
to $\omega Cat(I^{n+1},-)$. Then there
exists one and only one thick $(n+1)$-cubifica\-tion denoted by $f^*.\square$
such that for $1\leq i \leq n+1$,
$$(f^*.\square) i_{n+1}= f^* \square$$
where $i_{n+1}$ is the canonical natural transformation from $\tr_n$ to
$\tr_{n+1}$.
\eth

\bpf One has \beas \de _i^{\alpha} (f^*.\square) &=& \de _i^{\alpha}
(f^*.\square) d_{n}^{(-)^i} \\ &=& \de _i^{\alpha} (f^*.\square)
i_{n+1} d_{n}^{(-)^i}\\ &=& \de _i^{\alpha} f^* \square d_{n}^{(-)^i}
\eeas Therefore if $x\in \C_{n+1}$ for some $\omega$-category $\C$,
then $\de _i^{\alpha} (f^*.\square)x=\de_i^{\alpha} f^* \square
d_{n}^{(-)^i} x$ for $1\leq i \leq n+1$ and we obtain a fillable
$n$-shell in the sense of Proposition~\ref{remplissage}.  \epf

\bcor Let $u$ be an $\omega$-functor from $I^n$ to $2_n$ which
maps $R(0_n)$ to the unique $n$-mor\-phism of $2_n$ (we will say that $u$ is
thick because the corresponding cubification is also thick). Let $f$ be an $\omega$-functor
from $I^{n+1}$ to $I^n$ which maps $R(0_{n+1})$ to $R(0_n)$. Then there
exists one and only one thick $\omega$-functor $v$ from $I^{n+1}$ to
$2_{n+1}$ such that the following diagram commutes :
$$\xymatrix{{I^{n+1}}\fr{v} \fd{f} & {2_{n+1}}\fd{} \\ {I^n} \fr{u} & {2_n}}$$
the arrow from $2_{n+1}$ to $2_n$ being the unique $\omega$-functor which
sends the $(n+1)$-cell of $2_{n+1}$ to the $n$-cell of $2_n$.
\ecor

If $\square$ is a $n$-cubification and $f_i$ thick $\omega$-functors from
$I^{n+i+1}$ to $I^{n+i}$ for $0\leq i \leq p$ then we can denote without ambiguity
by $f_p.f_{p-1}.\dots .f_0.\square$ the $(n+p)$-cubification
$f_p.(f_{p-1}.(\dots .f_0.\square))$. Let us denote by  $\square_0$ the
unique $0$-cubification. We have the following formulas :

\bp\label{exemple_degenerate} Let $x\in \C$ be a $p$-dimensional morphism with $p \geq 1$
and let $n\geq p$. Then
$$\square_n^- x=\Gamma_{n-1}^-\dots \Gamma_{p}^- \square_p^- x$$
(by convention, the above formula is tautological for $n=p$)
\ep

\bpf We are going to show the formula by induction on $n$. The case $n=p$
is trivial.  If $i\leq n-1$, then $\de_i^\pm \square_{n+1}^- x=\Gamma_{n-1}^- \de_i^\pm \square_{n}^- x=\de_i^\pm \Gamma_{n}^- \Gamma_{n-1}^-\dots \Gamma_{p}^- \square_p^- x$.
And if $i\geq n$, then
$$\de_i^- \square_{n+1}^- x=\square_{n}^- x= \de_i^-
\Gamma_{n}^- \Gamma_{n-1}^-\dots \Gamma_{p}^- \square_p^- x$$
and
$$\de_i^+ \square_{n+1}^- x = \epsilon_n \de_i^+n \square_{n}^- x = \de_i^+ \Gamma_{n}^-
\Gamma_{n-1}^-\dots \Gamma_{p}^- \square_p^- x.$$ So the labelings
$\square_{n+1}^- x$ and $\Gamma_{n}^-\dots \Gamma_{p}^-
\square_p^- x$ of $I^{n+1}$ are the same ones. \epf

\bp For $n\geq 1$, we have
$\square_1^-=\epsilon_1.\square_0$ and
$$\square_n^-=\Gamma_{n-1}^- . \dots . \Gamma_{1}^-. \epsilon_1.\square_0$$
\ep

\bpf It is an immediate consequence of
Proposition~\ref{exemple_degenerate} and of the uniqueness of
Theorem~\ref{unicite}.  \epf

The converse of Theorem~\ref{unicite} is true. That is :

\bp Let $v$ be a thick $\omega$-functor from $I^{n+1}$ to $2_{n+1}$.
Then there exists an $\omega$-functor $f$ such that for any thick $\omega$-functor
$u$ from $I^n$ to $2_n$, the following diagram commutes\thinspace:
$$\xymatrix{{I^{n+1}}\fr{v} \fd{f} & {2_{n+1}}\fd{} \\ {I^n} \fr{u} & {2_n}}$$
\ep

\bpf Set $v(R(k_1\dots k_{n+1}))=d_{n_{k_1\dots
k_{n+1}}}^{\alpha_{k_1\dots k_{n+1}}}(A)$ where $\langle A \rangle=2_{n+1}$ and set
$\langle B \rangle=2_n$. By hypothesis, the equality $f(0_{n+1})=R(0_n)$ holds and
let
$$f(k_1\dots k_{n+1})=d_{n_{k_1\dots k_{n+1}}}^{\alpha_{k_1\dots
k_{n+1}}}(R(0_n))$$ Take any thick $\omega$-functor $u$ from $I^n$ to
$2_n$. Then \beas &&u\circ f(R(k_1\dots k_{n+1}))=u(d_{n_{k_1\dots
k_{n+1}}}^{\alpha_{k_1\dots k_{n+1}}}(R(0_n)))=d_{n_{k_1\dots
k_{n+1}}}^{\alpha_{k_1\dots k_{n+1}}} u(R(0_n))\\&& = d_{n_{k_1\dots
k_{n+1}}}^{\alpha_{k_1\dots k_{n+1}}}(B) \eeas
By
Proposition~\ref{remplissage}, it is clear that $f$ induces an
$\omega$-functor.  \epf

Here is an example of cubification : if the following picture
depicts the $3$-cube,

$$\xymatrix{
\fr{} & \fr{} & \\
\fu{} \fr{} & \fu{} \fr{} \ff{lu}{00-} & \ff{lu}{0+0}\fu{}\\
 & \fu{} \fr{}  & \ff{lu}{-00}\fu{}}
\Longrightarrow
\xymatrix{
\fr{} &  & \\
\fu{} \fr{} & \fu{} \fr{} \ff{lu}{+00} &\\
\fu{}\fr{} & \fu{} \fr{} \ff{lu}{0-0} & \ff{lu}{00+}\fu{}
}$$

we can represent a $3$-cubification $\square$ by indexing each face $k_1 k_2 k_3$
by the corresponding value of $ev_{k_1 k_2 k_3} \square i_3$ which is equal to $s_d$ or
$t_d$ for some $d$ between $0$ and $2$. So let us take $\square$ as follows :

$$\xymatrix{
\fr{t_0} & \fr{t_0} & \\
\fu{t_1} \fr{s_0} & \fu{t_1} \fr{t_1} \ff{lu}{t_1} & \ff{lu}{t_1}\fu{t_0}\\
 & \fu{s_0} \fr{s_0}  & \ff{lu}{s_2}\fu{s_1}}
\Longrightarrow
\xymatrix{
\fr{t_0} & & \\
\fu{t_0} \fr{t_0} & \fu{t_0} \fr{t_0} \ff{lu}{t_0} & \\
\fu{t_1}\fr{s_0} & \fu{t_1} \fr{s_1} \ff{lu}{t_1} & \ff{lu}{t_2}\fu{t_0}
}$$

We see that $\de_1^- \square i_3 = \Gamma_1^+ .\square_1$ and that
$\de_3^+ \square i_3 = \Gamma_1^- .\square_1$.

Now let us come back to our choice. It is not completely arbitrary
anyway because the operator $\square_n^-$ satisfies the following
important property : if $u$ is a $n$-morphism with $n\geqslant 2$,
then $\square_n^-(u)$ is a simplicial homotopy within the branching 
 nerve between $\square_{n-1}^-(s_{n-1}u)$ and
$\square_{n-1}^-(t_{n-1}u)$. Moreover, the family of cubifications
$(\square_n^-)_{n\geqslant 0}$ is the only family of cubifications
which satisfies this property because it is equivalent to defining a
$n$-shell for all $n$. However most of the results of the sequel can
be probably adapted to any family of $n$-cubification, provided that
they yield internal operations on the branching nerve (see
Conjecture~\ref{cubi1} and \ref{cubi2}).

\subsection{Characterization of the negative folding operators}

Now here is a useful property of the
folding operators :

\bth\label{caracterisation} Let $\C$ be an $\omega$-category. Let $x$ be an element of
$\mathcal{N}^\square_n(\C)$. Then the following two conditions are equivalent :
\begin{enumerate}
\item the equality $x=\Phi_n^-(x)$ holds
\item for $1\leqslant i \leqslant n$, one
has $\de_i^+ x\in Im(\epsilon_1^{n-1})$,
and for $1\leqslant i \leqslant n-2$, one
has $\de_i^- x\in Im(\Gamma_{n-2}^-\dots \Gamma_i^-)$.
\end{enumerate}
\eth

\bpf If $x=\Phi_n^-(x)$, then $x=\square_n^-(x(0_n))$ and by construction
of $\square_n^-$,
$$\de_i^-\square_n^-(x(0_n))=\Gamma_{n-2}^-\dots\Gamma_i^-\square_i^- d_i^{(-)^i}x(0_n)$$
for any $1\leqslant i \leqslant n-2$ and $\de^+_i x$ is
$0$-dimensional for any $1\leqslant i \leqslant n$.
For $n$ equal to $0$, $1$ or $2$, the converse is
obvious. Suppose the converse proved for $n-1\geqslant 2$ and let us
prove it by induction for $n\geqslant 3$. By hypothesis, as soon
as $+\in \{k_1,\dots,k_n\}$, then $x(k_1\dots k_n)$
is $0$-dimensional. For $1\leqslant i \leqslant n-3$, one
has
\beas
\de_i^-\de_{n}^-x&=&\de_{n-1}^- \de_i^- x\\
&=&\de_{n-1}^-\Gamma_{n-2}^-\dots
\Gamma_i^- Y_i\hbox{ for some $Y_i\in \mathcal{N}^\square_i(\C)$}\\
&=&\Gamma_{n-3}^-\dots \Gamma_i^- Y_i
\eeas
and
\beas
\de_i^-\de_{n-1}^-x&=&\de_{n-2}^- \de_i^- x\\
&=&\de_{n-2}^-\Gamma_{n-2}^-\dots
\Gamma_i^- Y_i\\
&=&\Gamma_{n-3}^-\dots \Gamma_i^- Y_i
\eeas
therefore $\de_{n-1}^-x$ and $\de_{n}^-x$ satisfy the induction
hypothesis. So $\de_{n-1}^-x=\Phi_{n-1}^-(\de_{n-1}^-x)$
and $\de_{n}^-x=\Phi_{n-1}^-(\de_{n}^-x)$. Since the $\de_i^-x$ are
thin $(n-1)$-cubes  for all $i$ between $1$ and $n-2$,
then $d_{n-1}^{(-)^{n-1}}(x(0_n))= \de_{n-1}^-x(0_{n-1})$
and $d_{n-1}^{(-)^{n}}(x(0_n))= \de_{n}^-x(0_{n-1})$.
Therefore $\de_{n-1}^-x=\square_{n-1}^-(d_{n-1}^{(-)^{n-1}}(x(0_n)))$
and $\de_{n}^-x=\square_{n-1}^-(d_{n-1}^{(-)^{n}}(x(0_n)))$.
For $1\leqslant i \leqslant n-2$, one has
\beas
Y_i&=&\de_i^-\dots \de_{n-2}^-\de_i^- x\\
&=&\de_i^- \dots \de_{n-1}^- x\\
&=&\de_i^- \dots \de_{n-2}^-
\square_{n-1}^-(d_{n-1}^{(-)^{n-1}}(x(0_n)))\\
&=&\de_i^- \dots \de_{n-3}^-
\square_{n-2}^-(d_{n-2}^{(-)^{n-2}}(x(0_n)))\\
&=&(\dots)\\
&=&\square_{i}^-(d_{i}^{(-)^{i}}(x(0_n)))
\eeas
therefore an easy calculation shows that $x=\square_n^-(x(0_n))$.
\epf

\bcor The folding operator $\Phi_n^-$ is  idempotent.
\ecor

The end of this section is devoted to the description of $\Phi_2^-$
and $\Phi_3^-$. Since $\de_1^-\square_2^- = \square_1 s_1$ and
$\de_2^-\square_2^- = \square_1 t_1$, then one has for any
$\omega$-functor $x$ from $I^2$ to $\C$
\[\Phi_2^-(x)=\left[\begin{array}{cc}
\gamn & \car\\
x & \gamn \\
\end{array}\right] \coin{2}{1}\]

If $x$ is an $\omega$-functor from $I^3$ to $\C$, then

$$s_2 (x(000))=\left[\begin{array}{cc}
\de_3^- x& \de_2^+ x\\
\gamp & \de_1^- x \\
\end{array}\right]\coin{2}{1}\  \ (00) $$

because the $2$-source of $R(000)$ in $I^3$ looks like
\[\left[\begin{array}{cc}
R(00-)& R(0+0)\\
\gamp & R(-00) \\
\end{array}\right] \]

and

$$t_2 (x(000))=\left[\begin{array}{cc}
\de_1^+ x& \gamn\\
\de_2^- x & \de_3^+ x \\
\end{array}\right]\coin{2}{1}\ \ \  (00) $$

because the $2$-target of $R(000)$ in $I^3$  looks like
\[\left[\begin{array}{cc}
R(+00)& \gamn\\
R(0-0) & R(00+) \\
\end{array}\right]\]

So by convention, an element $x$ of $\omega Cat(I^3,\C)$ will be represented as follows
\[x =
\begin{array}{|c|c|}
\hline
A & B\\
\hline
& C\\
\hline\end{array} \begin{array}{c}G \\\Longrightarrow\end{array}
\begin{array}{|c|c|}
\hline
D & \\
\hline
E & F\\
\hline\end{array}
\]
where $A=\de_3^-x$, $B=\de_2^+x$, $C=\de_1^-x$, $D=\de_1^+x$, $E=\de_2^-x$,
$F=\de_3^+x$ and $x(000)=G$\footnote{Beware of the fact that $A,\dots F$ are 
elements of the cubical singular nerve whereas $G$ is an element of the 
$\omega$-category we are considering.}.

With this convention, $\Gamma_1^- y$ for $y\in \omega Cat(I^2,\C)$ is
equal to
\[
\begin{array}{|c|c|}
\hline
\gamn & \degv\\
\hline
& y\\
\hline\end{array} \begin{array}{c}y(00) \\\Longrightarrow\end{array}
\begin{array}{|c|c|}
\hline
\degv & \\
\hline
y & \gamn \\
\hline\end{array}
\]
One has $\de_1^\pm \square_3^-=\Gamma_1^- \de_1^\pm \square_2^- s_2$,
$\de_2^-\square_3^- = \square_2^- t_2$, $\de_3^-\square_3^- = \square_2^- s_2$, $\de_1^+
\square_3^- =\de_2^+\square_3^- = \de_3^+\square_3^- = \square_2^-t_0$ by definition of
$\square_3^-$. Therefore
\beas
\de_i^+ \Phi_3^-(x)&=&\square_2^- t_0(G)\\
\de_2^- \Phi_3^-(x)&=&\square_2^- t_2(G)\\
\de_3^- \Phi_3^-(x)&=&\square_2^- s_2(G)
\eeas
and
\beas
&&\de_1^\pm \Phi_3^-(x)= \Gamma_1^- \de_1^\pm \square_2^- \left[\begin{array}{cc}
A & B\\
\gamp & C \\
\end{array}\right] \coin{2}{1}\\
&& = \Gamma_1^- \de_1^\pm
\left[\begin{array}{cccc}
\degh & \gamn & \car & \car \\
\gamn & \degv & \car & \car \\
A & B & \degh & \gamn \\
\gamp & C & \gamn & \degv \\
\end{array}\right] \coin{2}{1}\\
&&=\left\{ \begin{array}{c}
\Gamma_1^- (\de_1^- C+_1 \de_2^+ C +_1 \de_2^+ B)\hbox{ \hfill in the negative case}\\
\Gamma_1^- \de_1^+\de_1^+ B\hbox{ \hfill in the positive case} \end{array}\right.
\eeas
So if $x$ is the above $\omega$-functor from $I^3$ to $\C$, then
\[
\Phi_3^-(x)=\begin{array}{|c|c|}
\hline\begin{array}{cccc}
\degh & \gamn & \car & \car \\
\gamn & \degv & \car & \car \\
A & B & \degh & \gamn \\
\gamp & C & \gamn & \degv \\
\end{array}
& \car \\
\hline
& \gamn
\\
\hline\end{array} \begin{array}{c}G \\\Longrightarrow\end{array}
\begin{array}{|c|c|}
\hline
\car & \\
\hline
\begin{array}{cccc}

\gamn & \car  & \car \\
 D    & \gamn & \car \\
E   &  F   & \gamn \\
\end{array}
& \car\\
\hline\end{array}
\]

\section{Elementary moves in the cubical singular nerve}\label{elem_move}

In this section, the folding operators $\Phi_n^-$ are decomposed
in elementary moves. First of all, here is a definition.

\bd The elementary moves in the $n$-cube are one of the following
operators (with $1\leqslant i \leqslant n-1$ and $x\in \omega Cat(I^n,\C)$) :
\begin{enumerate}
\item $\ {}^v\psi^-_i x = \left[\begin{array}{c}\gamn \\ x\end{array}\right]\coin{i+1}{i}$
\item $\ {}^v\psi^+_i x = \left[\begin{array}{c}x \\ \gamp \end{array}\right]\coin{i+1}{i}$
\item $\ {}^h\psi^-_i x = \left[\begin{array}{cc} x & \gamn \end{array}\right]\coin{i+1}{i}$
\item $\ {}^h\psi^+_i x = \left[\begin{array}{cc} \gamp & x   \end{array}\right]\coin{i+1}{i}$
\end{enumerate}
\ed

\subsection*{Notation} 
One sets $\theta_i^-=\ {}^v\psi^-_{i+1} \ {}^v\psi^+_i$. This operator
plays a central r\^ole in the sequel.

Proposition~\ref{formule_3D}
expresses the elementary moves using the notation of the previous
paragraph (only the operators used in the sequel are calculated).

\bp\label{formule_3D} Let
\[x =
\begin{array}{|c|c|}
\hline
A & B\\
\hline
& C\\
\hline\end{array} \begin{array}{c}G \\\Longrightarrow\end{array}
\begin{array}{|c|c|}
\hline
D & \\
\hline
E & F\\
\hline\end{array}
\]
be an element of $\omega Cat(I^3,\C)$. Then one has
\[
\ {}^v\psi^+_1 x =
\begin{array}{|c|c|}
\hline
 \begin{array}{c}A\\\gamp\end{array} & \begin{array}{c}B\\C\end{array}\\
\hline
& \degv \\
\hline\end{array} \begin{array}{c}G \\\Longrightarrow\end{array}
\begin{array}{|c|c|}
\hline
\begin{array}{c}\degv\\D\end{array} & \\
\hline
\begin{array}{c}E\\\degv\end{array} & \begin{array}{c}F\\\gamp\end{array}\\
\hline\end{array}\]
\[\ {}^v\psi^+_2 x = \begin{array}{|c|c|}
\hline
 \begin{array}{cc}\degh & A\end{array} & B\\
\hline
& \begin{array}{c}C\\\gamp \end{array} \\
\hline\end{array} \begin{array}{c}G \\\Longrightarrow\end{array}
\begin{array}{|c|c|}
\hline
\begin{array}{c}D\\\gamp \end{array} & \\
\hline
\degh & \begin{array}{cc}E & F\end{array} \\
\hline\end{array}\]
\[\ {}^v\psi^-_2 x =
\begin{array}{|c|c|}
\hline
\begin{array}{cc} A & B\end{array}& \degh \\
\hline
& \begin{array}{c}\gamn \\ C\end{array}\\
\hline\end{array} \begin{array}{c}G \\\Longrightarrow\end{array}
\begin{array}{|c|c|}
\hline
\begin{array}{c}\gamn \\ D\end{array} & \\
\hline
E & \begin{array}{cc} F & \degh \end{array}\\
\hline\end{array}\]
\[\ {}^v\psi^-_1 x =
\begin{array}{|c|c|}
\hline
\begin{array}{c}\gamn \\ A \end{array} & \begin{array}{c}\degv \\ B  \end{array}\\
\hline
& C\\
\hline\end{array} \begin{array}{c}G \\\Longrightarrow\end{array}
\begin{array}{|c|c|}
\hline
\degv & \\
\hline
\begin{array}{c} D\\ E\end{array} & \begin{array}{c} \gamn \\ F\end{array}\\
\hline\end{array}\]
\[\ {}^h\psi_1^- x =
\begin{array}{|c|c|}
\hline
 \begin{array}{cc} A & \gamn \end{array}& \degv \\
\hline
& \begin{array}{c}B  \\ C\end{array}\\
\hline\end{array} \begin{array}{c}G \\\Longrightarrow\end{array}
\begin{array}{|c|c|}
\hline
\begin{array}{c}\degv  \\ D\end{array} & \\
\hline
E & \begin{array}{cc} F & \gamn \end{array}\\
\hline\end{array}\]
\[\ {}^h\psi_2^- x=\begin{array}{|c|c|}
\hline
A & \begin{array}{cc} B & \degh \end{array}\\
\hline
& \begin{array}{cc} C & \gamn \end{array}\\
\hline\end{array} \begin{array}{c}G \\\Longrightarrow\end{array}
\begin{array}{|c|c|}
\hline
\begin{array}{cc} D & \gamn \end{array} & \\
\hline
\begin{array}{cc} E & F  \end{array} & \degh \\
\hline\end{array}\]
\[\theta_1^- x= \begin{array}{|c|c|}
\hline
\begin{array}{cc} A & B\\ \gamp & C\end{array} & \begin{array}{c}\degh
\\ \degh \end{array}\\
\hline
& \begin{array}{c}\gamn
\\ \degv \end{array}\\
\hline\end{array} \begin{array}{c}G \\\Longrightarrow\end{array}
\begin{array}{|c|c|}
\hline
\begin{array}{c}\gamn
\\ D \end{array} & \\
\hline
\begin{array}{c}E
\\ \degv \end{array} & \begin{array}{cc}F & \degh \\ \gamp & \degh\end{array}\\
\hline\end{array}
\]
\ep

\bpf
One has $\ {}^v\psi^+_i x=\Gamma_i^+ \de_i^-x +_i x$. Therefore
\beas
&&\de_1^- \ {}^v\psi^+_1 x= \epsilon_1\de_1^-  \de_1^-x\\
&&\de_1^+ \ {}^v\psi^+_1(x)= \de_1^+ x\\
&&\de_2^- \ {}^v\psi^+_1 x= \epsilon_1 \de_1^- \de_1^-x +_1 \de_2^-x\\
&&\de_2^+ \ {}^v\psi^+_1 x = \de_1^-x +_1 \de_2^+ x\\
&&\de_3^\pm \ {}^v\psi^+_1 x= \de_3^\pm \Gamma_1^+ \de_1^-x +_1 \de_3^\pm x
\eeas
So one has
\[\ {}^v\psi^+_1 x =
\begin{array}{|c|c|}
\hline
 \begin{array}{c}A\\\gamp\end{array} & \begin{array}{c}B\\C\end{array}\\
\hline
& \degv \\
\hline\end{array} \begin{array}{c}G \\\Longrightarrow\end{array}
\begin{array}{|c|c|}
\hline
\begin{array}{c}\degv\\D\end{array} & \\
\hline
\begin{array}{c}E\\\degv\end{array} & \begin{array}{c}F\\\gamp\end{array}\\
\hline\end{array}
\]
And
\beas
&&\de_1^- \ {}^v\psi^+_2 x= \Gamma_1^+ \de_1^- \de_2^-x +_1 \de_1^- x\\
&&\de_2^- \ {}^v\psi^+_2 x= \epsilon_2 \de_2^- \de_2^-x\\
&&\de_3^- \ {}^v\psi^+_2 x= \epsilon_2 \de_2^- x +_2 \de_3^- x\\
&&\de_1^+ \ {}^v\psi^+_2 x= \Gamma_1^+ \de_1^+ \de_2^-x +_1 \de_1^+ x\\
&&\de_2^+ \ {}^v\psi^+_2 x= \de_2^+ x\\
&&\de_3^+ \ {}^v\psi^+_2 x = \de_2^-x +_2 \de_3^+ x
\eeas
Consequently one has
\[
\ {}^v\psi^+_2 x = \begin{array}{|c|c|}
\hline
 \begin{array}{cc}\degh & A\end{array} & B\\
\hline
& \begin{array}{c}C\\\gamp \end{array} \\
\hline\end{array} \begin{array}{c}G \\\Longrightarrow\end{array}
\begin{array}{|c|c|}
\hline
\begin{array}{c}D\\\gamp \end{array} & \\
\hline
\degh & \begin{array}{cc}E & F\end{array} \\
\hline\end{array}
\]
One has $\ {}^v\psi^-_i x=  x +_i \Gamma_i^-\de_i^+x$. Therefore
\beas
&&\de_1^- \ {}^v\psi^-_2 x = \de_1^- x +_1 \de_1^- \Gamma_2^- \de_2^+x\\
&&\de_1^+ \ {}^v\psi^-_2 x =\de_1^+ x +_1 \de_1^+ \Gamma_2^- \de_2^+x\\
&&\de_2^- \ {}^v\psi^-_2 x = \de_2^- x\\
&&\de_2^+ \ {}^v\psi^-_2 x=\epsilon_2\de_2^+\de_2^+ x\\
&&\de_3^- \ {}^v\psi^-_2 x = \de_3^- x+_2\de_2^+ x\\
&&\de_3^+ \ {}^v\psi^-_2 x = \de_3^+ x +_2 \epsilon_2\de_2^+\de_2^+ x
\eeas
So
\[\ {}^v\psi^-_2 x =
\begin{array}{|c|c|}
\hline
\begin{array}{cc} A & B\end{array}& \degh \\
\hline
& \begin{array}{c}\gamn \\ C\end{array}\\
\hline\end{array} \begin{array}{c}G \\\Longrightarrow\end{array}
\begin{array}{|c|c|}
\hline
\begin{array}{c}\gamn \\ D\end{array} & \\
\hline
E & \begin{array}{cc} F & \degh \end{array}\\
\hline\end{array}
\]
And
\beas
&&\de_1^- \ {}^v\psi^-_1 x = \de_1^- x \\
&&\de_1^+ \ {}^v\psi^-_1 x =\epsilon_1\de_1^+ \de_1^+ x\\
&&\de_2^- \ {}^v\psi^-_1 x = \de_2^- x+_1 \de_1^+ x\\
&&\de_2^+ \ {}^v\psi^-_1 x=\de_2^+ x +_1 \epsilon_1\de_1^+\de_1^+ x\\
&&\de_3^- \ {}^v\psi^-_1 x = \de_3^- x+_1\Gamma_1^- \de_2^-\de_1^+ x\\
&&\de_3^+ \ {}^v\psi^-_1 x = \de_3^+ x +_1 \Gamma_1^- \de_2^+\de_1^+ x
\eeas
therefore
\[\ {}^v\psi_1^- x=\begin{array}{|c|c|}
\hline
\begin{array}{c}\gamn \\ A \end{array} & \begin{array}{c}\degv \\ B  \end{array}\\
\hline
& C\\
\hline\end{array} \begin{array}{c}G \\\Longrightarrow\end{array}
\begin{array}{|c|c|}
\hline
\degv & \\
\hline
\begin{array}{c} D\\ E\end{array} & \begin{array}{c} \gamn \\ F\end{array}\\
\hline\end{array}
\]
One has $\ {}^h\psi_1^- x= x+_2 \Gamma_1^- \de_2^+ x$. Then 
\beas
&&\de_1^- \ {}^h\psi_1^- x= \de_1^- x +_1 \de_2^+ x\\
&&\de_1^+ \ {}^h\psi_1^- x= \de_1^+ x\\
&&\de_2^- \ {}^h\psi_1^- x= \de_2^- x\\
&&\de_2^+ \ {}^h\psi_1^- x= \epsilon_1\de_1^+ \de_2^+  x\\
&&\de_3^- \ {}^h\psi_1^- x= \de_3^- x +_2 \Gamma_1^- \de_2^-\de_2^+ x\\
&&\de_3^+ \ {}^h\psi_1^- x= \de_3^+ x +_2 \Gamma_1^- \de_2^+\de_2^+
x
\eeas
So
\[\ {}^h\psi_1^- x =
\begin{array}{|c|c|}
\hline
 \begin{array}{cc} A & \gamn \end{array}& \degv \\
\hline
& \begin{array}{c}B  \\ C\end{array}\\
\hline\end{array} \begin{array}{c}G \\\Longrightarrow\end{array}
\begin{array}{|c|c|}
\hline
\begin{array}{c}\degv  \\ D\end{array} & \\
\hline
E & \begin{array}{cc} F & \gamn \end{array}\\
\hline\end{array}
\]
One has $\ {}^h\psi_2^- x= x+_3 \Gamma_2^- \de_3^+ x$. Therefore
\beas
&&\de_1^- \ {}^h\psi_2^- x= \de_1^- x +_2\Gamma_1^- \de_1^- \de_3^+ x\\
&&\de_1^+ \ {}^h\psi_2^- x= \de_1^+ x +_2 \Gamma_1^- \de_1^+\de_3^+ x\\
&&\de_2^- \ {}^h\psi_2^- x= \de_2^- x +_2 \de_3^+x \\
&&\de_2^+ \ {}^h\psi_2^- x= \de_2^+x +_2 \epsilon_2\de_2^+ \de_3^+  x\\
&&\de_3^- \ {}^h\psi_2^- x= \de_3^- x \\
&&\de_3^+ \ {}^h\psi_2^- x= \de_3^+ \Gamma_2^- \de_3^+x = \epsilon_2 \de_2^+ \de_3^+ x
\eeas
so
\[
\ {}^h\psi_2^- =
\begin{array}{|c|c|}
\hline
A & \begin{array}{cc} B & \degh \end{array}\\
\hline
& \begin{array}{cc} C & \gamn \end{array}\\
\hline\end{array} \begin{array}{c}G \\\Longrightarrow\end{array}
\begin{array}{|c|c|}
\hline
\begin{array}{cc} D & \gamn \end{array} & \\
\hline
\begin{array}{cc} E & F  \end{array} & \degh \\
\hline\end{array}
\]
Now let us calculate $\theta_1^- x$. One has

\beas
&& \theta_1^-x = \ {}^v\psi^-_2
\ {}^v\psi^+_1 x\\
&& =\ {}^v\psi^-_2 \begin{array}{|c|c|}
\hline
 \begin{array}{c}A\\\gamp\end{array} & \begin{array}{c}B\\C\end{array}\\
\hline
& \degv \\
\hline\end{array} \begin{array}{c}G \\\Longrightarrow\end{array}
\begin{array}{|c|c|}
\hline
\begin{array}{c}\degv\\D\end{array} & \\
\hline
\begin{array}{c}E\\\degv\end{array} & \begin{array}{c}F\\\gamp\end{array}\\
\hline\end{array}\\
&&=\begin{array}{|c|c|}
\hline
\begin{array}{cc} A & B\\ \gamp & C\end{array} & \begin{array}{c}\degh
\\ \degh \end{array}\\
\hline
& \begin{array}{c}\gamn
\\ \degv \end{array}\\
\hline\end{array} \begin{array}{c}G \\\Longrightarrow\end{array}
\begin{array}{|c|c|}
\hline
\begin{array}{c}\gamn
\\ D \end{array} & \\
\hline
\begin{array}{c}E
\\ \degv \end{array} & \begin{array}{cc}F & \degh \\ \gamp & \degh\end{array}\\
\hline\end{array}
\eeas
\epf

The following proposition describes some of the commutation relations
satisfied by the previous operators, the differential maps and the
connection maps.

\bp\label{commutation} The following equalities hold (with $\alpha\in \{-,+\}$) :
\bea
&& \de_j^\alpha \ {}^v\psi_i^- =\left\{
\begin{array}{c} \ {}^v\psi_{i-1}^- \de_j^\alpha \hbox{ if $j<i$}\\ \ {}^v\psi_{i}^- \de_j^\alpha
\hbox{ if $j>i+1$}\end{array}\right. \\
&& \de_j^\alpha \ {}^h\psi_i^- =\left\{
\begin{array}{c} \ {}^h\psi_{i-1}^- \de_j^\alpha \hbox{ if $j<i$}\\ \ {}^h\psi_{i}^- \de_j^\alpha
\hbox{ if $j>i+1$}\end{array}\right. \\
&& \de_j^\alpha \theta_i^- =\left\{
\begin{array}{c} \theta_{i-1}^- \de_j^\alpha \hbox{ if $j<i$}\\ \theta_{i}^- \de_j^\alpha
\hbox{ if $j>i+2$}\end{array}\right.\\
&& \theta_i^- \Gamma_j^- = \left\{
\begin{array}{c} \Gamma_j^- \theta_{i-1}^- \hbox{ if $j<i$}\\
\Gamma_j^- \theta_{i}^- \hbox{ if $j>i+2$}
\end{array}\right.\\
&& \de_i^- \ {}^v\psi_i^- = \de_i^-\\
&& \de_i^+ \ {}^v\psi_i^- = \epsilon_i \de_i^+ \de_i^+ \\
&& \de_{i+1}^- \ {}^v\psi_i^- = \de_{i+1}^- +_i \de_i^+ \\
&& \de_{i+1}^+ \ {}^v\psi_i^- = \de_{i+1}^+\\
&& \de_i^- \ {}^h\psi_i^- = \de_i^- +_i \de_{i+1}^+\\
&& \de_i^+ \ {}^h\psi_i^- = \de_i^+ \\
&& \de_{i+1}^- \ {}^h\psi_i^- = \de_{i+1}^-  \\
&& \de_{i+1}^+ \ {}^h\psi_i^- = \epsilon_i \de_i^+ \de_{i+1}^+\\
&& \de_i^- \theta_i^- = \Gamma_i^- \de_i^- \de_i^- \\
&& \de_i^+ \theta_i^- = \ {}^v\psi_i^- \de_i^+ \\
&& \de_{i+1}^- \theta_i^- = \de_{i+1}^- \\
&& \de_{i+1}^+ \theta_i^- = \epsilon_{i+1} \de_{i+1}^+ \de_i^- +_i \epsilon_{i+1} \de_{i+1}^+ \de_{i+1}^+\\
&& \de_{i+2}^- \theta_i^- = \left[ \begin{array}{cc} \de_{i+2}^- & \de_{i+1}^+ \\ \gamp & \de_i^-\end{array}\right] \coin{i+1}{i} \\
&& \de_{i+2}^+ \theta_i^- = \ {}^v\psi_i^+ \de_{i+2}^+\\
&& \theta_i^- \Gamma_i^- = \Gamma_{i+1}^-\\
&& \theta_i^- \Gamma_{i+1}^- = \Gamma_{i+1}^-
\eea
\ep

\bpf Equalities (1), (2), (3) and (4) are obvious. Equalities from (5)
to (12) are immediate consequences of the definitions.  With
Proposition~\ref{formule_3D}, one sees that

\beas
&& \de_1^- \theta_1^- = \Gamma_1^- \de_1^- \de_1^- \\
&& \de_1^+ \theta_1^- = \ {}^v\psi_1^- \de_1^+ \\
&& \de_{2}^- \theta_1^- = \de_{2}^- \\
&& \de_{2}^+ \theta_1^- = \epsilon_{2} \de_{2}^+ \de_1^- +_1 \epsilon_{2} \de_{2}^+ \de_{2}^+\\
&& \de_{3}^- \theta_1^- = \left[ \begin{array}{cc} \de_{3}^- & \de_{2}^+ \\ \gamp & \de_1^-\end{array}\right] \coin{2}{1} \\
&& \de_{3}^+ \theta_1^- = \ {}^v\psi_1^+ \de_{3}^+
\eeas

For a given $x$, the above equalities are equalities in the free
cubical $\omega$-category generated by $x$. Therefore, they 
depend only on the relative position of the indices $1$, $2$ and $3$
with respect to one another. Therefore, we can replace each index $1$
by $i$, each index $2$ by $i+1$ and each index 
$3$ by $i+2$ to obtain the required formulae.

In the same way, it suffices to prove the last two formulae in
lower dimension and for $i=1$. One has
\beas
\theta_1^- \Gamma_1^- x &=&\theta_1^- \begin{array}{|c|c|}
\hline
\gamn & \degv\\
\hline
& x\\
\hline\end{array} \begin{array}{c}x(00) \\\Longrightarrow\end{array}
\begin{array}{|c|c|}
\hline
\degv & \\
\hline
x & \gamn \\
\hline\end{array}\\
&=&\begin{array}{|c|c|}
\hline
\begin{array}{cc} \gamn & \degv\\ \gamp & x\end{array} & \begin{array}{c}\degh
\\ \degh \end{array}\\
\hline
& \begin{array}{c}\gamn
\\ \degv \end{array}\\
\hline\end{array} \begin{array}{c}x(00) \\\Longrightarrow\end{array}
\begin{array}{|c|c|}
\hline
\begin{array}{c}\gamn
\\ \degv \end{array} & \\
\hline
\begin{array}{c}x
\\ \degv \end{array} & \begin{array}{cc}\gamn & \degh \\ \gamp & \degh\end{array}\\
\hline\end{array}\\
&=&\Gamma_2^- x
\eeas
and
\beas
 \theta_1^- \Gamma_2^- x &=&\theta_1^- \begin{array}{|c|c|}
\hline
x & \degh\\
\hline
& \gamn\\
\hline\end{array} \begin{array}{c}x(00) \\\Longrightarrow\end{array}
\begin{array}{|c|c|}
\hline
\gamn & \\
\hline
x & \degh \\
\hline\end{array}\\
&=&\begin{array}{|c|c|}
\hline
\begin{array}{cc} x & \degh\\ \gamp & \gamn\end{array} & \begin{array}{c}\degh
\\ \car \end{array}\\
\hline
& \begin{array}{c}\gamn
\\ \degv \end{array}\\
\hline\end{array} \begin{array}{c}x(00) \\\Longrightarrow\end{array}
\begin{array}{|c|c|}
\hline
\begin{array}{c}\car
\\ \gamn \end{array} & \\
\hline
\begin{array}{c}x
\\ \degv \end{array} & \begin{array}{cc}\degh & \degh \\ \car & \car\end{array}\\
\hline\end{array}\\
&=&\Gamma_2^- x
\eeas

\epf

\bth\label{aspiration_face_plus} Set $\ {}^v\Psi_k^-=\ {}^v\psi_k^- \dots \ {}^v\psi_1^-$ and
$\ {}^h\Psi_k^-=\ {}^h\psi_k^- \dots \ {}^h\psi_1^-$. Then for $n\geqslant 2$ and
$1\leqslant i\leqslant n$,
one has
$$\de_i^+({}^v\Psi_1^- \ {}^h\Psi_1^-)\dots ({}^v\Psi_{n-1}^- \ {}^h\Psi_{n-1}^-)=\epsilon_1^{n-1}(\de_1^+)^n.$$
\eth

\bpf It is obvious for $n=2$. We are going to make an induction on $n$. Let $n\geqslant 2$
and $1\leqslant i\leqslant n$. Then
\beas
&& \de_i^+ ({}^v\Psi_1^- \ {}^h\Psi_1^-)\dots ({}^v\Psi_{n}^- \ {}^h\Psi_{n}^-)\\
&&=\epsilon_1^{n-1}(\de_1^+)^n({}^v\Psi_{n}^- \ {}^h\Psi_{n}^-)\\
&&=\epsilon_1^{n-1}(\de_1^+)^{n-1} \ {}^v\Psi_{n-1}^- \epsilon_1 \de_1^+ \de_1^+ \ {}^h\Psi_{n}^-\\
&&=\epsilon_1^{n-1}(\de_1^+)^{n-2}\ {}^v\Psi_{n-2}^-(\epsilon_1 \de_1^+ \de_1^+)^2 \ {}^h\Psi_{n}^-\\
&&=(\dots)\\
&&=\epsilon_1^{n-1}(\epsilon_1 \de_1^+ \de_1^+)^n \ {}^h\Psi_{n}^-\\
&&=\epsilon_1^{n}(\de_1^+)^{n+1}\ {}^h\Psi_{n}^-
\eeas
The equality
$$\de_i^+ ({}^v\Psi_1^- \ {}^h\Psi_1^-)\dots ({}^v\Psi_{n}^-
\ {}^h\Psi_{n}^-)x=\epsilon_1^{n}(\de_1^+)^{n+1}\ {}^h\Psi_{n}^-x$$ makes
sense if $x$ is a $(n+1)$-cube. And in this case,
$\epsilon_1^{n}(\de_1^+)^{n+1}\ {}^h\Psi_{n}^-x$ is $0$-dimensional and
$\epsilon_1^{n}(\de_1^+)^{n+1}\ {}^h\Psi_{n}^-x=\epsilon_1^{n}(\de_1^+)^{n+1}x$. This equality
holds in the free cubical $\omega$-category generated by $x$, and therefore
$$\epsilon_1^{n}(\de_1^+)^{n+1}\ {}^h\Psi_{n}^-=\epsilon_1^{n}(\de_1^+)^{n+1}.$$
Now suppose that $i=n+1$. Then
\beas
&&\de_{n+1}^+ ({}^v\Psi_1^- \ {}^h\Psi_1^-)\dots ({}^v\Psi_{n}^- \ {}^h\Psi_{n}^-)\\
&&=({}^v\Psi_1^- \ {}^h\Psi_1^-)\dots ({}^v\Psi_{n-1}^- \ {}^h\Psi_{n-1}^-)\de_{n+1}^+({}^v\Psi_{n}^- \ {}^h\Psi_{n}^-)\\
&&=({}^v\Psi_1^- \ {}^h\Psi_1^-)\dots ({}^v\Psi_{n-1}^- \ {}^h\Psi_{n-1}^-)({}^v\Psi_{n-1}^-\de_{n+1}^+\ {}^h\Psi_{n}^-)\\
&&=({}^v\Psi_1^- \ {}^h\Psi_1^-)\dots ({}^v\Psi_{n-1}^- \ {}^h\Psi_{n-1}^-)\ {}^v\Psi_{n-1}^-\de_{n+1}^+({}^h\psi_{n}^-\dots\ {}^h\psi_1^-)\\
&&=({}^v\Psi_1^- \ {}^h\Psi_1^-)\dots ({}^v\Psi_{n-1}^- \ {}^h\Psi_{n-1}^-)\ {}^v\Psi_{n-1}^- \epsilon_n\de_n^+ \de_{n+1}^+ ({}^h\psi_{n-1}^-\dots\ {}^h\psi_1^-)\\
&&=({}^v\Psi_1^- \ {}^h\Psi_1^-)\dots ({}^v\Psi_{n-1}^- \ {}^h\Psi_{n-1}^-)\ {}^v\Psi_{n-1}^- \epsilon_n\de_n^+({}^h\psi_{n-1}^-\dots\ {}^h\psi_1^-) \de_{n+1}^+\\
&&=({}^v\Psi_1^- \ {}^h\Psi_1^-)\dots ({}^v\Psi_{n-1}^- \ {}^h\Psi_{n-1}^-)\ {}^v\Psi_{n-1}^- \epsilon_n \epsilon_{n-1} \de_{n-1}^+ \de_n^+ ({}^h\psi_{n-2}^-\dots\ {}^h\psi_1^-) \de_{n+1}^+\\
&&=({}^v\Psi_1^- \ {}^h\Psi_1^-)\dots ({}^v\Psi_{n-1}^- \ {}^h\Psi_{n-1}^-)\ {}^v\Psi_{n-1}^- \epsilon_n \epsilon_{n-1} \de_{n-1}^+({}^h\psi_{n-2}^-\dots\ {}^h\psi_1^-)\de_n^+\de_{n+1}^+\\
&&=(\dots)\\
&&=({}^v\Psi_1^- \ {}^h\Psi_1^-)\dots ({}^v\Psi_{n-1}^- \ {}^h\Psi_{n-1}^-)\ {}^v\Psi_{n-1}^-\epsilon_n \dots \epsilon_1 \de_1^+\dots \de_{n+1}^+\\
&&=(\epsilon_1)^n(\de_1^+)^{n+1} \hbox{ for the same reason as above}
\eeas

\epf

\subsection*{Why does the proof of Theorem~\ref{aspiration_face_plus} work}

The principle of the proof of Theorem~\ref{aspiration_face_plus} is
the following observation (see in \cite{equiv_glob_cub})\thinspace: let
$f_1,\dots,f_n$ be $n$ operators such that (the product notation means
the composition)
\begin{enumerate}
\item for any $i$, one has $f_i  f_i = f_i$ (the operators $f_i$ are idempotent)
\item $|i-j|\geqslant 2$ implies $f_i f_j=f_j f_i$
\item $f_i f_{i+1}f_i =f_{i+1}f_i f_{i+1}$ for any $i$
\end{enumerate}
Then the operator $F=f_1(f_2 f_1)\dots (f_n f_{n-1}\dots f_1)$
satisfies $f_i F=F$ for any $i$. This means that $F$ enables to apply
all $f_i$ a maximal number of times. It turns out that the operators
$\ {}^v\psi_i^\pm$ and $\ {}^h\psi_i^\pm$ satisfy the above relations :

\bp The operators $\ {}^v\psi_i^\alpha$ and  $\ {}^h\psi_j^\beta$ are idempotent.
Moreover for any $i\geq 1$ and any $j\geq 1$, with $|i-j|\geq 2$, the following
equalities hold :
\bea
&& \ {}^v\psi_i^\alpha \ {}^h\psi_j^\beta = \ {}^h\psi_j^\beta \ {}^v\psi_i^\alpha \hbox{ for $\alpha\in\{-,+\}$} \label{evident1}\\
&& \ {}^v\psi_i^\alpha \ {}^h\psi_i^\alpha = \ {}^h\psi_i^\alpha \ {}^v\psi_i^\alpha\hbox{ for $\alpha\in\{-,+\}$}\label{evident2}\\
&& \ {}^h\psi_{i+1}^\alpha \ {}^v\psi_i^\alpha = \ {}^v\psi_i^\alpha \ {}^h\psi_{i+1}^\alpha\hbox{ for $\alpha\in\{-,+\}$}\label{moinsclair1}\\
&& \ {}^a\psi_i^\alpha \ {}^a\psi_{i+1}^\alpha\ {}^a\psi_i^\alpha = \ {}^a\psi_{i+1}^\alpha \ {}^a\psi_i^\alpha \ {}^a\psi_{i+1}^\alpha \hbox{ for $a \in\{v,h\}$ and $\alpha\in\{-,+\}$}
\label{interessant}
\eea
\ep

\bpf Equalities~\ref{evident1} and ~\ref{evident2} are obvious.

For the sequel, one can suppose $\alpha=-$. In the cubical singular
nerve of an $\omega$-category, two elements $A$ and $B$ of the same
dimension $n$ are equal if and only if $A(0_n)=B(0_n)$ and for $1\leq
k\leq n$ and $\alpha\in\{-,+\}$, one has $\de_k^\alpha A=\de_k^\alpha
B$.

Now we want to prove Equality~\ref{moinsclair1}.  Since
$({}^v\psi_i^\alpha x)(0_n)=({}^h\psi_i^\alpha x)(0_n)=x(0_n)$,
then $\ {}^h\psi_{i+1}^\alpha \ {}^v\psi_i^\alpha x = \
{}^v\psi_i^\alpha \ {}^h\psi_{i+1}^\alpha x$ for any $x$ of
dimension $n$ ($P_n$) is equivalent to $\de_k^\beta \
{}^h\psi_{i+1}^\alpha \ {}^v\psi_i^\alpha x = \de_k^\beta \
{}^v\psi_i^\alpha \ {}^h\psi_{i+1}^\alpha x$ for $1\leq k\leq n$
and $\beta\in\{-,+\} (E_{k,n})$. Proposition~\ref{commutation}
implies that $P_{n-1}\Longrightarrow E_{k,n}$ for $k<i$ or
$k>i+2$. For $k\in\{i,i+1,i+2\}$, proving Equality $E_{k,n}$ is
equivalent to proving it for the case $i=1$ and to replacing each
index $1$ by $i$, each index $2$ in by $i+1$ and
each index $3$ by $i+2$. And in the case $i=1$, the
equality is a calculation in the free cubical $\omega$-category
generated by $x$. So we can suppose that $x$ is of dimension as
low as possible. In our case, this equality makes sense if $x$ is
$3$-dimensional.  Therefore it suffices to verify
Equality~\ref{moinsclair1} in dimension $3$ for $i=1$.  And one
has \beas &&\ {}^h\psi_{2}^- \ {}^v\psi_1^-
\begin{array}{|c|c|} \hline
A & B\\
\hline
& C\\
\hline\end{array} \begin{array}{c}G \\\Longrightarrow\end{array}
\begin{array}{|c|c|}
\hline
D & \\
\hline
E & F\\
\hline\end{array}\\
&&= \ {}^h\psi_{2}^- \begin{array}{|c|c|}
\hline
\begin{array}{c}\gamn \\ A \end{array} & \begin{array}{c}\degv \\ B  \end{array}\\
\hline
& C\\
\hline\end{array} \begin{array}{c}G \\\Longrightarrow\end{array}
\begin{array}{|c|c|}
\hline
\degv & \\
\hline
\begin{array}{c} D\\ E\end{array} & \begin{array}{c} \gamn \\ F\end{array}\\
\hline\end{array}\\
&&=\begin{array}{|c|c|}
\hline
\begin{array}{c}\gamn \\ A \end{array} & \begin{array}{cc}\degv&\car  \\ B&\degh  \end{array}\\
\hline
& \begin{array}{cc}C & \gamn \end{array}\\
\hline\end{array} \begin{array}{c}G \\\Longrightarrow\end{array}
\begin{array}{|c|c|}
\hline
\begin{array}{cc}\degv&\car\end{array} & \\
\hline
\begin{array}{cc} D&\gamn \\ E& F\end{array} & \begin{array}{c} \car \\ \degh\end{array}\\
\hline\end{array}\\
&&= \ {}^v\psi_1^-
\begin{array}{|c|c|}
\hline
A & \begin{array}{cc} B & \degh \end{array}\\
\hline
& \begin{array}{cc} C & \gamn \end{array}\\
\hline\end{array} \begin{array}{c}G \\\Longrightarrow\end{array}
\begin{array}{|c|c|}
\hline
\begin{array}{cc} D & \gamn \end{array} & \\
\hline
\begin{array}{cc} E & F  \end{array} & \degh \\
\hline\end{array}\\
&&=\ {}^v\psi_1^- \ {}^h\psi_{2}^- \begin{array}{|c|c|}
\hline
A & B\\
\hline
& C\\
\hline\end{array} \begin{array}{c}G \\\Longrightarrow\end{array}
\begin{array}{|c|c|}
\hline
D & \\
\hline
E & F\\
\hline\end{array}
\eeas

In the same way, to prove Equality~\ref{interessant}, it suffices
to prove it for $i=1$ and in the $3$-dimensional case. And one has
\beas
&&\ {}^v\psi_1^- \ {}^v\psi_{2}^-\ {}^v\psi_1^- \begin{array}{|c|c|}
\hline
A & B\\
\hline
& C\\
\hline\end{array} \begin{array}{c}G \\\Longrightarrow\end{array}
\begin{array}{|c|c|}
\hline
D & \\
\hline
E & F\\
\hline\end{array}
\\
&&= \ {}^v\psi_1^- \ {}^v\psi_{2}^- \begin{array}{|c|c|}
\hline
\begin{array}{c}\gamn \\ A \end{array} & \begin{array}{c}\degv \\ B  \end{array}\\
\hline
& C\\
\hline\end{array} \begin{array}{c}G \\\Longrightarrow\end{array}
\begin{array}{|c|c|}
\hline
\degv & \\
\hline
\begin{array}{c} D\\ E\end{array} & \begin{array}{c} \gamn \\ F\end{array}\\
\hline\end{array}\\
&&=\ {}^v\psi_1^- \begin{array}{|c|c|}
\hline
\begin{array}{cc}\gamn & \degv \\ A &B\end{array} & \begin{array}{c}\car\\ \degh  \end{array}\\
\hline
& C\\
\hline\end{array} \begin{array}{c}G \\\Longrightarrow\end{array}
\begin{array}{|c|c|}
\hline
\begin{array}{c}\gamn \\\degv \end{array} & \\
\hline
\begin{array}{c} D\\ E\end{array} & \begin{array}{cc} \gamn & \car \\ F&\degh\end{array}\\
\hline\end{array}\\
&&=\begin{array}{|c|c|}
\hline
\begin{array}{cc}\degh & \gamn \\ \car& \degv\\\gamn & \degv \\ A &B\end{array} & \begin{array}{c}\car\\\car\\\car\\ \degh  \end{array}\\
\hline
& C\\
\hline\end{array} \begin{array}{c}G \\\Longrightarrow\end{array}
\begin{array}{|c|c|}
\hline
\begin{array}{c}\car \\\car \end{array} & \\
\hline
\begin{array}{c} \gamn\\ \degv\\ D\\ E\end{array} &
\begin{array}{cc} \car & \car \\ \car & \car \\ \gamn & \car \\ F&\degh\end{array}\\
\hline\end{array}
\eeas
and
\beas
&&\ {}^v\psi_2^- \ {}^v\psi_{1}^-\ {}^v\psi_2^- \begin{array}{|c|c|}
\hline
A & B\\
\hline
& C\\
\hline\end{array} \begin{array}{c}G \\\Longrightarrow\end{array}
\begin{array}{|c|c|}
\hline
D & \\
\hline
E & F\\
\hline\end{array}
\\
&&= \ {}^v\psi_2^- \ {}^v\psi_{1}^- \begin{array}{|c|c|}
\hline
\begin{array}{cc} A & B\end{array}& \degh \\
\hline
& \begin{array}{c}\gamn \\ C\end{array}\\
\hline\end{array} \begin{array}{c}G \\\Longrightarrow\end{array}
\begin{array}{|c|c|}
\hline
\begin{array}{c}\gamn \\ D\end{array} & \\
\hline
E & \begin{array}{cc} F & \degh \end{array}\\
\hline\end{array}\\
&&=\ {}^v\psi_2^- \begin{array}{|c|c|}
\hline
\begin{array}{cc} \degh & \gamn \\ \gamn & \degv \\ A & B\end{array}&
\begin{array}{c}\car \\ \car \\ \degh\end{array} \\
\hline
& \begin{array}{c}\gamn \\ C\end{array}\\
\hline\end{array} \begin{array}{c}G \\\Longrightarrow\end{array}
\begin{array}{|c|c|}
\hline
\begin{array}{c}\degv \\ \degv\end{array} & \\
\hline
\begin{array}{c} \gamn\\ D\\ E\end{array} & \begin{array}{cc} \car & \car \\ \gamn & \car \\F & \degh \end{array}\\
\hline\end{array}\\
&&= \begin{array}{|c|c|}
\hline
\begin{array}{cc}\degh & \gamn \\ \car& \degv\\\gamn & \degv \\ A &B\end{array} & \begin{array}{c}\car\\\car\\\car\\ \degh  \end{array}\\
\hline
& C\\
\hline\end{array} \begin{array}{c}G \\\Longrightarrow\end{array}
\begin{array}{|c|c|}
\hline
\begin{array}{c}\car \\\car \end{array} & \\
\hline
\begin{array}{c} \gamn\\ \degv\\ D\\ E\end{array} &
\begin{array}{cc} \car & \car \\ \car & \car \\ \gamn & \car \\ F&\degh\end{array}\\
\hline\end{array}
\eeas
In the same way, one can verify that

\beas
&& \ {}^h\psi_2^- \ {}^h\psi_{1}^-\ {}^h\psi_2^- \begin{array}{|c|c|}
\hline
A & B\\
\hline
& C\\
\hline\end{array} \begin{array}{c}G \\\Longrightarrow\end{array}
\begin{array}{|c|c|}
\hline
D & \\
\hline
E & F\\
\hline\end{array}\\
&&=\ {}^h\psi_1^- \ {}^h\psi_{2}^-\ {}^h\psi_1^- \begin{array}{|c|c|}
\hline
A & B\\
\hline
& C\\
\hline\end{array} \begin{array}{c}G \\\Longrightarrow\end{array}
\begin{array}{|c|c|}
\hline
D & \\
\hline
E & F\\
\hline\end{array}\\
&&=\begin{array}{|c|c|}
\hline
\begin{array}{cc} A&\gamn \end{array} & \begin{array}{ccc} \degv &\car &\car\end{array}\\
\hline
& \begin{array}{ccc} \degv &\car &\car\\B&\degh & \gamn \\C& \gamn & \degv \end{array}\\
\hline\end{array} \begin{array}{c}G \\\Longrightarrow\end{array}
\begin{array}{|c|c|}
\hline
\begin{array}{ccc} D & \gamn & \car\end{array} & \\
\hline
\begin{array}{ccc} E & F & \gamn\end{array} & \car\\
\hline\end{array}
\eeas

\epf

\bth\label{composition} For any $n\geqslant 2$, $\Phi_n^-$ is a composition of
$\ {}^v\psi_i^-$,
$\ {}^h\psi_i^-$ and $\theta_i^-$. \eth

\bpf It is easy to see that
$\Phi_2^-= \ {}^v\psi_1^- \ {}^h\psi_1^- = \ {}^h\psi_1^- \ {}^v\psi_1^-$.
 Now we suppose that $n\geqslant 3$. Set $\Theta_k^{n-2}= \theta_k^- \dots \theta_{n-2}^-$.

We are going to prove that
$$\Phi_n^- = \Theta_{n-2}^{n-2}\Theta_{n-3}^{n-2}\dots \Theta_1^{n-2}
({}^v\Psi_1^- \ {}^h\Psi_1^-)\dots ({}^v\Psi_{n-1}^- \ {}^h\Psi_{n-1}^-)$$
by verifying that the second member satisfies the characterization of
Theorem~\ref{caracterisation}. Let $x\in \omega Cat(I^n,\C)^-$.
Theorem~\ref{aspiration_face_plus}
implies that
for $1\leqslant i \leqslant n$, the dimension of
$$\de_i^+ ({}^v\Psi_1^- \ {}^h\Psi_1^-)\dots ({}^v\Psi_{n-1}^- \ {}^h\Psi_{n-1}^-) x$$
is zero (or equivalently that it belongs
to the image of $\epsilon_1^{n-1}$). With Proposition~\ref{commutation}, one gets
$$
\de_i^+\Theta_{n-2}^{n-2}\Theta_{n-3}^{n-2}\dots \Theta_1^{n-2} ({}^v\Psi_1^- \ {}^h\Psi_1^-)\dots ({}^v\Psi_{n-1}^- \ {}^h\Psi_{n-1}^-) x \in Im(\epsilon_1^{n-1})$$
for $1\leqslant i \leqslant n$. It remains to prove that for $1\leqslant k\leqslant n-2$,
$$
\de_k^-\Theta_{n-2}^{n-2}\Theta_{n-1}^{n-2}\dots \Theta_1^{n-2} y\in Im(\Gamma_{n-2}^-\dots \Gamma_{k}^-)$$
for any $y\in \omega Cat(I^n,\C)^-$. One has
\beas
&& \de_k^-\Theta_{n-2}^{n-2}\Theta_{n-3}^{n-2}\dots \Theta_1^{n-2} y\\
&&= \left(\Theta_{n-3}^{n-3} \dots \Theta_{k}^{n-3}\right) \de_k^- \Theta_{k}^{n-2} z \hbox{ with $z=\Theta_{k-1}^{n-2}\dots \Theta_1^{n-2} y$}\\
&&= \left(\Theta_{n-3}^{n-3} \dots \Theta_{k}^{n-3}\right) \Gamma_k^- \de_k^- \de_k^-\Theta_{k+1}^{n-2} z\\
&&= \Gamma_{n-2}^-\left(\Theta_{n-4}^{n-4} \dots \Theta_{k}^{n-4}\right) \de_k^- \Theta_{k}^{n-3} \de_k^- z \\
&&= (\dots)\\
&&= \left(\Gamma_{n-2}^- \dots \Gamma_{k+1}^-\right)\de_k^- \theta_k^- (\de_k^-)^{n-2-k} z \\
&&=\left(\Gamma_{n-2}^- \dots \Gamma_{k}^-\right)(\de_k^-)^{n-k} z
\eeas

\epf

The operators $\ {}^v\psi_i^\pm$, $\ {}^h\psi_i^\pm$ and $\theta_i^-$ for
$1\leqslant i\leqslant n-1$ and $\Phi_n^-$ induce natural
transformations of set-valued functors from $\omega Cat(I^n,-)^-$ to
itself.

\subsection{Conjecture}\label{cubi1} 
\textit{Let $f$ be an $\omega$-functor from $I^n$ to
  itself such that $f(0_n)=0_n$ and such that the corresponding
  natural transformation from $\omega Cat(I^n,-)$ to itself induces a
  natural transformation $\Phi^-$ from $\omega Cat(I^n,-)^-$ to
  itself. Then $\Phi^-$ is a composition of $\ {}^v\psi_i^-$,
$\ {}^h\psi_i^-$ and $\theta_i^-$ for $1\leqslant i\leqslant n-1$.}

\subsection{Conjecture}\label{cubi2} 
\textit{Let $\Phi$ be a natural transformation from
  $\omega Cat(I^n,-)$ to itself such that the corresponding functor
  $(\Phi)^*$ from $I^n$ to itself satisfies
  $(\Phi)^*(0_n)=0_n$. Then $\Phi$ is a composition of $\ {}^v\psi_i^\pm$ and 
$\ {}^h\psi_i^\pm$ for $1\leqslant i\leqslant  n-1$.}

By Yoneda, the operators $\ {}^v\psi_i^\pm$ and $\ {}^h\psi_i^\pm$ for
$1\leqslant i\leqslant n-1$ induce $\omega$-functors from $I^n$ to
itself denoted by $({}^v\psi_i^\pm)^*$ and $({}^h\psi_i^\pm)^*$.  The
dual conjecture is then

\subsection{Conjecture} \textit{
Let $f$ be an $\omega$-functor from $I^n$ to itself such that
  $f(0_n)=0_n$.  Then $f$ is a composition of $({}^v\psi_i^\pm)^*$ and
  $({}^h\psi_i^\pm)^*$.}

\section{Comparison of $x$ and $\Phi_n^-(x)$ in the reduced branching complex}
\label{id_phi}

This section is devoted to proving that for any $x\in\omega
Cat(I^n,\C)^-$, $x$ and $\Phi_n^-(x)$ are T-equivalent.

\bp For any $i\geqslant 1$ and any $n\geqslant 2$, if $x\in\omega Cat(I^n,\C)^-$, then
$\ {}^h\psi_i^-(x)$ and $x$ are T-equivalent.
\ep

\bpf First let us make the proof for $i=1$ and $n=2$. Let us consider
the following $\omega$-functor from $I^3$ to $\C$ :

\[y_1 =
\begin{array}{|c|c|}
\hline
\begin{array}{cc} x & \gamn \end{array} & \car\\
\hline
& \begin{array}{c}\degh\\ \gamn \end{array}\\
\hline\end{array} \begin{array}{c}x(00) \\\Longrightarrow\end{array}
\begin{array}{|c|c|}
\hline
\gamn  & \\
\hline
x & \gamn\\
\hline\end{array}
\]

Then $\de^- y_1= \ {}^h\psi_1^-(x) -x + t_1$ where $t_1$ is a thin element. Therefore
$x$ and $\ {}^h\psi_1^-(x)$ are T-equivalent.

We claim that the above construction is sufficient to prove that
$x$ and $\ {}^h\psi_1^-(x)$ are T-equivalent for any $x\in\omega
Cat(I^n,\C)^-$ and for any $n\geqslant 2$. The labeled $3$-cube
$y_1$ is actually a certain thin $3$-dimensional element of the
cubical $\omega$-category $\mathcal{N}^\square(\C)$ and it
corresponds to the filling of a thin $2$-shell. So
$$y_1=f_1(\epsilon_1 x, \epsilon_2 x, \epsilon_3 x, \Gamma_1^- x, \Gamma_2^- x,
\Gamma_1^+ x, \Gamma_2^+ x)$$
where $f_1$ is a function which only uses  the operators $+_1$, $+_2$,
and $+_3$. In this particular case, $f_1$ could be of course
calculated. But it will not be always possible in the sequel to make
such a calculation : this is the reason why no explicit formula is used
here. And one has $\de_2^- f_1(x)= x$, $\de_3^- f_1(x)=
\ {}^h\psi_1^-(x)$ and all other $2$-faces $\de_i^\alpha f_1(x)$ are
(necessarily) thin $2$-faces. The equalities $\de_2^- f_1(x)= x$ and
$\de_3^- f_1(x)= \ {}^h\psi_1^-(x)$ do not depend on the dimension of
$x$. Therefore for any $x\in\omega Cat(I^n,\C)^-$ and for any $n\geq
2$, one gets $\de^- y_1=\ {}^h\psi_1^-(x) -x + t$ where $t$ is a linear
combination of thin elements.

Now we want to explain that the above construction is also
sufficient to prove that $x$ and $\ {}^h\psi_i^-(x)$ are
T-equivalent for any $i\geq1$ and any $x\in\omega Cat(I^n,\C)^-$
and for any $n\geqslant 2$. The equalities $\de_2^- f_1(x)= x$
and $\de_3^- f_1(x)= \ {}^h\psi_1^-(x)$ do not depend on the
absolute values $1$, $2$ $3$. But only on the relative values
$1=3-2$, $2=3-1$ and $3=3-0$. So let us introduce a labeled
$(n+1)$-cube $y_i=f_i(x)$ by replacing in $f_1$ any index $1$ in
by $i$, any index $2$ by $i+1$ and any index $3$
by $i+2$. Then one gets a thin $(n+1)$-cube
$y_i=f_i(x)$ such that $\de_{i+1}^- f_i(x)= x$ and $\de_{i+2}^-
f_i(x)= \ {}^h\psi_i^-(x)$.

If the reader does not like this proof and prefers explicit calculations, it
suffices to notice that $y_1= \ {}^h\psi_1^- \Gamma_2^- x$ by
Proposition~\ref{formule_3D}. Set $y_i=
\ {}^h\psi_i^- \Gamma_{i+1}^- x$. Then
\beas &&\de^-
(y_i)=\sum_{j<i}(-1)^{j+1}\ {}^h\psi_{i-1}^- \Gamma_{i}^-\de_j^- x +
(-1)^{i+1}(\Gamma_{i}^-\de_i^- x+_i \epsilon_{i+1}\de_{i+1}^+ x) +\\&&
(-1)^{i+2} (x-\ {}^h\psi_i^- x) + \sum_{j>i+2} (-1)^{j+1}\ {}^h\psi_i^-
 \Gamma_{i+1}^- \de_{j-1}^- x \eeas and one completes the proof by an
 easy induction on the dimension of $x(0_n)$.
\epf

\bp For any $i\geqslant 1$ and any $n\geqslant 2$, if $x\in\omega Cat(I^n,\C)^-$, then
$\ {}^v\psi_i^-(x)$ and $x$ are T-equivalent.
\ep

\bpf It suffices to make the proof for $i=1$ and $n=2$. And to consider
the following thin $3$-cube
\[
\begin{array}{|c|c|}
\hline
\begin{array}{cc}\gamn & \degv\end{array} & \car\\
\hline
& \begin{array}{c} \gamn\\ x\end{array}\\
\hline\end{array} \begin{array}{c}x(00) \\\Longrightarrow\end{array}
\begin{array}{|c|c|}
\hline
\gamn  & \\
\hline
x & \gamn\\
\hline\end{array}
\]
Notice that the above $3$-cube is exactly $\ {}^v\psi_2^- \Gamma_1^- x$ by
Proposition~\ref{formule_3D}.
\epf

\bp For any $i\geqslant 1$ and any $n\geqslant 3$, if $x\in\omega Cat(I^n,\C)^-$, then
$\theta_i^-(x)$ and $x$ are T-equivalent. \ep

\bpf It suffices to make the proof for $i=1$ and $n=3$. Set
\[x =
\begin{array}{|c|c|}
\hline
A & B\\
\hline
& C\\
\hline\end{array} \begin{array}{c}G \\\Longrightarrow\end{array}
\begin{array}{|c|c|}
\hline
D & \\
\hline
E & F\\
\hline\end{array}
\]
One has already seen that
\[\theta_1^- x =\begin{array}{|c|c|}
\hline
\begin{array}{cc} A & B\\ \gamp & C\end{array} & \begin{array}{c}\degh
\\ \degh \end{array}\\
\hline
& \begin{array}{c}\gamn
\\ \degv \end{array}\\
\hline\end{array} \begin{array}{c}G \\\Longrightarrow\end{array}
\begin{array}{|c|c|}
\hline
\begin{array}{c}\gamn
\\ D \end{array} & \\
\hline
\begin{array}{c}E
\\ \degv \end{array} & \begin{array}{cc}F & \degh \\ \gamp & \degh\end{array}\\
\hline\end{array}
\]

\begin{figure}
\begin{center}
\xymatrix{
&  {\begin{array}{|c|c|c|}\hline{++00}&&\\\hline{+0-0}&{+00+}&\\\hline&{0-0+}&{00++}\\\hline\end{array}} & 
\\
{\begin{array}{|c|c|c|}\hline{+0-0}&&\\\hline{0--0}&{00-+}&{0+0+}\\\hline&&{-00+}\\\hline\end{array}} \ar@3{->}[ru]^-{000+} && 
{\begin{array}{|c|c|c|}\hline{+00-}&{+0+0}&\\\hline&{+-00}&\\\hline&{0--0}&{0-0+}\\\hline\end{array}} \ar@3{->}[lu]_-{+000}
\\
{\begin{array}{|c|c|c|}\hline&{++00}&\\\hline{00--}&{0+-0}&{0+0+}\\\hline&{-0-0}&\\\hline\end{array}} \ar@3{->}[u]^-{00-0}
&&{\begin{array}{|c|c|c|}\hline&{+0+0}&\\\hline{0-0-}&{0-+0}&{00++}\\\hline&{--00}&\\\hline\end{array}}
\ar@3{->}[u]_-{0-00}
\\
{\begin{array}{|c|c|c|}\hline {0+0-}&{0++0}&\\\hline&{-+00}&\\\hline &{-0-0}&{-00+}\\\hline\end{array}}
\ar@3{->}[u]^-{0+00}&&{\begin{array}{|c|c|c|}\hline {+00-}&&\\\hline {0-0-}&{00+-}&{0++0}\\\hline&&{-0+0}\\\hline\end{array}}
\ar@3{->}[u]_-{00+0}
\\
&{\begin{array}{|c|c|c|}\hline {00--}&{0+0-}&\\\hline&{-00-}&{-0+0}\\\hline&&{--00}\\\hline\end{array}}
\ar@3{->}[ul]^-{-000}\ar@3{->}[ur]_-{000-}&
}
\end{center}
\caption{2 -categorical representation of the  $4$ -cube}
\label{4cube_symb}
\end{figure}
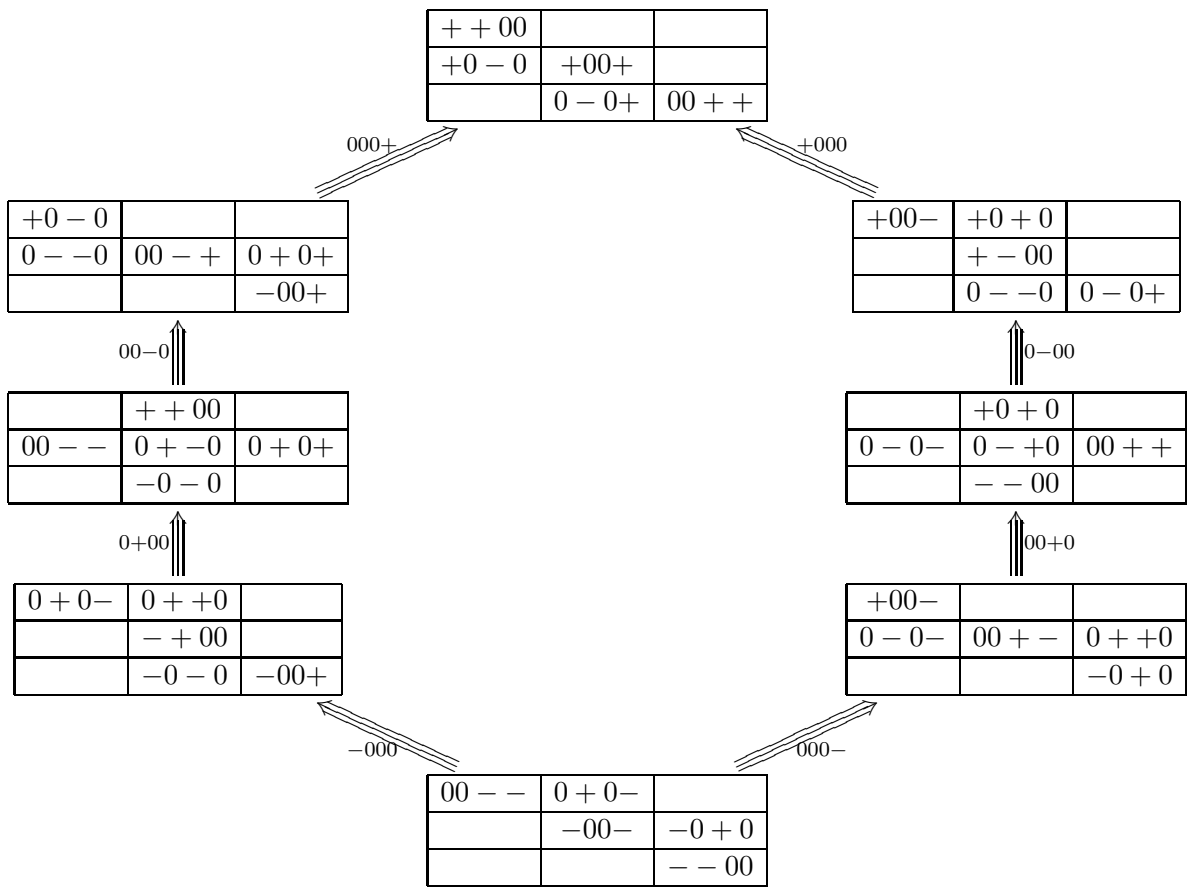

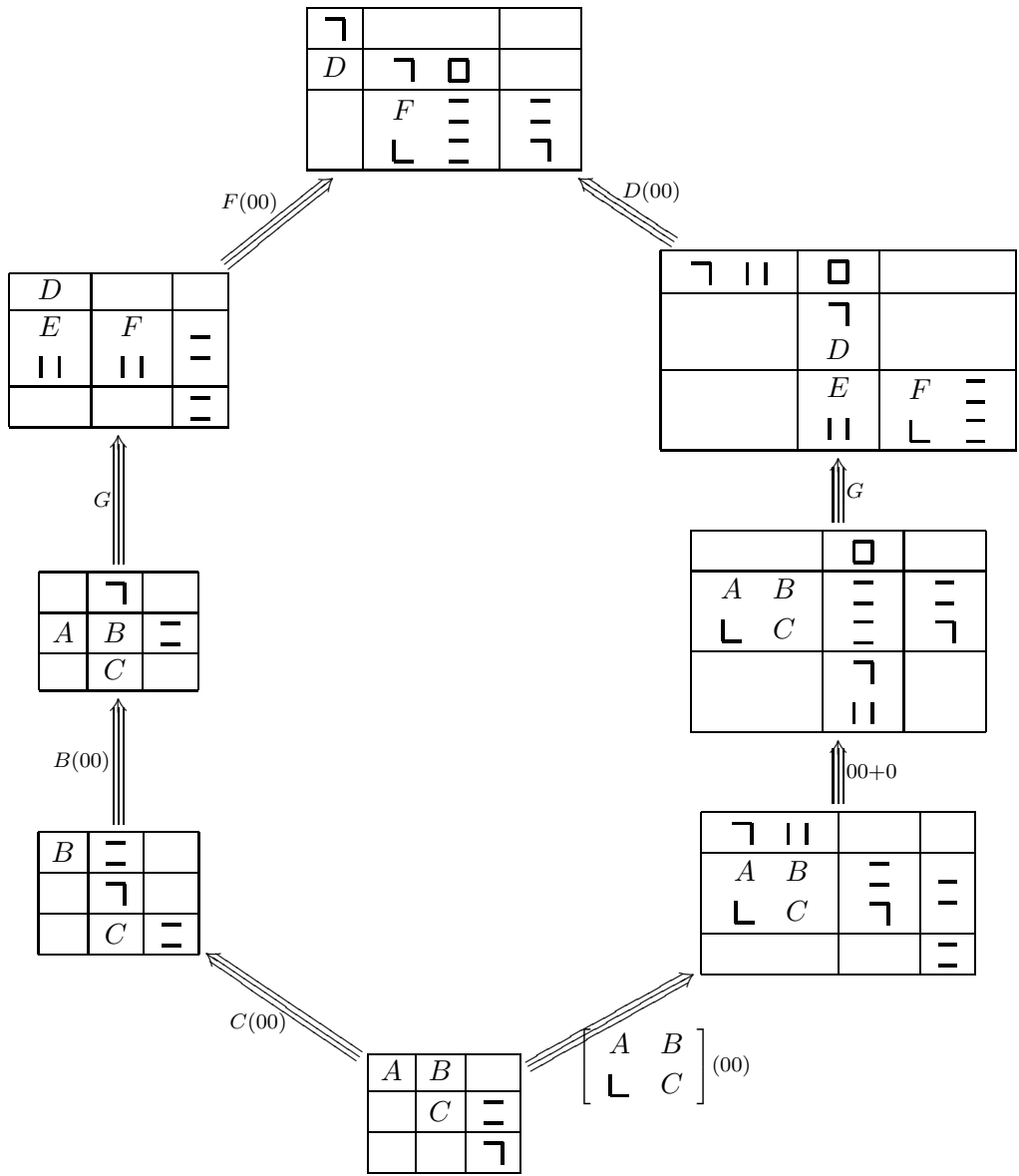
\begin{figure}
\begin{center}
\[
\xymatrix{
& {\begin{array}{|c|c|c|}\hline{\gamn}&&\\\hline{D}&{\begin{array}{cc}\gamn&\car\end{array}}&\\\hline&{\begin{array}{cc}F & \degh \\ \gamp & \degh\end{array}}&{\begin{array}{c}\degh\\\gamn\end{array}}\\\hline\end{array}} & 
\\
{\begin{array}{|c|c|c|}\hline{D}&&\\\hline{\begin{array}{c}E
\\ \degv \end{array}}&{\begin{array}{c}F\\\degv\end{array}}&{\degh}\\\hline&&{\degh}\\\hline\end{array}}
\ar@3{->}[ur]^-{F(00)} & & {\begin{array}{|c|c|c|}\hline{\begin{array}{cc}\gamn &\degv \end{array}}&{\car}&\\\hline&{\begin{array}{c}\gamn
\\ D \end{array}}&\\\hline&{\begin{array}{c}E
\\ \degv \end{array}}&{\begin{array}{cc}F & \degh \\ \gamp & \degh\end{array}}\\\hline\end{array}}
\ar@3{->}[ul]_-{D(00)}
\\
{\begin{array}{|c|c|c|}\hline&{\gamn}&\\\hline{A}&{B}&{\degh}\\\hline&{C}&\\\hline\end{array}}
\ar@3{->}[u]^-{G} && {\begin{array}{|c|c|c|}\hline&{\car}&\\\hline{\begin{array}{cc} A & B\\ \gamp & C\end{array}}&{\begin{array}{c}\degh
\\ \degh \end{array}}&{\begin{array}{c}\degh\\\gamn\end{array}}\\\hline&{\begin{array}{c}\gamn
\\ \degv \end{array}}&\\\hline\end{array}}
\ar@3{->}[u]_-{G} 
\\
{\begin{array}{|c|c|c|}\hline {B}&{\degh}&\\\hline&{\gamn}&\\\hline &{C}&{\degh}\\\hline\end{array}}
\ar@3{->}[u]^-{B(00)} && {\begin{array}{|c|c|c|}\hline {\begin{array}{cc}\gamn &\degv \end{array}}&&\\\hline {\begin{array}{cc} A & B\\ \gamp & C\end{array}}&{\begin{array}{c}\degh\\\gamn\end{array}}&{\degh}\\\hline&&{\degh}\\\hline\end{array}}
\ar@3{->}[u]_-{00+0}\\
& {\begin{array}{|c|c|c|}\hline {A}&{B}&\\\hline&{C}&{\degh}\\\hline&&{\gamn
}\\\hline\end{array}}
\ar@3{->}[ul]^-{C(00)}\ar@3{->}[ur]_-{\left[\begin{array}{cc}A&B\\\gamp &C\end{array}\right](00)} &
}\]
\end{center}
\caption{A labeled $4$-cube} \label{thin4}
\end{figure}

It suffices to construct a thin $4$-cube $y$ such that $\de_3^-
y=x$ and $\de_2^- y=\theta_1^- x$. If the $4$-cube is
conventionally represented by Figure~\ref{4cube_symb}, the thin
<<<<<<< coin.tex
labeled $4$-cube of Figure~\ref{thin4} with 
$00\!+\!0\mapsto (\de_2^+\de_1^- x +_1 \de_2^+ \de_2^+ x)(0)$ meets the requirement.
=======
labeled $4$-cube of Figure~\ref{thin4} with 
$00\!+\!0\mapsto (\de_2^+\de_1^- x +_1 \de_2^+ \de_2^+ x)(0)$ 
meets the requirement.
>>>>>>> 1.19
The latter labeled $4$-cube can be defined as the unique thin
$4$-cube $\omega(x)$ which fills the $3$-shell defined by \beas
&&\de_1^-\omega(x)=\Gamma_2^- \de_1^- x\\
&& \de_2^-\omega(x)=\theta_1^-(x)\\
&&\de_3^-\omega(x)=x\\
&&\de_4^-\omega(x)=\left[\begin{array}{cc} \Gamma_2^- \de_3^- x & \epsilon_2
\de_2^+ x\\ \Gamma_1^- \Gamma_1^+ \de_1^- \de_2^+ x +_2 \epsilon_1 \Gamma_1^- \de_1^-
\de_2^+ x & \Gamma_1^- \de_1^- x\end{array}\right]\coin{3}{1}\\
&&\de_1^+ \omega(x)= \ {}^v\psi_2^- \Gamma_1^- \de_1^+ x\\
&&\de_2^+ \omega(x)= \Gamma_2^- \de_2^+ x\\
&& \de_3^+ \omega(x)= \epsilon_3(\Gamma_1^- \de_2^+ \de_1^- x +_1 \epsilon_2\de_2^+ \de_2^+ x)\\
&&\de_4^+ \omega(x)=\ {}^v\psi_2^+ \Gamma_2^- \de_3^+ x
\eeas

\epf

\bcor\label{identification} For any $n\geqslant 2$, for any
$x\in\omega Cat(I^n,\C)^-$, $x$ and
$\Phi_n^-(x)$ are T-equivalent and $\Phi_n^-$ is the identity map on
the reduced branching complex.
\ecor

We have proved that for any $x\in \omega Cat(I^n,\C)^-$, there exists
$t_1\in M_n$ and $t_2\in M_{n+1}$ such that $\Phi_n^-(x)-x=t_1 +\de^-
t_2$.  The proofs of this section use only calculations in the free
cubical $\omega$-category generated by $x$. This means that $t_1$ and
$t_2$ can be formulated in terms of expressions in the same cubical
$\omega$-category. And so this means that $t_1$ and $t_2$ are linear
combinations of expressions which use only $x$ as variable and the
operators $\de_i^\pm $, $\Gamma_i^\pm $, $\epsilon_i$ and $+_i$. With
Theorem~\ref{composition} which allows to consider $\Phi_n^-$ like an
operator defined in any cubical $\omega$-category, one sees that
Corollary~\ref{identification} does make sense in an appropriate
cubical setting. Moreover the terms $t_1$ and $t_2$ being elements 
of the free cubical $\omega$-category generated by $x$, then $t_1$ and 
$t_2$ depend in a functorial way on $x$.

\section{Folding operations and composition maps}\label{folding_composition}

\bth\label{composition_glob_folding} Let $x$ and $y$ be two $n$-morphisms of
$\C$ with $n\geqslant 2$.
\begin{enumerate}
\item if $x *_{n-1} y$ exists, then $\square_n^-(x *_{n-1}y)-\square_n^-(x)
-\square_n^-(y)$ is a boundary in the normalized chain complex of the
branching simplicial nerve of $\C$. Moreover, $\square_n^-(x *_{n-1}y)$ is
T-equivalent to $\square_n^-(x)+\square_n^-(y)$.
\item if $1\leqslant p \leqslant n-2$, then
$\square_n^-(x *_{p}y)$ is T-equivalent to $\square_n^-(x)+\square_n^-(y)$.
\end{enumerate}
\eth

\bpf Let us denote by $P(h)$ the following property :

``for any $n\geqslant 2$ and with $p=n-h\geqslant 1$, for any $n$-morphisms $x$
and $y$ of any $\omega$-category $\C$ such that $x *_p y$ exists,
there exists a thin $n$-cube $A^n_p(x,y)$ and a thin $(n+1)$-cube
$B^n_{p}(x,y)$ which lie in the cubical singular nerve of the free
globular $\omega$-category
generated by $x$ and $y$, and even in its  branching  nerve, such
that
$$\square_n^-(x *_{p}y)= \square_n^-(x)+\square_n^-(y)+
A^n_p(x,y)+\de^- B^n_{p}(x,y)$$ in the normalized branching
complex (i.e. the equality holds modulo degenerate elements of the
branching simplicial nerve) and such that for any $(n+1)$-morphisms $u$
and $v$, $A^n_p(s_n u,s_n v)=A^n_p(t_n u,t_n v)$.''

Since \beas
&&\de^-\left(\square_n^-(x *_{n-1}y)-\square_n^-(x)-\square_n^-(y)\right)\\
&&=\square_{n-1}^-(s_{n-1}x)-\square_{n-1}^-(t_{n-1}y)\\&&-\square_{n-1}^-(s_{n-1}x)
+ \square_{n-1}^-(t_{n-1}x)-\square_{n-1}^-(s_{n-1}y)
+ \square_{n-1}^-(t_{n-1}y)\\
&&=\square_{n-1}^-(t_{n-1}x)-\square_{n-1}^-(s_{n-1}y)=0 \eeas in the
normalized chain complex of the branching simplicial nerve, then
$\square_n^-(x *_{n-1}y)-\square_n^-(x)-\square_n^-(y)$ is a cycle in
the branching homology of the free globular $\omega$-category $\D$ generated by
two $n$-morphisms such that $t_{n-1}x=s_{n-1}y$. The $\omega$-category
$\D$ is of length at most one and non-contracting. Therefore its
branching nerve coincides with the simplicial nerve of $\P\D$, the
latter being the globular $\omega$-category freely generated by the
composable pasting scheme whose total composition is $X *_{n-2} Y$
where $X$ and $Y$ are two $(n-1)$-dimensional cells. Therefore this
simplicial nerve is contractible. Consequently there exists
$B^n_{n-1}(x,y)$ lying in the cubical singular nerve of $\D$ (and also
in its branching nerve) such that 
$$\square_n^-(x *_{n-1}y)-\square_n^-(x)-\square_n^-(y)=\de^- B^n_{n-1}(x,y).$$
The
$(n+1)$-cube $B^n_{n-1}(x,y)$ is necessarily thin because there is no
morphism of dimension $n+1$ in $\D$. By setting $A^n_{n-1}(x,y)=0$, we
obtain $P(1)$. We are going to prove $P(h)$ by induction on $h$.
Suppose $P(h)$ proved for $h\geqslant 1$. Then \beas &&\de^-(
\square_n^-(x *_{n-h-1} y)-\square_n^-(x)-\square_n^-(y)-
B_{n-h-1}^{n-1}(s_{n-1}x , s_{n-1}y) \\&&+ B_{n-h-1}^{n-1}(t_{n-1}x , t_{n-1}y))\\
&=&\left(\square_n^-(s_{n-1} x *_{n-h-1} s_{n-1}
  y)-\square_n^-(s_{n-1}x)-\square_n^-(s_{n-1}y)
  -\de^-B_{n-h-1}^{n-1}(s_{n-1}x,
  s_{n-1}y)\right)\\&&-\left(\square_n^-(t_{n-1} x *_{n-h-1} t_{n-1}
  y)-\square_n^-(t_{n-1}x)-\square_n^-(t_{n-1}y)
  -\de^-B_{n-h-1}^{n-1}(t_{n-1}x , t_{n-1}y)\right)\\
&=&A_{n-h-1}^{n-1}(s_{n-1}x, s_{n-1}y)-A_{n-h-1}^{n-1}(t_{n-1}x, t_{n-1}y)\hbox{ by induction hypothesis}\\
&=&0\hbox{ again by induction hypothesis} \eeas Therefore we can set
$A_{n-h-1}^{n}(x,y)=B_{n-h-1}^{n-1}(s_{n-1}x , s_{n-1}y) -
B_{n-h-1}^{n-1}(t_{n-1}x , t_{n-1}y)$ and we have \beas
&&A_{n-h-1}^{n}(s_{n}u,s_{n}v)-A_{n-h-1}^{n}(t_{n}u,t_{n}v)\\
&=&B_{n-h-1}^{n-1}(s_{n-1}s_{n}u , s_{n-1}s_{n}v) - B_{n-h-1}^{n-1}(t_{n-1}s_{n}u , t_{n-1}s_{n}v)  \\&&-  B_{n-h-1}^{n-1}(s_{n-1}t_{n}u , s_{n-1}t_{n}v) + B_{n-h-1}^{n-1}(t_{n-1}t_{n}u , t_{n-1}t_{n}v)\\
&=&0 \eeas because of the globular equations. So we get a thin
$n$-cube $A_{n-h-1}^{n}(x,y)$ such that $$\square_n^-(x *_{n-h-1}
y)-\square_n^-(x)-\square_n^-(y)-A_{n-h-1}^{n}(x,y)$$
is a cycle in
the normalized chain complex associated to the branching simplicial
nerve of $\C$. This cycle lies in the branching nerve of the free
$\omega$-category generated by two $n$-morphisms $x$ and $y$ such that
$t_{n-h-1}x=s_{n-h-1}y$. This $\omega$-category is of length at most
one and non-contracting.  Therefore its branching nerve is isomorphic to the simplicial
nerve of the globular $\omega$-category freely generated by the
composable pasting scheme whose total composition is $X *_{n-h-2} Y$
where $X$ and $Y$ are two $(n-1)$-dimensional cells. Therefore it is
contractible.  Therefore there exists $B_{n-h-1}^{n}(x,y)$ such that
$$\square_n^-(x *_{n-h-1}
y)-\square_n^-(x)-\square_n^-(y)-A_{n-h-1}^{n}(x,y)=\de^-
B_{n-h-1}^{n}(x,y).$$ The cube $B_{n-h-1}^{n}(x,y)$ is necessarily
thin because there is no morphism of dimension $n+1$ in the cubical
sub-$\omega$-category generated by $x$ and $y$. And $P(h+1)$ is
proved.
\epf

It turns out that the $(n+1)$-cube $B^n_{n-1}(x,y)$ can be
explicitly calculated. One can easily verify that
\[B^n_{n-1}(x,y)_h^-=\Gamma_{n-1}^- \Gamma_{n-2}^-\dots \Gamma_{h}^- \square_h^- d_h^{(-)^h} x\]
for $1\leq h \leqslant n-2$ (observe that in this case, $d_h^{(-)^h} x=d_h^{(-)^h} y$),
\beas
&& B^n_{n-1}(x,y)_{n-1}^-=\square_n^- y \\
&& B^n_{n-1}(x,y)_n^-=\square_n^- (x *_{n-1} y) \\
&& B^n_{n-1}(x,y)_{n+1}^-=\square_n^- x
\eeas
and for all $i$ between
$1$ and $n+1$, 
\[B^n_{n-1}(x,y)_i^+=\square_n^- t_0 x\] is a
solution.  It suffices to prove that 
$(B^n_{n-1}(x,y)_i^\pm)_{1\leqslant i\leqslant n+1}$ is a thin $n$-shell.

\bth\label{folding_plus} Let $\C$ be a non-contracting $\omega$-category. Let $x$ and $y$
be two elements of
$\omega Cat(I^n,\C)$ such that $x+_j y$ exists for some $j$ between
$1$ and $n$ and such that $dim(x(0_n))\geq 1$, $dim(y(0_n))\geq 1$ and 
$dim((x+_j y)(0_n))\geq 1$. Then $\Phi_n^-(x+_j y)$ is T-equivalent to $\Phi_n^-(x)$ or
$\Phi_n^-(y)$ or to $\Phi_n^-(x)+\Phi_n^-(y)$. If $x$ is itself in the
branching complex, then $\Phi_n^-(x+_j y)$ is T-equivalent to
$x$. \eth

\subsection*{Remark} The hypotheses about the dimension of $x(0_n)$, 
$y(0_n)$ and $(x+_j y)(0_n)$ are  only to ensure that $\Phi_n^-(x)$, 
$\Phi_n^-(y)$ and $\Phi_n^-(x+_j y)$ are in the branching nerve. The hypothesis 
about the dimension of $(x+_j y)(0_n)$ is necessary because we do not assume that 
$1$-morphisms in non-contracting $\omega$-categories are not invertible. In dimension 
$1$, the case $x(0)*_0 y(0)=(x +_1 y)(0)\in \C_0$ may happen.

\bpf By definition, one has $\Phi_n^-(x+_j y)=\square_n^-((x+_j y)(0_n))$.
If $\C$ was equal to the globular sub-$\omega$-category generated by
$$X=\{x(k_1\dots k_n),k_1\dots k_n\in \underline{cub}^n\}\cup \{y(k_1\dots
k_n),k_1\dots k_n\in \underline{cub}^n\}$$ then $x+_j y$ still would exist in
the cubical singular nerve. Therefore, $(x+_j y)(0_n)$ can be written
as an expression using only the composition laws $*_n$ of $\C$ and the
variables of $X$ and moreover, the variables $x(0_n)$ and $y(0_n)$ can
appear at most once. By Theorem~\ref{composition_glob_folding},
$\square_n^-((x+_j y)(0_n))$ is therefore T-equivalent to
$\square_n^-(x(0_n))$, $\square_n^-(y(0_n))$ or
$\square_n^-(x(0_n))+\square_n^-(y(0_n))$.

Now suppose that $x\in\omega Cat(I^n,\C)^-$.  Let $z=\Gamma_j^- x +_j
\epsilon_{j+1} y \in \omega Cat(I^{n+1},\C)^-$.  Then $\de_j^- z=x$,
$\de_{j+1}^- z = x+_j y$ and $\de_{j+1}^+ z=y$. Since $z$ is a thin
element, then all other faces $\de_k^\pm z$ are thin (this can be
verified directly by easy calculations). Therefore $\de^-z$ is
T-equivalent to $\pm (x+_jy -x)$. As illustration, let us notice that
for $j=1$ and $n=2$, $z$ is equal to
\[
\begin{array}{|c|c|}
\hline
\begin{array}{c}\degh\\\gamn \end{array}
& \begin{array}{c}y\\\degv\end{array}\\
\hline
& x\\
\hline\end{array} \begin{array}{c}(x+_1y)(00) \\\Longrightarrow\end{array}
\begin{array}{|c|c|}
\hline
\degv & \\
\hline
\begin{array}{c}y\\x\end{array} & \begin{array}{c}\degh\\\gamn\end{array}\\
\hline\end{array}
\]
\epf

\bth Let $x$ and $y$ be two morphisms of a non-contracting 
$\omega$-category $\C$
such that $x*_0y$ exists such that $x$ and $x*_0 y$ are  of 
dimension lower than $n$ and of dimension strictly greater than $0$. Then
$\square_n^-(x*_0y)$ is T-equivalent to $\square_n^-(x)$.
\eth

\bpf We need, only for this proof, the operator $\square_n^+$ introduced
in \cite{Gau}. One has
\[\de_1^+\square_n^-(x)=\epsilon_1^{n-1}\square_0(t_0x)\] and
\[\de_1^-\square_n^+(y)=\epsilon_1^{n-1}\square_0(s_0y).\]Therefore
$\square_n^-(x)+_1 \square_n^+(y)$ exists and is T-equivalent to
$\square_n^-(x)$ by Theorem~\ref{folding_plus}. If we work in the
$\omega$-category generated by $x$ and $y$, then $\square_n^-(x)+_1
\square_n^+(y)$ is a well-defined element of the branching simplicial nerve
of $\D$. And $\D$ is the free $\omega$-category generated by a
composable pasting scheme whose total composition is $x *_0y$. 
Since union means composition in such a
$\omega$-category, then necessarily $\left(\square_n^-(x)+_1
\square_n^+(y)\right)(0_n)=x*_0 y$. Since $\Phi_n^-$ is the identity
map on the reduced branching complex, then $\square_n^-(x)+_1
\square_n^+(y)$ is T-equivalent to $\square_n^-(x*_0y)$. \epf

The preceding formulae suggest another way of defining the
reduced branching homology.

\bp\label{coin_formel} Set $$\CF^-_n(\C)= \Z\C_n/\{x *_0y=x, x*_1y=x+y, \dots , x
*_{n-1}y=x+y\hbox{ mod }\Z\tr_{n-1}\C\}.$$ Then $s_{n-1}-t_{n-1}$ from $\CF^-_n(\C)$ to
$\CF^-_{n-1}(\C)$ for $n\geqslant 2$ and $s_0$ from $\CF^-_1(\C)$
to $\CF^-_0(\C)$ induce a differential map $\de_f^-$ on
the $\N$-graded group $\CF^-_*(\C)$ and the chain complex one gets is
called the formal branching complex. The associated homology is
denoted by $\HF_n^-(\C)$ and is called the formal branching homology.  \ep

\bpf Obvious. \epf

A relation like $x *_0y=x\hbox{ mod }\Z\tr_{n-1}\C$ means that if $x$
is for example a $p$-morphism for $p<n$ and $y$ a $n$-morphism such
that $x*_0 y$ exists, then in $\CF^-_n(\C)$, $x*_0y=0$.

\bp Let $\C$ be a non-contracting $\omega$-category. 
The linear map $\square_n^-$ from $\Z\C_n$ to
$\CR^-_n(\C)$ induces a surjective morphism of chain complexes and
therefore a morphism from $\HF_*^-(\C)$ to $\HR_*^-(\C)$. \ep

\bpf One has in the reduced branching complex $\square_n^-(x *_0y)=\square_n^-(x)$ and
$\square_n^-(x *_py-x-y)=0$ therefore $\square_n^-$ induces a linear
map from $\CF_n^-(\C)$ to $\CR_n^-(\C)$. And
$\square_{n-1}^-(\de_f^-(x))=\square_{n-1}^-(s_{n-1}-t_{n-1})(x)=\de^-\square_n^-(x)$.
Since $\Phi_n^-$ is the identity map on $\CR_n^-(\C)$, then
$\CR_n^-(\C)$ is generated by the $\square_n^-(x)$ where $x$ runs over
$\C_n$. Therefore the induced morphism of chain complexes is surjective.
\epf

\subsection{Question}\label{formelle}
\textit{When is the preceding map a quasi-isomorphism ?}

The meaning of the results of this section is that one homology class
in branching  homology corresponds really to one branching area. Here are
some simple examples to understand this fact.

\begin{figure}
\begin{center}
\includegraphics[width=7cm]{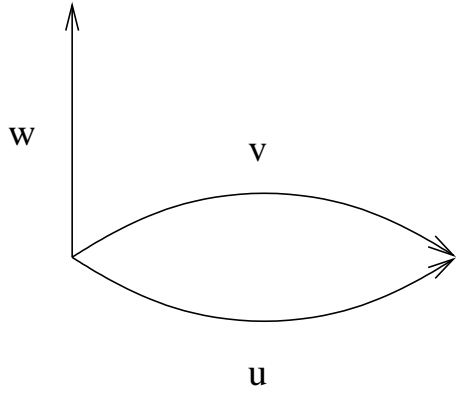}
\end{center}
\caption{Another $1$-dimensional branching area}
\label{ex2}
\end{figure}

Figure~\ref{ex1} represents a $1$-dimensional branching area. This branching
area corresponds to one element in the reduced branching  homology, that is
$$\square_1(u)-\square_1(w)=\square_1(u*_0v)-\square_1(w)=\square_1(u)-\square_1(w*_0x)=
\square_1(u*_0v)-\square_1(w*_0x)$$ in homology. In fact, it even corresponds  to one
cycle in the reduced branching  complex. The reason why it is more
appropriate to work anyway with cycles modulo boundaries, and not only with
cycles modulo boundaries of thin elements is illustrated in
Figure~\ref{ex2}.  The two cycles $\square_1(u)-\square_1(v)$ and
$\square_1(u)-\square_1(w)$ are in the same homology class as soon as
$u$ is homotopic to $v$.

\begin{figure}
\begin{center}
\includegraphics[width=7cm,height=7cm]{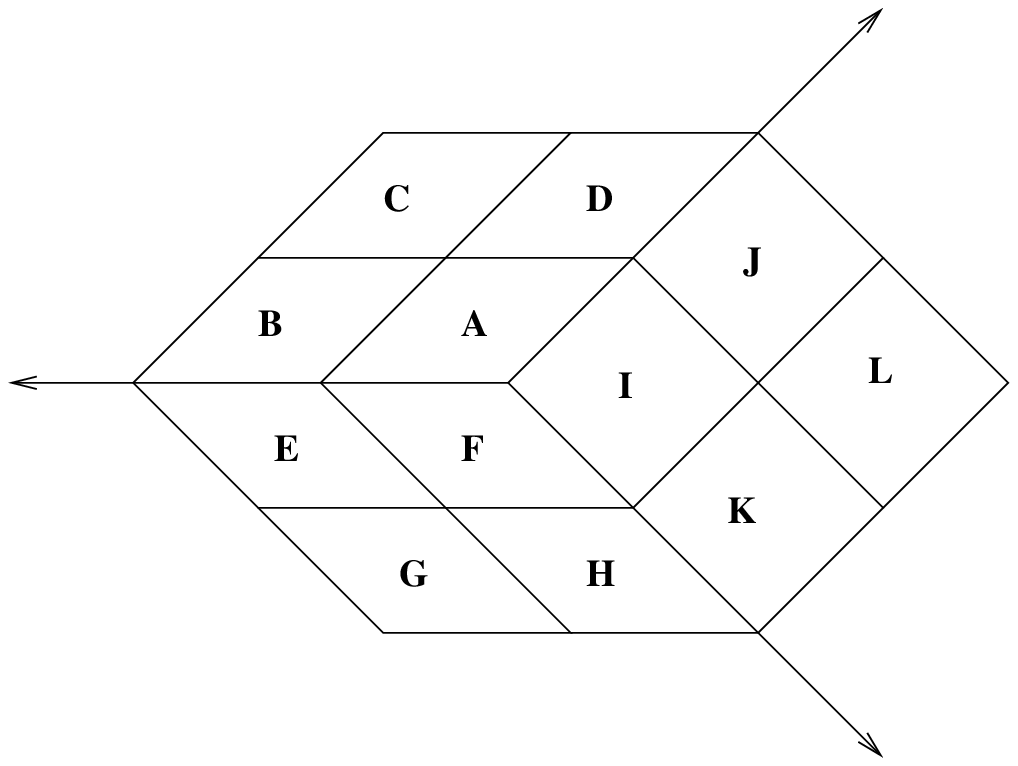}
\end{center}
\caption{A $2$-dimensional branching area}
\label{ex3}
\end{figure}

These observations can be generalized in higher dimension but they are more
difficult to draw. If $u$ is a $n$-morphism, then, by definition of
$\square_n^-$, $\square_n^-(u)$ is an homotopy between
$\square_{n-1}^-s_{n-1}u$ and $\square_{n-1}^-t_{n-1}u$ in the
branching simplicial nerve. Figure~\ref{ex3} is an analogue of
Figure~\ref{ex1} in dimension $2$. Figure~\ref{ex3} represents a
$2$-dimensional branching area. In the branching complex, it corresponds
to the cycles $(A)-(F)+(I)$, $(A,B,C,D)-(E,F,G,H)+(I,J,K,L)$,
$(A)-(F,H)+(I,K)$, etc. In the reduced branching complex, there are even
more possible cycles which correspond to this branching area. For
example $(A,D)-(E,F)+(I,J,K,L)$, $(A)-(E,F)+(I,J)$, etc. In the branching
homology, all these cycles are equivalent and therefore there is
really one homology class which corresponds to one branching area. Or
in other terms, the homology class does not depend on a cubification of
the HDA.

\section{Folding operations and differential map}\label{diff_section}

Now we explore the relations between the folding operators and the
differential map of the branching complex.

\bp Let $x$ be an element of $\omega Cat(I^n,\C)^-$. Then
$$\square_{n-1}^-(s_{n-1}-t_{n-1})(x(0_n))=\square_{n-1}^-(\de^-x)(0_{n-1})=\sum_{p=1}^{n}(\de_p^-x)(0_{n-1})$$
in $\CR_{n-1}^-(\C)$.
\ep

\bpf Since $\Phi_n^-$ induces the identity map on $\CR_n^-(\C)$, then $\Phi_{n-1}^- \de^-=\de^- \Phi_n^-=\de^-$. Therefore
$$\square_{n-1}^-(\de^-x)(0_{n-1})=\Phi_{n-1}^- \de^- x=\de^- \Phi_n^-x=\de^- \square_n^-(x(0_n))=\square_{n-1}^-(s_{n-1}-t_{n-1})(x(0_n)).$$
\epf

\bp In the reduced branching  homology of a given $\omega$-category $\C$,
one has
\begin{enumerate}
\item if $x\in \omega Cat(I^2,\C)^-$, then $\square_1^- (s_1 x(00))=\square_1^- x(-0)$
and $\square_1^- (t_1 x(00))=\square_1^- x(0-)$
\item if $x\in \omega Cat(I^3,\C)^-$, then
\beas
&&\square_2^- (s_2 x(000))=\square_2^- x(-00)+\square_2^- x(-0-)\\
&&\square_2^- (t_2 x(000))=\square_2^- x(0-0).\eeas
\end{enumerate}
\ep

\bpf One has $$\square_1^- (s_1 x(00))=\square_1^- (s_1 (x(-0)*_0
x(0+)))=\square_1^- x(-0)$$ and $$\square_1^- (t_1 x(00))=\square_1^-
(s_1 (x(0-)*_0 x(+0)))=\square_1^- x(0-).$$ Now suppose that
$x\in\omega Cat(I^3,\C)^-$. Then
\beas
&&\square_2^- (s_2 x(000)) \\
&&= \square_2^- \left(
(x(-00)*_0 x(0++)) *_1 (x(-0-)*_0 x(0+0)) *_1 (x(00-)*_0 x(++0))
\right))\\
&&= \square_2^-(x(-00)*_0 x(0++)) + \square_2^-(x(-0-)*_0 x(0+0))+ \square_2^-(x(00-)*_0 x(++0))
\eeas

So $\square_2^- (s_2 x(000))=\square_2^-(x(-00))+\square_2^-(x(00-))$.

In the same way, one has
\beas
&&\square_2^- (t_2 x(000)) \\
&&=\square_2^- \left(
(x(--0)*_0 x(00+)) *_1 (x(0-0)*_0 x(+0+)) *_1 (x(0--)*_0 x(+00))
\right))\\
&&=\square_2^-(x(--0)*_0 x(00+)) + \square_2^-(x(0-0)*_0 x(+0+)) + \square_2^-(x(0--)*_0 x(+00))
\\
&&=\square_2^-(x(0-0))
\eeas
\epf

The preceding propositions can be in fact generalized as follows :

\bth\label{diff} Let $x$ be an element of $\omega Cat(I^n,\C)^-$ with $n\geq
2$. Then in the reduced branching complex, one has \beas &&
\square_{n-1}^-(s_{n-1}x(0_n))=\sum_{1\leqslant 2i+1\leqslant n} \square_{n-1}^-
((\de_{2i+1}^-x)(0_{n-1}))\\ &&
\square_{n-1}^-(t_{n-1}x(0_n))=\sum_{1\leqslant 2i\leqslant n} \square_{n-1}^-
((\de_{2i}^-x)(0_{n-1}))\eeas \eth

\bpf For all $n$, we have seen that $\Phi_n^-$ induces the identity map on
the reduced branching complex. Therefore for all $x\in \omega Cat(I^n,\C)^-$,
$\Phi_{n-1}^- \de^- x=\de^- \Phi_n^- x$.
The latter equality can be translated into
$$\sum_{1\leqslant 2i+1\leqslant n} \Phi_{n-1}^- (\de_{2i+1}^- x) -
\sum_{1\leqslant 2i\leqslant n} \Phi_{n-1}^- (\de_{2i}^- x)=
\square_{n-1}^- s_{n-1} x(0_n) - \square_{n-1}^- t_{n-1} x(0_n).$$ If
the above equality was in $\Z\omega Cat(I^{n-1},\C)^-$, the proof
would be complete.  Unfortunately, we are working in the reduced
branching chain complex, and so there exists $t_1\in M_{n-1}^-$ and $t_2\in
M_n^-$ such that, in $\Z\omega Cat(I^{n-1},\C)^-$ 

\[\sum_{1\leqslant 2i+1\leqslant n} \Phi_{n-1}^- (\de_{2i+1}^- x) -\sum_{1\leqslant 2i\leqslant n} \Phi_{n-1}^- (\de_{2i}^- x)=\square_{n-1}^- s_{n-1} x(0_n) - \square_{n-1}^- t_{n-1}x(0_n)+t_1+\de^- t_2.\]


Set $t_2=\sum_{i\in I} \lambda_i T_i$
where $T_i$ are thin elements of $\omega Cat(I^n,\C)^-$. Each $T_i$
corresponds to a thin $(n-1)$-cube in the free cubical
$\omega$-category generated by the $n$-cube $x$ which will be denoted
in the same way (see the last paragraph of Section~\ref{id_phi}). One
can suppose that each $T_i(0_n)$ is $(n-1)$-dimensional. In the free
cubical $\omega$-category generated by $x$, either $T_i$ is in the
cubical $\omega$-category generated by the $\de_i^{(-)^i}x$ for
$1\leqslant i \leqslant n$ (let us denote this fact by $T_i \leqslant
s_{n-1}x(0_n)$), or $T_i$ is in the cubical $\omega$-category
generated by the $\de_i^{(-)^{i+1}}x$ for $1\leqslant i \leqslant n$
(let us denote this fact by $T_i \leq t_{n-1}x(0_n)$). Therefore one
has
$$t_2=\sum_{i\in I,T_i \leqslant s_{n-1}x(0_n)} \lambda_i T_i+ \sum_{i\in I,T_i \leqslant t_{n-1}x(0_n)}\lambda_i T_i.$$
and
\beas
&&\de^- t_2=\\
&&\sum_{\begin{array}{c}i\in I,T_i \leqslant s_{n-1}x(0_n)\\1\leqslant j\leqslant n, \de_j^- T_i \hbox{ thin}\end{array}} (-1)^{j+1}\lambda_i \de_j^- T_i + \sum_{\begin{array}{c}i\in I,T_i \leqslant s_{n-1}x(0_n)\\1\leqslant j\leqslant n, \de_j^- T_i \hbox{ non-thin}\end{array}} (-1)^{j+1}\lambda_i \de_j^- T_i
\\&&+
\sum_{\begin{array}{c}i\in I,T_i \leqslant t_{n-1}x(0_n)\\1\leqslant j\leqslant n, \de_j^- T_i \hbox{ thin}\end{array}} (-1)^{j+1}\lambda_i \de_j^- T_i + \sum_{\begin{array}{c}i\in I,T_i \leqslant t_{n-1}x(0_n)\\1\leqslant j\leqslant n, \de_j^- T_i \hbox{ non-thin}\end{array}} (-1)^{j+1}\lambda_i \de_j^- T_i
\eeas
Because of the freeness of $\Z\omega Cat(I^{n-1},\C)^-$, one gets
\beas
&& \sum_{1\leqslant 2i+1\leqslant n} \Phi_{n-1}^- (\de_{2i+1}^- x)=\square_{n-1}^- s_{n-1} x(0_n) + \sum_{\begin{array}{c}i\in I,T_i \leqslant s_{n-1}x(0_n)\\1\leqslant j\leqslant n, \de_j^- T_i \hbox{ non-thin}\end{array}} (-1)^{j+1}\lambda_i \de_j^- T_i\\
&&
\sum_{1\leqslant 2i\leqslant n} \Phi_{n-1}^- (\de_{2i}^- x)= \square_{n-1}^- t_{n-1} x(0_n) - \sum_{\begin{array}{c}i\in I,T_i \leqslant t_{n-1}x(0_n)\\1\leqslant j\leqslant n, \de_j^- T_i \hbox{ non-thin}\end{array}} (-1)^{j+1}\lambda_i \de_j^- T_i\\
&& -t_1 = \sum_{\begin{array}{c}i\in I,T_i \leqslant s_{n-1}x(0_n)\\1\leqslant j\leqslant n, \de_j^- T_i \hbox{ thin}\end{array}} (-1)^{j+1}\lambda_i \de_j^- T_i + \sum_{\begin{array}{c}i\in I,T_i \leqslant t_{n-1}x(0_n)\\1\leqslant j\leqslant n, \de_j^- T_i \hbox{ thin}\end{array}} (-1)^{j+1}\lambda_i \de_j^- T_i
\eeas

\epf

\section{Some consequences for the reduced branching homology}\label{invariance_result}

The following result generalizes the invariance result of \cite{Gau}
for the branching homology theory.

\bp\label{invariance} Let $f$ and $g$ be two non-contracting $\omega$-functors from
$\C$ to $\D$ satisfying the following conditions :
\begin{itemize}
\item for any $0$-morphism $x$, $f(x)=g(x)$
\item for any $n$-morphism $x$, $f(x)$ and $g(x)$ are
two homotopic morphisms (and so of the same dimension).
\end{itemize}
Then for any $n\geqslant 0$,
$\HR_n^\pm (f)=\HR_n^\pm (g)$.  \ep

\bpf Consider
the case of the reduced branching homology. Let $x\in \CR_n^-(\C)$.
If $$dim(f(x(0_n)))=dim(g(x(0_n))<n,$$ then $f(x)$
and $g(x)$ are two thin elements of $\omega Cat(I^n,\C)^-$. Therefore
$f(x)=g(x)$ in the reduced branching complex of $\D$. Now suppose that
$$dim (f(x(0_n)))=dim (g(x(0_n)))=n.$$ By hypothesis, there exists $z\in
\D_{n+1}$ such that $f(x(0_n))-g(x(0_n))=(s_n-t_n)(z)$. Therefore in the
reduced branching complex, one has
$f(x)-g(x)=\square_n^-((s_n-t_n)(z))=\de^- \square_{n+1}^-(z)$. So
$f(x)-g(x)$ is a boundary. \epf

We end up this section with another invariance result for the reduced
branching homology and with some results related to
Question~\ref{formelle}.

\bth\label{invariance2}
Let $\C$ and $\D$ be two $\omega$-categories. Let $f$ and $g$ be two
non $1$-contracting $\omega$-functors from $\C$ to $\D$ which coincide for
the $0$-morphisms and such that for
any $n\geqslant 1$, there exists a linear map $h_n$ from $\CF_n^-(\C)$
to $\CF_{n+1}^-(\D)$ such that for any $x\in \CF_n^-(\C)$,
$h_{n-1}(s_{n-1}-t_{n-1})+(s_n-t_n)h_n(x)=f(x)-g(x)$.  Then
$\HR_n^-(f)=\HR_n^-(g)$ for any $n\geqslant 0$.
\eth

\bpf Set $h_n^-x=\square_{n+1}^-h_n(x(0_n))$
for any $x\in \omega Cat(I^n,\C)^-$. It is clear that $h_n^-(M_n^-(\C))=\{0\}$
in $\CR_{n+1}^-(\D)$. Now suppose that $x=\de^-y$ for some $y\in M_{n+1}^-(\C)$.

We already mentioned that $I^n[-_n,+_n]$ is the free
$\omega$-category generated by a composable pasting scheme in the
proof of Corollary~\ref{prem_cit_perm}. It turns out that
$s_n(R(0_{n+1}))$ and $t_n(R(0_{n+1}))$ belong to $I^n[-_n,+_n]$
and it is possible thereby to use the explicit combinatorial
description of \cite{symetrique_cube}.

Set $I=\{1,2,\dots,n\}$ equipped with the total order $1<2<\dots<n$.
Let $C(I,k)$ (or $C(n,k)$) be the set of all subsets of $I$ of
cardinality $k$.  Let $\mathcal{P}$ an arbitrary subset of $C(I,k)$.
There is a lexicographical order on $C(I,k)$ usually defined as
follows : if $J=(j_1,\dots,j_k)$ with $j_1<\dots<j_k$ and
$J'=(j'_1,\dots,j'_k)$ with $j'_1<\dots<j'_k$, then $J\leqslant J'$
means that either $j_1<j'_1$, or $j_1=j'_1$ and $j_2<j'_2$, etc.  If
$K\in C(I,k+1)$, a $K$-packet is a set like $P(K)=\{J,J\in
C(I,k),J\subset K\}$. If $K=(i_1,\dots,i_{k+1})$ with $i_j<i_{j+1}$,
then $P(K)$ consists of the sets $K_{\widehat{a}}=K-\{i_a\}$ for
$a=1,\dots,k+1$. We have lexicographically
$$K_{\widehat{k+1}}<K_{\widehat{k}}<\dots<K_{\widehat{1}}.$$ A total
order $\sigma$ on $C(I,k)$ will be denoted by $\sigma=J_1J_2\dots J_N$
for $N=\left(\begin{array}{c}n\\k\end{array}\right)$, that is
$J_i\sigma J_j$ for $i<j$. A total order is called \textit{admissible}
by Manin and Schechtman if on each packet it induces either a
lexicographical order or the inverse lexicographical order. The set of
admissible orders of $C(I,k)$ is denoted by $A(I,k)$ (or
$A(n,k)$). Two total orders $\sigma$ and $\sigma'$ of $A(I,k)$ are
called \textit{elementary equivalent} if they differ by an interchange
of two neighbours which do not belong to a common packet. The quotient
of $A(I,k)$ by this equivalence relation is denoted by $B(I,k)$ (or
$B(n,k)$). Suppose that for some $K\in C(I,k+1)$, the members of the
packet $P(K)$ form a chain with respect to an admissible order
$\sigma$ of $A(I,k)$, i.e. any element of $C(I,k)$ lying between two
elements of $P(K)$ belongs to $P(K)$. Define $p_K(\sigma)$ the
admissible order in which this chain is reversed while all the rest
elements conserve their positions.  Then $p_K(\sigma)$ is still an
admissible order and $p_K$ passes to the quotient $B(I,k)$. The lemma on 
page 300 claims that
$A(n,n-1)=B(n,n-1)=\{K_{\widehat{n}}\dots K_{\widehat{1}},K_{\widehat{1}}\dots K_{\widehat{n}}\}$
where $K=(1,\dots,n)$. And the poset $B(n,n-2)$ is described by the
following picture :

\[
\xymatrix{ & {r_{min}} \ar@{->}[ld]_{p_{K_{\widehat{n}}}}\ar@{->}[rd]^{p_{K_{\widehat{1}}}}& \\
\ar@{->}[d]_{p_{K_{\widehat{n-1}}}}&& \ar@{->}[d]^{p_{K_{\widehat{2}}}}\\
{\vdots}\ar@{->}[d]_{p_{K_{\widehat{2}}}} & & {\vdots}\ar@{->}[d]^{p_{K_{\widehat{n-1}}}} \\
\ar@{->}[rd]_{p_{K_{\widehat{1}}}} && \ar@{->}[ld]^{p_{K_{\widehat{n}}}}\\
& {r_{max}} &\\}
\]

It turns out that in the picture $B(n,n-2)$, the vertices are exactly
the $(n-2)$-morphisms of $I^n[-_n,+_n]$ and the arrows are exactly the
$(n-1)$-morphisms of $I^n[-_n,+_n]$. This explicit description shows therefore that
$s_n(R(0_{n+1}))$ is equal to a composition $X_1 *_{n-1}\dots *_{n-1}
X_{n+1}$ where the only morphism of dimension $n$ contained in $X_j$
is $R(\delta_j^{(-)^j}(0_n))$.  And the same description shows that
$t_n(R(0_{n+1}))$ is equal to a composition $Y_{n+1} *_{n-1}\dots *_{n-1}
Y_{1}$ where the only morphism of dimension $n$ contained in $Y_j$
is $R(\delta_j^{(-)^{j+1}}(0_n))$. And one has
\beas
&&s_n(y(0_{n+1}))=y (s_n(0_{n+1}))=y(X_1)*_{n-1}\dots  *_{n-1}y(X_{n+1})\\
&& t_n(y(0_{n+1}))=y (t_n(0_{n+1}))=y(Y_{n+1}) *_{n-1} \dots *_{n-1} y(Y_1)
\eeas
Since $y$ is thin, $s_n(y(0_{n+1}))=t_n(y(0_{n+1}))$. Since $h_n$ is a
map from $\CF_n^{-}(\C)$ to $\CF_{n+1}^{-}(\D)$, then
$$\sum_{p=1}^{n+1}h_n(y(X_{p}))=\sum_{p=1}^{n+1}h_n(y(Y_{p}))$$
in $\CF_{n+1}^{-}(\D)$.

Since $I^{n+1}$ is the free $\omega$-category generated by the pasting
scheme $\underline{cub}^{n+1}$, then for any $p$ between $1$ and $n+1$, $X_p$
is a composition of $R(\delta_p^{(-)^p}(0_n))$ with other $R(k_1\dots
k_{n+1})$ of dimension strictly lower than $n$. Suppose that $p$ is
odd. There exists $X_p^{(1)}$ and $X_p^{(1)}{'}$ such that
$X_p=X_p^{(1)}*_{i_1}X_p^{(1)}{'}$ for some $0\leqslant i_p \leqslant
n-2$. If $i_p>0$, then only one of the $X_p^{(1)}$ or $X_p^{(1)}{'}$
is of dimension $n$ therefore $y(X_p)=y(X_p^{(1)})$ or
$y(X_p)=y(X_p^{(1)}{'})$. If $i_p=0$, then since $s_0 X_p= s_0
X_p^{(1)}=s_0 R(\delta_p^{(-)^p}(0_n))$ then in this case $ X_p^{(1)}$
is $n$-dimensional and $X{'}_p^{(1)}$ is of dimension strictly lower
than $n$. Therefore in this case $h_n (y(X_p))=h_n(y(X_p^{(1)}))$. By
repeating as many times as necessary the process, the number of cells
$R(k_1\dots k_n)$ included in $y(X_p)$ decreases. And we obtain
$$h_n(y(X_p))=h_n(y(\delta_p^{(-)^p}(0_n)))=h_n((\de_p^-y)(0_n)).$$
Now suppose that $p$ is even.  Since $R(-_n)=s_0 (X_p)\neq s_0
R(\delta_p^{(-)^p}(0_n))$, then necessarily at one step of the
process, we have $i_h=0$. Take the last $h$ such that $i_h=0$. Then
$h_n(y(X_p))=h_n(y(X_p^{(h)}))$ and $X_p^{(h)}=X_p^{(h+1)} *_0
X_p^{(h+1)}{'}$. Since $s_0 X_p^{(h)}= s_0 X_p^{(h+1)} \neq s_0
R(\delta_p^{(-)^p}(0_n))$, then for this $h$, $X_p^{(h+1)}$ is of
dimension strictly lower than $n$ and $X_p^{(h+1)}{'}$ is of dimension
$n$. Therefore $h_n(y(X_p))=0$.

In the same way, $h_n(y(Y_p))=0$ if $p$ is odd and
$$h_n(y(Y_p))=h_n(y(\delta_p^{(-)^{p+1}}(0_n)))=h_n((\de_p^-y)(0_n))$$
if $p$ is even. So 
$$h_n^-(x)=\square_{n+1}^-h_n(x(0_n))=\sum_{p=1}^{n+1}(-1)^{p+1}\square_{n+1}^-h_n((\de_p^-y)(0_n))=0$$
in $\CR_{n+1}^-(\D)$ by
Theorem~\ref{composition_glob_folding}. Therefore $h_n^-$ induces a
linear map from $\CR_n^-(\C)$ to $\CR_{n+1}^-(\D)$ still denoted by
$h_n^-$.
Take $x\in \omega Cat(I^n,\C)^-$. Then in $\CR_n^-(\D)$, one has
\beas
&&\de^-h_n^-(x)+ h_{n-1}^-\de^-(x)\\
&&=\de^-\square_{n+1}^- h_n(x(0_n))+h_{n-1}^-\de^-\Phi_{n}^-(x)\hfill\hbox{  since $\Phi_{n}^-$ is the identity map}\\
&&=\square_{n}^-(s_n-t_n)h_n(x(0_n))+h_{n-1}^-\de^-\square_{n}^- x(0_{n})\hbox{ by definition of $\Phi_n^-$}\\
&&=\square_{n}^-(s_n-t_n)h_n(x(0_n))+h_{n-1}^-\square_{n-1}^-(s_{n-1}x(0_n)-t_{n-1}x(0_n))\\
&&=\square_{n}^-(s_n-t_n)h_n(x(0_n))+\square_{n}^-h_{n-1}(s_{n-1}x(0_n)-t_{n-1}x(0_n))\hbox{ by definition of $h_n^-$}\\
&&=\square_{n}^-(f(x)(0_n)-g(x)(0_n))\hbox{ by hypothesis on $h_*$}\\
&&=\Phi_n^-(f(x)-g(x))\hbox{ by definition of $\Phi_n^-$}\\
&&=f(x)-g(x)\hbox{ since $\Phi_{n}^-$ is the identity map}
\eeas

\epf

The proof of Theorem~\ref{invariance2} provides another way of proving
Theorem~\ref{diff} and also establishes that Theorem~\ref{diff} is
still true for the formal branching homology.

\bp Let $p\geqslant 1$ and let $2_p$ be the $\omega$-category
generated by a $p$-morphism $A$. Then $\HF_n^-(2_p)=\HR_n^-(2_p)=0$ for
$n>0$ and $\HF_0^-(2_p)=\HR_0^-(2_p)=\Z$. \ep

\bpf The assertions concerning the formal branching homology are obvious.
Since the negative folding operator induces the identity on the
reduced branching complex, then $\CR_n^-(2_p)$ is equal to $0$ for $n>p$
and is generated by $\square_n^-(s_n A)$ and $\square_n^-(t_n A)$ for
$0\leqslant n\leqslant p$. The point is to prove that there is no
relations between $\square_n^-(s_n A)$ and $\square_n^-(t_n A)$ for
$1\leqslant n<p$, that is $\CR_n^-(2_p)=\Z \square_n^-(s_n A)\oplus \Z
\square_n^-(t_n A)=\CF_n^-(2_p)$. Suppose that there exists a linear combination
of thin $n$-cubes
$t_1$ and a linear combination of thin $(n+1)$-cubes $t_2$ such that
for some integers $\lambda$ and $\mu$,
$$\lambda \square_n^-(s_n A)+\mu \square_n^-(t_n A)=t_1+\de^- t_2$$ in
$C_n^-(2_p)$. Then $s_n t_2(0_{n+1})= t_n t_2(0_{n+1})$ and so
$\de^- t_2$ is necessarily a linear combination of thin $n$-cube
therefore $\lambda=\mu=0$.

Another possible proof of this proposition is to use
Theorem~\ref{invariance2} and to use the homotopy equivalence of
\cite{Gau} Proposition 8.5 between $2_p$ and $2_1$.
\epf

\bp Let $p\geqslant 1$ and let $G_p\AB$ be the $\omega$-category
generated by two non-homotopic $p$-morphisms $A$ and $B$. Then
$\HF_n^-(G_p\AB)=\HR_n^-(G_p\AB)=0$ for $0<n<p$ and
$\HF_0^-(G_p\AB)=\HR_0^-(G_p\AB)=\Z=\HF_p^-(G_p\AB)=\HR_p^-(G_p\AB)$. \ep

\bpf Analogous to the previous proof. \epf

\bp Let $n\geq 0$.  $\HF_0^-(I^n)=\Z$ and for $p>0$,
$\HF_p^-(I^n)=0$. \ep

\bpf We know that
$$\CF_p^-(I^n)=\bigoplus_{R(k_1\dots k_n) \hbox{ of dimension $p$}}
\Z\square_p^-(R(k_1\dots k_n)).$$ And the differential maps is
also completely known. In the formal branching complex, one
has $$\de^- \square_p^-(R(k_1\dots k_n))=\sum_{1\leqslant
j\leqslant p} (-1)^{j+1} \square_p^-(R(k_1\dots[-]_{n_j}\dots
k_n))$$ where $k_{n_1},\dots,k_{n_p}$ are the $0$'s appearing in
the word $k_1\dots k_n$ with $n_1<\dots<n_p$. It follows that
this chain complex can be splitted depending on the position and
the number of the $+$'s, and that these positions and numbers are  not 
modified by the
differential maps. If the number of the $+$ signs is $N$, we are
reduced to calculating the simplicial homology of the
$(n-N)$-simplex which is known to vanish in dimension strictly
greater than $0$. \epf

As for the calculation of $\HR_*^-(I^n)$, the point is to prove as
above for $2_p$ and $G_p\AB$ that there is no additional relations
between the $\square_p^-(R(k_1\dots k_n))$ in the reduced branching
complex. Unfortunately, for a thin $(n+1)$-cube $t_2$ of the 
branching nerve of $I^n$, $\de^-t_2$ is not necessarily a linear
combination of thin $n$-cube. For example if $a$ and $b$ are two
$1$-morphisms of $I^n$ such that $a*_0b$ exists, then let $t_2$ the
thin $2$-cube such that $\de_1^- t_2=\square_1(a*_0b)$, $\de_2^+
t_2=\square_1(t_0 b)$, $\de_2^- t_2=\square_1(a)$ and $\de_1^+
t_2=\square_1(b)$. Then $\de^- t_2=\square_1(a*_0b)-\square_1(a)$.

\section{Acknowledgments}

I would like to thank the London Mathematical Society for supporting
the workshop ``Geometric methods in Computer Science'' held in Bangor
(UK) on September 8-10th 1999 and the participants for their interest
in this work.

\end{document}